%% file: ms.tex
\pdfoutput=1

\documentclass[a4paper, 11pt, abstract=true, titlepage, oneside]{scrartcl}
\makeatletter

\input{resources/packages.tex}
\input{resources/configuration.tex}

\newboolean{bluechanges}
\setboolean{bluechanges}{false} 

\newcommand{\revnew}[1]{%
    \ifthenelse{\boolean{bluechanges}}{\textcolor{blue}{#1}}{#1}%
}
\begin{document}
\emergencystretch 3em
\input{./resources/acronyms}

\input{./resources/titlepage.tex}
\input{./contents/1.Introduction}
\input{./contents/2.ProblemStatement}
\input{./contents/3.Methodology}
\input{./contents/4.CaseStudy}
\input{./contents/5.Conclusion}


%
\singlespacing{
\begingroup
\scriptsize 
\bibliographystyle{model5-names}
\bibliography{references.bib}
\endgroup
} 
\newpage
\onehalfspacing
\begin{appendices}
	\normalsize
	\input{./contents/6.Appendix}

\end{appendices}
\end{document}

%% file: resources/packages.tex
\usepackage[T1]{fontenc}
\usepackage[utf8]{inputenc}
\usepackage[onehalfspacing]{setspace}
\usepackage[headsepline]{scrlayer-scrpage}
\usepackage{layout}
\usepackage{geometry}
\usepackage{adjustbox}
\usepackage{authblk}
\usepackage[titletoc,title]{appendix}
\usepackage[hyphens]{url}
\usepackage{color}

\usepackage{amsmath,amssymb,amsfonts,amsthm, mathtools, bm}
\usepackage{mathrsfs}
\usepackage{accents}
\usepackage{thmtools, thm-restate}
\usepackage{nicefrac}
\usepackage{array}

\usepackage[normal]{threeparttable}
\usepackage{threeparttablex}
\usepackage{tabularx}
\usepackage{longtable}
\usepackage{booktabs}
\usepackage{colortbl}
\usepackage{multirow}
\usepackage{makecell}

\usepackage{algorithm,algpseudocode}
\usepackage[bottom,hang,multiple]{footmisc}
\usepackage{algorithm}
\usepackage{algpseudocode}
\usepackage{tikz}
\usepackage{comment}
\usepackage{units}
\usepackage{accents}
\usetikzlibrary{backgrounds}
\pgfdeclarelayer{background}
\pgfdeclarelayer{foreground}
\pgfsetlayers{background,main,foreground} 
\usepackage{subfigure}
\usepackage{graphicx}
\usepackage{mathtools, amsmath}
\usepackage{amsthm}
\newcommand{\halmos}{\qed}
\usepackage{import}
\usepackage[acronym,nomain,nonumberlist, section=subsection]{glossaries}
\usepackage{tabularx}
\usepackage{booktabs}
\usepackage{pgfplots}
\usepackage{import}
\usepackage{setspace}
\usepackage[bottom]{footmisc}
\usepgfplotslibrary{groupplots,dateplot}
\usetikzlibrary{patterns,shapes.arrows}
\pgfplotsset{compat=newest}
\newtheorem{result}{Result}
\usepackage[acronym]{glossaries}

\usepackage{subcaption}
\usepackage{import}
\usepackage{soul}
\usepackage{graphicx}
\usepackage{xcolor} 

\usepackage{placeins}

\usepackage{tikz}
\usetikzlibrary{arrows.meta}
\usetikzlibrary{math}
\usetikzlibrary{shapes}
\usetikzlibrary{patterns,shapes.arrows}
\usepackage{pgfplots}
\pgfplotsset{compat=1.18} 
\usepgfplotslibrary{groupplots}
\usepgfplotslibrary{groupplots,dateplot}
\pgfplotsset{compat=newest}
\usepgfplotslibrary{dateplot}
\usetikzlibrary{backgrounds}
\pgfdeclarelayer{background}
\pgfdeclarelayer{foreground}
\pgfsetlayers{background,main,foreground} 
\usepackage{tikzscale}

\usepackage{enumerate}
\usepackage{multicol}
\usepackage[author=Patrick]{pdfcomment}
\usepackage{xparse}
\usepackage{twoopt}
\usepackage{etoolbox}

\usepackage{eurosym}
\usepackage{ellipsis}
\usepackage{natbib}
\usepackage{paralist}
\usepackage{ulem}
\usepackage{overpic} 

\usetikzlibrary{shapes.geometric, arrows.meta, positioning}
\tikzstyle{startstop} = [ellipse, draw=black, thick, fill=green!20,
    minimum height=1cm, minimum width=2.3cm, text width=2.8cm, align=center]
\tikzstyle{process} = [rectangle, rounded corners, minimum height=1.1cm,
    text centered, draw=black, thick, fill=blue!10, text width=3.2cm]
\tikzstyle{subprocess} = [rectangle, rounded corners, text centered,
    draw=black, thick, fill=teal!15, text width=3.3cm, minimum height=1cm]
\tikzstyle{fixstep} = [rectangle, text centered, draw=black, thick,
    fill=gray!10, text width=2.8cm, minimum height=1cm]
\tikzstyle{decision} = [diamond, draw=black, thick, fill=orange!20,
    text centered, inner sep=1pt, aspect=2.3, text width=3.1cm]
\tikzstyle{arrow} = [thick,->,>=stealth]

\makeatletter
\renewenvironment{enumerate}
  {\ifnum \@enumdepth >3 \@toodeep\else
   \advance\@enumdepth\@ne
   \edef\@enumctr{enum\romannumeral\the\@enumdepth}%
   \list{\csname label\@enumctr\endcsname}%
     {\usecounter{\@enumctr}%
      \setlength{\topsep}{0pt}%
      \setlength{\itemsep}{0pt}%
      \setlength{\parsep}{0pt}%
      \setlength{\partopsep}{0pt}%
      \setlength{\leftmargin}{1.5em}}%
   \fi}
  {\endlist}
\makeatother

%% file: resources/configuration.tex
\AtBeginDocument{
	\newgeometry{
		textheight=9in,
		textwidth=6.5in,
		top=1in,
		headheight=16.49675pt,
		headsep=25pt,
		footskip=30pt
	}
}
\newsavebox{\abstractbox}
\renewenvironment{abstract}
{\begin{lrbox}{0}\begin{minipage}{\textwidth}
			\begin{center}\normalfont\sectfont\abstractname\end{center}\quotation}
		{\endquotation\end{minipage}\end{lrbox}%
	\global\setbox\abstractbox=\box0 }

\makeatletter
\expandafter\patchcmd\csname\string\maketitle\endcsname
{\vskip\z@\@plus3fill}
{\vskip\z@\@plus2fill\box\abstractbox\vskip\z@\@plus1fill}
{}{}
\makeatother

\allowdisplaybreaks

\addtokomafont{captionlabel}{\footnotesize\bfseries} 
\addtokomafont{caption}{\footnotesize\bfseries} 

\bibpunct[, ]{(}{)}{,}{a}{}{,}%
\newtheoremstyle{upright} 
  {3pt}                    
  {3pt}                    
  {\upshape}               
  {}                       
  {\bfseries}              
  {.}                      
  {.5em}                   
  {}                       

\theoremstyle{upright}

\newtheorem{theorem}{Theorem}

\let\proof\relax 

\let\theoremstyle\relax

\counterwithin*{equation}{section}

\newenvironment{myfont}{\fontfamily{ccr}\selectfont}{\par}
\DeclareTextFontCommand{\textmyfont}{\myfont}

\AtBeginDocument{%
	\abovedisplayskip=7.5pt plus 0pt minus 0pt
	\abovedisplayshortskip=7.5pt plus 0pt minus 0pt
	\belowdisplayskip=7.5pt plus 0pt minus 0pt
	\belowdisplayshortskip=7.5pt plus 0pt minus 0pt
}

\newcolumntype{L}[1]{>{\raggedright\let\newline\\\arraybackslash\hspace{0pt}}p{#1}}
\newcolumntype{C}[1]{>{\centering\let\newline\\\arraybackslash\hspace{0pt}}p{#1}}
\newcolumntype{R}[1]{>{\raggedleft\let\newline\\\arraybackslash\hspace{0pt}}p{#1}}
\renewcommand{\emph}[1]{\textit{#1}}

\algtext*{EndIf} 
\algtext*{EndWhile} 
\algtext*{EndFor} 

\definecolor{legendfill}{rgb}{0.9,0.9,0.9}  

%% file: resources/acronyms.tex
\newcommand{\SolutionCost}[1]{C\left({#1}\right)}
\newcommand{\DetourDelay}[1]{D^{'}\left(#1\right)}
\newcommand{\CongestionDelay}[1]{D^{''}\left(#1\right)}
\newcommand{\TotalDelayRemoved}[0]{D^{\mathbin{-}}}
\newcommand{\TotalDelayExogenous}[0]{D^{\text{B}}}
\newcommand{\TotalDelayFleet}[0]{D^{\text{AMoD}}}
\newcommand{\TotalDelay}[1]{D{\left(#1\right)}}
\newcommand{\AvgTotalDelay}[0]{\bar{D}}
\newcommand{\AvgTotalDelayExogenous}[0]{\bar{D}^{\text{B}}}
\newcommand{\AvgTotalDelayFleet}[0]{\bar{D}^{\text{AMoD}}}
\newcommand{\F}[1]{F\left(#1\right)}
\newcommand{\Infeasibility}[1]{I({#1})}
\newcommand{\TripInfeasibility}[1]{I^{#1}}
\newcommand{\Job}{J}
\newcommand{\NumPieces}{K}
\newcommand{\Loss}[3]{L\left(#1,#2,#3\right)}
\newcommand{\BigM}[1]{M_{#1}}
\newcommand{\Queue}[0]{\text{Q}}
\newcommand{\SortedTripsOnArc}[1]{R_{#1}}
\newcommand{\SortedInactiveTripsOnArc}[1]{R^{\text{I}}_{#1}}
\newcommand{\SortedActiveTripsOnArc}[1]{R^{\text{A}}_{#1}}
\newcommand{\Sim}[2]{\operatorname{Sim}(#1, #2)}
\newcommand{\Schedule}{S}
\newcommand{\TimeHorizon}{T}

\newcommand{\TotalTravelTime}[1]{Z^{#1}}
\newacronym{amod}{AMoD}{Autonomous mobility-on-demand}
\newcommand{\Arc}{a}
\newcommand{\FirstArc}[1]{\ubar{a}_{#1}}
\newcommand{\LastArc}[1]{\bar{a}_{#1}}
\newacronym{co}{CO}{combinatorial optimization}
\newacronym{coml}{CO-enriched ML}{combinatorial optimization- enriched machine learning}
\newacronym{dta}{DTA}{dynamic traffic assignment}
\newacronym{dso}{DSO}{dynamic system optimum}
\newacronym{due}{DUE}{dynamic user equilibrium}
\newacronym{dfo}{DFO}{dynamic fleet optimum}
\newcommand{\TotalDelayTrip}[1]{d^{#1}}
\newcommand{\DelayOnArc}[2]{d_{#1}^{#2}}
\newcommand{\DelayFunction}[0]{d(\Flow{\Arc}{\Trip})}
\newcommand{\TripEarliestDeparture}[1]{e^{#1}}
\newacronym{fo}{FO}{fleet optimum}
\newcommand{\Flow}[2]{f_{#1}^{#2}}
\newcommand{\InactiveFlow}[2]{\ubar{f}_{#1}^{#2}}
\newcommand{\ActiveFlow}[2]{\bar{f}_{#1}^{#2}}
\newcommand{\ThFlow}[2]{\hat{f}_{#1}^{#2}}
\newacronym{fifo}{FIFO}{first-in first-out}
\newacronym{ga}{GA}{genetic algorithm}
\newcommand{\NeuralNetwork}{g_w}
\newacronym{hc}{HC}{high congestion}
\newcommand{\PostProcessing}{h}
\newacronym{kspwlo}{kSPwLO}{$k$-shortest paths with limited overlap}
\newacronym{lns}{LNS}{large neighborhood search}
\newacronym{lc}{LC}{low congestion}
\newcommand{\TripLatestArrival}[1]{\ell^{#1}}
\newcommand{\TripTimeWindow}[1]{\TripLatestArrival{#1} - \TripEarliestDeparture{#1}}
\newacronym{ml}{ML}{machine learning}
\newacronym{mfd}{MFD}{macroscopic fundamental diagram}
\newacronym{milp}{MILP}{mixed integer linear program}
\newcommand{\Move}[1]{\left(\TripPath{#1},\StartTime{#1}\right)}
\newcommand{\PopulationSize}[0]{m}
\newcommand{\SolutionPaths}[0]{p}
\newcommand{\Route}[2]{p^{#1}_{#2}}
\newcommand{\TripPath}[1]{p^{#1}}
\newcommand{\FirstArcRoute}[1]{\ubar{p}^{#1}_{i}}
\newcommand{\LastArcRoute}[1]{\bar{p}^{\hspace{.1em}#1}_{i}}
\newcommand{\NextRoute}[0]{p^{\Trip}_{+}}

\newacronym{rduo}{RDUO}{reactive dynamic user optimum}
\newcommand{\Trip}{r}
\newcommand{\SecondTrip}{r'}
\newacronym{srdtc}{SRDTC}{simultaneous route and departure time choice}
\newacronym{so}{SO}{system optimum}
\newcommand{\SolutionTimes}{s}
\newacronym{sl}{SL}{structured learning}
\newacronym{spt}{SPT}{shortest processing time first}
\newacronym{sgd}{SGD}{stochastic gradient descent}
\newcommand{\Departure}[2]{\ubar{s}_{#1}^{#2}}
\newcommand{\Arrival}[2]{\bar{s}_{#1}^{\hspace{1pt}#2}}
\newcommand{\StartTime}[1 ]{s^{#1}}
\newcommand{\LatestPathDeparture}[0]{\bar{s}\left(\TripPath{\Trip}\right)}
\newcommand{\LatestPathDepartureTested}[0]{s\left(\TripPath{\Trip}\right)}
\newcommand{\DecisionTime}{t}
\newacronym{ue}{UE}{user equilibrium}
\newcommand{\PredictedStatParams}{\hat{w}}
\newcommand{\StatParams}{w}
\newcommand{\TravelTimeVar}[2]{x_{#1}^{#2}}
\newcommand{\PathSelected}[2]{y_{#1}^{#2}}
\newcommand{\Arcs}{\mathcal{A}}
\newcommand{\ArcSetTrip}{\mathcal{A}_\Trip}
\newcommand{\TrainingSet}{\mathcal{D}}
\newcommand{\DiGraph}{\mathcal{G}}
\newcommand{\Population}{\mathcal{K}}
\newcommand{\Moves}{\mathcal{M}}
\newcommand{\Normal}{\mathcal{N}}
\newcommand{\Offspring}{\mathcal{O}}
\newcommand{\SetPaths}[1]{\mathcal{P}^{#1}}
\newcommand{\SetTrips}{\mathcal{R}}
\newcommand{\SetTripsArc}[1]{\mathcal{R}^{#1}}
\newcommand{\SetTripsTime}[1]{\mathcal{R}^{#1}}
\newcommand{\TripsInterval}[0]{\SetTrips^{\text{I}}}
\newcommand{\SetTripsDestroyed}{\mathcal{R}^{\text{D}}}
\newcommand{\SetOfArcArrivals}[1]{\mathcal{T}_{#1}}
\newcommand{\Uniform}[2]{\mathcal{U}(#1,#2)}
\newcommand{\DiscreteUniform}[2]{\bar{\mathcal{U}}(#1,#2)}
\newcommand{\Nodes}{\mathcal{V}}
\newcommand{\StatParamsSpace}{\mathcal{W}}
\newcommand{\SetInstances}{\mathcal{X}}
\newcommand{\SetSystemStates}{\mathcal{X}}
\newcommand{\EasySolutionSpace}[1]{\mathcal{Y}\left(#1\right)}
\newcommand{\ImprovementTracker}{\Delta}
\newcommand{\Perturbation}{Z}
\newcommand{\WeightsSpace}[1]{\Theta(#1)}
\newcommand{\SolutionSpace}[0]{\Pi}
\newcommand{\Gradient}[2]{\nabla_{#1} #2}
\newcommand{\LocalSearchInterval}[0]{\alpha}
\newcommand{\DestroyPercentage}[0]{\beta}
\newcommand{\MinimumImprovement}[0]{\gamma}
\newcommand{\SmallConstant}{\varepsilon}
\newcommand{\SmallProbability}{\epsilon}
\newcommand{\ProcessingTime}{\zeta}
\newcommand{\NumIterations}{\eta}
\newcommand{\SimilarityThs}{\theta}
\newcommand{\Individual}{\kappa}
\newcommand{\BestIndividual}{\kappa^{\star}}
\newcommand{\Penalization}{\lambda}
\newcommand{\PWLSlope}[2]{\mu_{#1}^{#2}}
\newcommand{\CrossoverPoint}{\xi}
\newcommand{\Solution}{\pi}
\newcommand{\OptimalSolution}{\pi^{\ast}}
\newcommand{\BackwardShift}[1]{\rho^{\text{#1}}}
\newcommand{\Staggering}[1]{\sigma^{#1}}
\newcommand{\TripMaxStaggering}[1]{\bar{\sigma}^{#1}}
\newcommand{\SigmaUp}[0]{\TripMaxStaggering{}\hspace{-.2em}\uparrow}
\newcommand{\SigmaDown}[0]{\TripMaxStaggering{}\hspace{-.2em}\downarrow}
\newcommand{\NominalTravelTime}[1]{\tau_{#1}}
\newcommand{\FreeFlowShortestRoute}[1]{\ubar{\tau}^{\Trip}}
\newcommand{\FreeFlowRoute}[1]{\ubar{\tau}^{\Trip}(#1)}
\newcommand{\TravelTimeRoute}[1]{{\tau}^{\Trip}(#1)}
\newcommand{\TravelTimeFunction}[2]{\tau_{#1}^{#2}(\Flow{#1}{#2})}
\newcommand{\TravelTime}[2]{\tau_{#1}^{#2}}
\newcommand{\Fitness}[1]{\phi(#1)}
\newcommand{\FractionControlled}{\varphi}
\newcommand{\EX}[1]{\mathbb{E}\left[#1\right]}
\newcommand{\Reals}{\mathbb{R}}

\newcommand{\ubar}[1]{\underaccent{\bar}{#1}}
\newcommand{\mathdefault}[1]{#1}

\newcommand{\SolutionStructureHeight}[0]{4.8cm}
\newcommand{\SolutionStructureWidth}[0]{0.2\columnwidth}
\newcommand{\SolutionStructureSpaceBetween}[0]{-1 cm}
\newcommand{\SolutionStructureXYLabel}{-0.8}
\newcommand{\SolutionStructureXYRouteLabel}{-0.25}
\newcommand{\DelayReductionWidth}[0]{.5\columnwidth}
\newcommand{\DelayReductionHeight}[0]{4cm}

\newcommand{\ViolinHeight}[0]{5.5cm}
\newcommand{\ViolinWidth}[0]{4cm}

\newcommand{\HeightFigureOne}[0]{5cm}
\newcommand{\RouteFigureWidth}[0]{\columnwidth}
\newcommand{\RouteFigureHeight}[0]{4.1cm}
\newcommand{\Width}[0]{7cm}
\newcommand{\HeightFlowControl}[0]{.3\linewidth}
\newcommand{\WidthFlowControl}[0]{.35\linewidth}

\newcommand{\InstancesHistogramsWidth}[0]{4cm}
\newcommand{\InstancesHistogramsHeight}[0]{4cm}

\newcommand{\PathAnalysisBoxplotsWidth}[0]{.9\linewidth}
\newcommand{\PathAnalysisBoxplotsHeight}[0]{4.6cm}

\newcommand{\MathComparisonWidth}[0]{\columnwidth}
\newcommand{\MathComparisonHeight}[0]{3cm}

\newcommand{\MathAppendixWidth}[0]{\columnwidth}
\newcommand{\MathAppendixHeight}[0]{3cm}

\newcommand{\HeatmapHorSep}{0.2cm}
\newcommand{\HeatmapVerSep}{.8cm}
\newcommand{\HeatmapWidth}[0]{.95\linewidth}
\newcommand{\HeatmapHeight}[0]{.95\linewidth}
\newcommand{\heatmappath}[1]{figures/appendix/parameter_search/heatmap-grid-relative-improvement_dp_#1_tex/}
\newcommand{\XCbar}{1.1}
\newcommand{\YCbar}{-0.83}
\newcommand{\HeightCbar}{2.1\columnwidth}

%% file: resources/titlepage.tex

\title{\large Integrated Balanced and Staggered Routing in Autonomous
Mobility-on-Demand Systems}
\author[1]{\normalsize Antonio Coppola}
\author[2]{\normalsize Gerhard Hiermann}
\author[3]{\normalsize Dario Paccagnan}
\author[4]{\normalsize Michel Gendreau}
\author[1,5]{\\ \normalsize Maximilian Schiffer}

\affil[1]{\small School of Management, Technical University of Munich, Germany

\texttt{antonio.coppola@tum.de}}

\affil[2]{\small Institute of Production and Logistics Management, Johannes Kepler University Linz, Austria

\texttt{gerhard.hiermann@jku.at}}

\affil[3]{\small Department of Computing, Imperial College London, U.K.

\texttt{d.paccagnan@imperial.ac.uk}}

\affil[4]{\small Department of Mathematical and Industrial Engineering, Polytechnique Montréal, Canada

\texttt{michel.gendreau@polymtl.ca}}

\affil[5]{\small Munich Data Science Institute, Technical University of Munich, Germany

\texttt{schiffer@tum.de}}

\date{}

\lehead{\pagemark}
\rohead{\pagemark}

\begin{abstract}
\begin{singlespace}
{\small\noindent \input{./contents/0.Abstract.tex}
\smallskip}

\end{singlespace}
\end{abstract}

\maketitle

%% file: contents/0.Abstract.tex
Autonomous mobility-on-demand (AMoD) systems — centrally coordinated fleets of self-driving vehicles — offer a promising alternative to traditional ride-hailing by improving traffic flow and reducing operating costs. Centralized control in AMoD systems enables two complementary routing strategies: balanced routing, which distributes traffic across alternative routes to ease congestion, and staggered routing, which delays departures to smooth peak demand over time. In this work, we introduce a unified framework that jointly optimizes both route choices and departure times to minimize system travel times. We formulate the problem as an optimization model and show that our congestion model yields an unbiased estimate of travel times derived from a discretized version of Vickrey’s bottleneck model. To solve large-scale instances, we develop a custom metaheuristic based on a large neighborhood search framework. To assess our method, we conduct a case study on the Manhattan street network using real-world taxi data. In a setting with exclusively centrally controlled AMoD vehicles, our approach reduces total traffic delay by up to 25\% and mitigates network congestion by up to 35\% compared to selfish routing. We also consider mixed-traffic settings with both AMoD and conventional vehicles, comparing a welfare-oriented operator that minimizes total system travel time with a profit-oriented one that optimizes only the fleet’s travel time. Independent of the operator's objective, the analysis reveals a win–win outcome: across all control levels, both autonomous and non-autonomous traffic benefit from the implementation of balancing and staggering strategies.

%% file: contents/1.Introduction.tex
\section{Introduction}\label{sec:introduction}
The ride-hailing market has recently experienced substantial growth, fueled by evolving consumer preferences and technological innovations. Projections estimate that by 2028, ride-booking platforms will serve nearly 100 million users, driving the U.S. market to a value of approximately \$55 billion \citep{Statista2024}.

While ride-hailing services offer convenient transportation options, they have been criticized for exacerbating urban congestion \citep{DiaoKongEtAl2021}, primarily due to uncoordinated operations \citep{KondorBojicEtAl2022} and drivers favoring personal route choices \citep{ZanardiSessaEtAl2023}, which result in suboptimal travel patterns and overloaded networks. This inefficiency carries a dual cost: from a societal standpoint, it worsens congestion and its externalities, including prolonged delays, air and noise pollution, and broader public health impacts; from the operator's standpoint, it drives up operational costs through increased travel time, fuel use, and vehicle underutilization.

\gls{amod} systems, centrally controlled fleets of self-driving vehicles offering door-to-door mobility solutions \citep{Pavone2015}, promise to overcome these limitations. Their central governance enables the coordination of vehicle flow in response to travel demand and traffic conditions, mitigating the inefficiencies associated with self-interested, uncoordinated fleets of human-driven vehicles. However, without effective routing policies, \gls{amod} systems risk replicating the issues observed in traditional ride-hailing services, potentially further straining the street network \citep{OhSeshadriEtAl2020}. To address congestion and enhance operational efficiency, \gls{amod} routing policies can employ two main coordination strategies: balanced and staggered routing. Balanced routing distributes \gls{amod} flow across alternative road segments to alleviate congestion on popular routes and smooth overall traffic flow. Staggered routing complements this approach by strategically delaying vehicle departures to spread demand more evenly over time, alleviating temporal imbalances in network utilization.

Although complementary, these operational strategies lack a unified framework for integrated implementation in \gls{amod} systems. On one hand, most existing studies that address \gls{amod} routing under congestion rely on static models \citep[cf.][]{ZardiniLanzettiEtAl2022}, which are well-suited for long-term planning but fall short in supporting \gls{amod} operations for two key reasons: first, they cannot capture the real-time nature of fleet management, which must continuously adapt to evolving travel demand; second, they cannot accommodate departure time assignments, which require explicit temporal modeling. On the other hand, previous work on optimal staggering of \gls{amod} departures lacks a routing control layer that leverages the additional travel time reductions enabled by staggering, thus treating route and departure time assignments in \emph{isolation}.

In contrast, an \emph{integrated} approach combines these strategies by moving the departure time decision to an earlier control stage. This shift influences route selection and can potentially outperform balanced and staggered routing methods when applied in isolation. However, this integrated optimization introduces additional combinatorial complexity and challenges the tractability of existing solution approaches. Against this background, we explore the untapped potential of integrating balanced and staggered routing strategies to coordinate autonomous fleets more efficiently under high travel demand. It introduces scalable algorithmic solutions to jointly optimize route assignment and departure timing, addressing the computational and operational complexity of large-scale, dynamic fleet coordination while accounting for congestion.

\subsection{Contribution}\label{sec:contribution}
This paper makes the following contributions. 
First, we formalize the problem of minimizing travel time in a network operated by a \gls{amod} system using an integrated approach that combines balanced and staggered routing. 
Second, we formally show that the congestion model employed serves as an unbiased estimator of travel times under a discretized version of Vickrey's bottleneck model — a well-established framework in dynamic traffic assignment. 
{
Third, given that the staggered routing problem—optimizing \gls{amod} departure times along fixed routes—is already NP-hard~\citep{CoppolaHiermannEtAl2024}, we develop a scalable metaheuristic based on a \gls{lns} framework to solve the more general problem of jointly optimizing routes and departure times.
We relax the assumption of piecewise-linear delay constraints and enable our algorithm to handle any convex, non-decreasing delay function.}
Fourth, we go beyond prior implicit representations of baseload traffic in the staggered routing literature by explicitly modeling a realistic number of trips, approximating real-world congestion conditions. Existing approaches typically rely on a \gls{milp} formulation tightened through trip time windows; however, the large number of resulting constraints limits the number of trips that can be explicitly tracked. To circumvent this, prior work approximates baseload traffic by imposing arc capacity limits instead of modeling individual trips — a limitation we overcome by no longer utilizing a \gls{milp} formulation.
Finally, we conduct a numerical study using real-world data from Manhattan, New York City, to examine the trade-offs between adjusting vehicle routes and departure times to minimize travel times and study network congestion. Our case study initially assumes the \gls{amod} operator has full traffic control, showing that our algorithmic solutions reduce fleet delays by up to 25\% and congestion delays by up to 35\% compared to a selfish baseline.
{We then transition to a mixed-traffic setting, where \gls{amod} vehicles serve only a fraction of the demand. This allows us to examine how varying the level of control affects overall network performance. In both welfare- and profit-oriented settings, centralized balanced and staggered routing consistently produce a win--win outcome, reducing travel times for both \gls{amod} and non-\gls{amod} traffic across all scenarios. Notably, controlling just 10\% of the system’s vehicles already achieves 25\% of the maximum delay reduction, while 50\% control yields 75\%, with further gains becoming more marginal.}

\section{Related literature}\label{sec:related_literature}
Our work builds upon previous research aimed at mitigating congestion through strategic routing and timing of traffic. 
First, adopting a theoretical perspective, we highlight the intersections of this work with the field of \gls{dta}. Second, from an applicative perspective, we focus on control and operational strategies for \gls{amod} systems.

\subparagraph{Overview of \gls{dta} models} The primary objective of traffic assignment is to model vehicular movement patterns based on predefined behavioral rules. A \gls{dta} model typically consists of a route and/or departure choice model combined with a traffic congestion model, so we review both components in sequence. Based on the available route and departure time choices, \gls{dta} models fall into one of three main categories: (1) pure departure time choice models \citep[see, e.g.,][]{Vickrey1969, HendricksonKocur1981, LindseyVandenBergEtAl2012}, (2) pure route choice models \citep[see, e.g.,][]{MounceCarey2011, CareyWatling2012, LongHuangEtAl2013, Levin2017}, and (3) \gls{srdtc} models \citep[see, e.g.,][]{JauffredBernstein1996, RanHallEtAl1996, ZiliaskopoulosRao1999, HuangLam2002, HanFrieszEtAl2013, JalotaPaccagnanEtAl2023}.

Focusing on \gls{srdtc} models, which are closest in nature to our work, route and departure time choices typically follow either the dynamic user optimal principle—where travelers select routes and departure times to minimize their own travel costs \citep{Nie2010, FrieszKimEtAl2011, UkkusuriHanEtAl2012}—or the dynamic system optimal principle, which seeks to minimize total system-wide travel cost \citep{LongWangEtAl2018}. Efforts to bridge the gap between user equilibrium and system optimality in \gls{srdtc} models include tolling schemes \citep{ArnottPalmaEtAl1990, YangMeng1998, LiuNie2011}, demand management and route guidance systems \citep{MakridisMenelaouEtAl2024, MenelaouTimotheouEtAl2021}, and signal control \citep{HuMahmassani1997}. Despite these efforts, research on achieving dynamic system optimality or dynamic fleet optimality \citep{BattifaranoQian2023} through centralized control of departures and routes remains limited.

A complete \gls{dta} model combines travelers' behavioral choices with a congestion model that captures how traffic evolves over time along routes connecting origin-destination pairs. Commonly adopted congestion models include extensions of Vickrey’s bottleneck model, which represent traffic formation and dissipation as a vertical queue at the bottleneck entrance, offering analytically tractable formulations \citep{LiHuangEtAl2020}, and cell transmission models, which can capture complex dynamics such as spillback effects \citep{AdacherTiriolo2018}. However, applying these models to large-scale networks remains computationally demanding \citep{ZiliaskopoulosWallerEtAl2004}. To address this, researchers have proposed macroscopic fundamental diagram-based models as aggregated tools for analyzing large-scale traffic behavior \citep{AghamohammadiLaval2020}. In this work, we introduce a discretized congestion model inspired by Vickrey’s bottleneck framework, enabling explicit vehicle-level routing and timing while preserving computational tractability.

\subparagraph{\gls{amod} control \& operations}
Most mesoscopic and macroscopic studies examining the impact of congestion in \gls{amod} systems focus on strategic control settings, optimizing the routes of customer-carrying and rebalancing trips under steady-state conditions \citep[see, e.g.,][]{RossiZhangEtAl2018, RossiIglesiasEtAl2020, EstandiaSchifferEtAl2021, SalazarLanzettiEtAl2019, BahramiRoorda2020, BangMalikopoulos2022}. While these approaches are valuable for long-term evaluation, they overlook the dynamic nature of \gls{amod} operations, limiting their applicability to our setting.

\noindent In contrast, research that considers the dynamic aspects of \gls{amod} operations, including studies on vehicle dispatching \citep[see, e.g.,][]{ZhangHuEtAl2017, XuLiEtAl2018, LiQinEtAl2019, TangQinEtAl2019, ZhouJinEtAl2019, LiangWenEtAl2021, SadeghiEshkevariTangEtAl2022, EndersHarrisonEtAl2023}, joint dispatching and rebalancing \citep[see, e.g.,][]{TsaoIglesiasEtAl2018, JungelParmentierEtAl2023}, request pooling \citep[see, e.g.,][]{AlonsoMoraWallarEtAl2017, TsaoMilojevicEtAl2019}, and fleet rebalancing \citep[see, e.g.,][]{LiuSamaranayake2022, IglesiasRossiEtAl2018b, JiaoTangEtAl2021, GammelliYangEtAl2021, SkordilisHouEtAl2021} rarely explicitly model congestion, as many rely on precise forecasts or deterministic assumptions for travel times.

In this work, we develop a scalable framework that explicitly models congestion and supports balanced and staggered routing for \gls{amod} systems. Our model enables dynamic assignment of routes and departure times at the vehicle level, using a congestion model inspired by Vickrey’s bottleneck formulation to provide a tractable representation of traffic dynamics.

\subsection{Roadmap}\label{sec:roadmap}
The remainder of this paper is organized as follows. Section~\ref{sec:problem_setting} introduces the problem setting for the integrated balanced and staggered routing problem and shows the relation between our congestion model and Vickrey's bottleneck model. Section~\ref{sec:methodology} describes our metaheuristic algorithm developed to solve real-world instances. In Section~\ref{sec:case_study}, we outline our experimental design, and in Section~\ref{sec:results}, we present the results from our numerical study. Finally, Section~\ref{sec:conclusion} concludes the paper with remarks and future research directions.

%% file: contents/2.ProblemStatement.tex
\section{Problem setting}\label{sec:problem_setting}

We consider an \gls{amod} operator tasked with coordinating a fleet of autonomous vehicles to serve a given set of trips within a fixed planning horizon. We focus on an offline setting in which the operator has full knowledge of all trip requests in advance. Each trip utilizes a single vehicle and can represent either on-duty activities (e.g., customer delivery) or off-duty activities (e.g., rebalancing). Our objective is to minimize total travel time by jointly optimizing route selection and trip departure times. This involves balancing trips across alternative (potentially longer in distance) routes and staggering departures beyond their originally requested times to prevent congestion and improve fleet efficiency.

We model the road network within which the fleet operates as a directed graph ${\DiGraph=(\Nodes, \Arcs)}$, where the nodes~$\Nodes$ represent trip origins, destinations, and road intersections, and the arcs~$\Arcs$ correspond to road segments connecting these nodes. We consider a set of trips $\SetTrips$, each traveling between two distinct network nodes. Let a quadruple $(\SetPaths{\Trip}, \TripEarliestDeparture{\Trip}, \TripLatestArrival{\Trip}, \TripMaxStaggering{\Trip})$ define a trip $\Trip \in \SetTrips$, with $\SetPaths{\Trip} = \{\Route{\Trip}{1}...\Route{\Trip}{k}\}$ denoting the alternative $k$ routes -- sequences of network arcs -- that the vehicle can travel, $\TripEarliestDeparture{\Trip}$ being the earliest departure from the trip's origin, $\TripLatestArrival{\Trip}$ being the latest feasible arrival at the destination, and $\TripMaxStaggering{\Trip}$ being the maximum staggering time applicable to a trip's departure. 
A trip $\Trip$ traversing arc $\Arc$ experiences a travel time $\TravelTimeFunction{\Arc}{\Trip}$, which consists of two components: the free-flow travel time $\NominalTravelTime{\Arc}$ and a congestion-induced delay $\DelayFunction$. The delay depends on the number of trips $\Flow{\Arc}{\Trip}$ present on arc $\Arc$ at the time $\Trip$ begins its traversal:
\begin{align}\label{eq:travel_time}
\TravelTimeFunction{\Arc}{\Trip} = \NominalTravelTime{\Arc} + \DelayFunction,
\end{align}
where $\DelayFunction$ is a convex, non-decreasing function of flow. 
A tuple $\Solution = (\SolutionPaths, \SolutionTimes)$ defines a solution to this problem, where $\SolutionPaths$ and $\SolutionTimes$ denote the route and departure time assigned to each trip $\Trip$, respectively. As we assume that vehicles do not idle once a trip has started, this information is sufficient to determine the travel time along each selected route and, consequently, the arrival time for every trip.
We require a solution $\Solution$ to satisfy the following constraints, which we collect in the feasible set $\SolutionSpace$:
\begin{enumerate}
\item[i)] The route $\TripPath{\Trip}$ assigned to each trip $\Trip$ must belong to $\SetPaths{\Trip}$.
\item[ii)] The start time $\StartTime{\Trip}$ of each trip $\Trip$ must occur between $\TripEarliestDeparture{\Trip}$ and $\TripEarliestDeparture{\Trip} + \TripMaxStaggering{\Trip}$.
\item[iii)] The combination of start times and routes in $\Solution$ must ensure that the arrival at the destination of each trip $\Trip$ occurs before $\TripLatestArrival{\Trip}$.
\end{enumerate}

In this setting, we seek for a solution $\OptimalSolution \in \SolutionSpace$ that minimizes the total fleet travel time, defined as the sum of travel times incurred by each trip across the arcs of its assigned route:
\begin{align}\label{eq:problem}
\min_{\Solution \in \SolutionSpace} \; \TotalTravelTime{\Solution} = \sum\limits_{\TripPath{\Trip} \in \Solution} \sum\limits_{\Arc \in \TripPath{\Trip}} \TravelTime{\Arc}{\Trip}(\Solution),
\end{align}
where, with slight abuse of notation, we make explicit that the travel time $\TravelTime{\Arc}{\Trip}$ experienced by trip $\Trip$ on arc $\Arc$ depends on all route and departure time assignments.

\subparagraph{Discussion}%
A few comments on our modeling approach are in order.
First, we assume that only \gls{amod} trips traverse the network, meaning the operator has full control over all trips in $\SetTrips$ to keep this basic problem setting concise. To extend the model to include background traffic, one can designate a subset of trips with routes and departure times not subject to optimization. This extension also enables the definition of two distinct objective functions: one that minimizes total system travel time, and another that minimizes only the travel time of the controlled trips. We present such an extended formulation in Section~\ref{sec:results}.

Second, incorporating route choice enhances behavioral realism beyond pure staggering models but increases computational complexity. To ensure scalability, we precompute a set of alternative routes—an approach aligned with route-based methods commonly used in dynamic traffic assignment \citep[cf.][]{LoSzeto2002}. As detailed in Sections~\ref{sec:methodology}~and~\ref{sec:case_study}, we design these routes to support effective balancing while preserving realistic travel patterns within the network.

Finally, although stochastic demand modeling would better capture real-world uncertainty, we assume full knowledge of demand at the central controller. This assumption aligns with the exploratory focus of our study, which aims to quantify the \emph{potential} benefits of integrated balanced and staggered routing in \gls{amod} systems and to provide an upper bound on travel time reductions that can serve as a benchmark for online methods.

\subparagraph{Relation to the discretized Vickrey’s bottleneck model}

Consider an arc with nominal travel time $\tau_a$. Let $f(t)$ denote the flow on this arc at time $t$, and let $\tau(t)$ represent the corresponding travel time, defined as:
\begin{equation}\label{eq:travel_time_linear}
    \tau(t) = \tau_a + (\phi \cdot \tau_a) \, f(t),
\end{equation}
where $\phi \in (0,1]$ is a scalar parameter that quantifies the sensitivity of travel time to congestion. We claim that, under Poisson arrivals, there exists a unique choice of \( \phi \) that renders \eqref{eq:travel_time_linear} an unbiased estimator of travel times in the discretized Vickrey bottleneck model \citep[cf.][]{OtsuboRapoport2008}.
To establish this claim, we first review the key features of the discretized Vickrey bottleneck model and show that, under a stable Poisson arrival process, it is equivalent to an M/D/1 queueing system \citep[cf.][]{Thomopoulos2012}. We then map our congestion model onto this queueing framework, allowing us to identify the unique value of $\phi$ that ensures an unbiased estimation of travel times. 

\begin{figure}[!t]
\centering
\subfigure[$t = t_1$]{%
    \def\svgwidth{0.2\linewidth}
    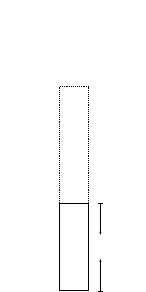}
\hspace{1cm}
\subfigure[$t = t_2$]{%
    \def\svgwidth{0.2\linewidth}
    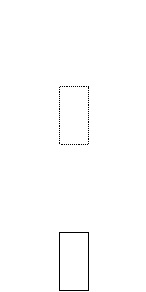}
\hspace{1cm}
\subfigure[$t = t_3$]{%
    \def\svgwidth{0.2\linewidth}
    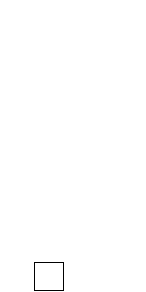}
\caption{Evolution of a bottleneck at three distinct time steps, \( t_1 \), \( t_2 \), and \( t_3 \), corresponding to the departure times of three distinct trips.}
\label{fig:bottleneck}
\end{figure}

Consider a finite set of trips traveling along a bottleneck road that connects a common origin to a common destination. The bottleneck has a nominal travel time $\tau_a$ and a capacity that restricts the number of vehicles simultaneously traversing it. We assume unit capacity, such that at most one trip can pass through the bottleneck at a time, while the remaining trips must wait in a \gls{fifo} queue. Figure~\ref{fig:bottleneck} schematizes the evolution of the system. The travel time for a trip departing at time \( t \) follows:
\begin{align}\label{eq:vickrey_behavior}
\tau^{\text{V}}(t)  = \begin{cases} 
\tau_a, & \text{if $f(t) =  0$} \\ 
\Delta \tau_a(t) + f(t) \cdot \tau_a, & \text{otherwise} 
\end{cases} 
\end{align} 
where \( \Delta \tau_a(t) \) represents the remaining time for the traveling trip to leave the bottleneck. 
Assume Poisson arrivals at the bottleneck with rate $\lambda \leq \tau_a^{-1}$. In this setting, \eqref{eq:vickrey_behavior} functions as a stable M/D/1 queueing system since the service time is constant and equal to $\tau_a$.
Let $\rho = \lambda \cdot \tau_a$ denote the traffic intensity of the system. The expected travel time at the bottleneck then follows:
\begin{equation}\label{eq:expected_travel_time_vickrey}
\mathbb{E}[\tau^V(t)] = \tau_a + \frac{\tau_a \cdot \rho}{2 (1 - \rho)}.
\end{equation}
We are now ready to show that, with an appropriate choice of parameters, Equation~\eqref{eq:travel_time_linear} yields expected travel times that match those of the discretized Vickrey bottleneck model in Equation~\eqref{eq:expected_travel_time_vickrey}. We refer the reader to Appendix~\ref{appendix:vickrey_proof} for the proof.
\begin{theorem}\label{thm:vickrey}
Consider an arc with nominal travel time $\tau_a$, where arrivals follow a Poisson process with rate $\lambda \leq \tau_a^{-1}$. Let the arc travel time be defined by Equation~\eqref{eq:travel_time_linear}. Then, by setting $\phi = (2 - \rho)^{-1}$, the expected travel time under Equation~\eqref{eq:travel_time_linear} coincides with the expected travel time given by Equation~\eqref{eq:expected_travel_time_vickrey}.
\end{theorem}
This theorem offers the following interpretation. When $\rho = 1$, the system operates at capacity, and the average traversal time $\Delta \tau_a$ equals $\tau_a$. In this case, setting $\phi = 1$ reproduces this behavior, ensuring that delay increases by $\tau_a$ per trip ahead. When $\rho < 1$, the average $\Delta \tau_a$ falls below $\tau_a$, so setting $\phi = (2 - \rho)^{-1}$ appropriately compensates for the reduced average delay.

%% file: figures/problem_setting/vickrey_example/vickrey_example_1.pdf_tex
\begingroup%
  \makeatletter%
  \providecommand\color[2][]{%
    \errmessage{(Inkscape) Color is used for the text in Inkscape, but the package 'color.sty' is not loaded}%
    \renewcommand\color[2][]{}%
  }%
  \providecommand\transparent[1]{%
    \errmessage{(Inkscape) Transparency is used (non-zero) for the text in Inkscape, but the package 'transparent.sty' is not loaded}%
    \renewcommand\transparent[1]{}%
  }%
  \providecommand\rotatebox[2]{#2}%
  \newcommand*\fsize{\dimexpr\f@size pt\relax}%
  \newcommand*\lineheight[1]{\fontsize{\fsize}{#1\fsize}\selectfont}%
  \ifx\svgwidth\undefined%
    \setlength{\unitlength}{69.84283447bp}%
    \ifx\svgscale\undefined%
      \relax%
    \else%
      \setlength{\unitlength}{\unitlength * \real{\svgscale}}%
    \fi%
  \else%
    \setlength{\unitlength}{\svgwidth}%
  \fi%
  \global\let\svgwidth\undefined%
  \global\let\svgscale\undefined%
  \makeatother%
  \begin{picture}(1,2.00263308)%
    \lineheight{1}%
    \setlength\tabcolsep{0pt}%
    \put(0.69883384,0.28727881){\makebox(0,0)[lt]{\lineheight{1.25}\smash{\begin{tabular}[t]{l}$\tau_a$\end{tabular}}}}%
    \put(0.49019921,0.28278335){\makebox(0,0)[lt]{\lineheight{1.25}\smash{\begin{tabular}[t]{l}1\end{tabular}}}}%
    \put(0,0){\includegraphics[width=\unitlength,page=1]{vickrey_example_1.pdf}}%
    \put(0.22687935,1.48536903){\makebox(0,0)[lt]{\lineheight{1.25}\smash{\begin{tabular}[t]{l}$f(t_1) = 0$\end{tabular}}}}%
  \end{picture}%
\endgroup%

%% file: figures/problem_setting/vickrey_example/vickrey_example_2.pdf_tex
\begingroup%
  \makeatletter%
  \providecommand\color[2][]{%
    \errmessage{(Inkscape) Color is used for the text in Inkscape, but the package 'color.sty' is not loaded}%
    \renewcommand\color[2][]{}%
  }%
  \providecommand\transparent[1]{%
    \errmessage{(Inkscape) Transparency is used (non-zero) for the text in Inkscape, but the package 'transparent.sty' is not loaded}%
    \renewcommand\transparent[1]{}%
  }%
  \providecommand\rotatebox[2]{#2}%
  \newcommand*\fsize{\dimexpr\f@size pt\relax}%
  \newcommand*\lineheight[1]{\fontsize{\fsize}{#1\fsize}\selectfont}%
  \ifx\svgwidth\undefined%
    \setlength{\unitlength}{69.84372116bp}%
    \ifx\svgscale\undefined%
      \relax%
    \else%
      \setlength{\unitlength}{\unitlength * \real{\svgscale}}%
    \fi%
  \else%
    \setlength{\unitlength}{\svgwidth}%
  \fi%
  \global\let\svgwidth\undefined%
  \global\let\svgscale\undefined%
  \makeatother%
  \begin{picture}(1,2.00260766)%
    \lineheight{1}%
    \setlength\tabcolsep{0pt}%
    \put(0.70564254,0.18500932){\makebox(0,0)[lt]{\lineheight{1.25}\smash{\begin{tabular}[t]{l}$\Delta \tau_a(t_2)$\end{tabular}}}}%
    \put(0.49045417,0.16947452){\makebox(0,0)[lt]{\lineheight{1.25}\smash{\begin{tabular}[t]{l}1\end{tabular}}}}%
    \put(0,0){\includegraphics[width=\unitlength,page=1]{vickrey_example_2.pdf}}%
    \put(0.69908647,0.68813424){\makebox(0,0)[lt]{\lineheight{1.25}\smash{\begin{tabular}[t]{l}$\tau_a$\end{tabular}}}}%
    \put(0.48584297,0.68363884){\makebox(0,0)[lt]{\lineheight{1.25}\smash{\begin{tabular}[t]{l}2\end{tabular}}}}%
    \put(0,0){\includegraphics[width=\unitlength,page=2]{vickrey_example_2.pdf}}%
    \put(0.21366695,1.48197444){\makebox(0,0)[lt]{\lineheight{1.25}\smash{\begin{tabular}[t]{l}$f(t_2) = 1$\end{tabular}}}}%
  \end{picture}%
\endgroup%

%% file: figures/problem_setting/vickrey_example/vickrey_example_3.pdf_tex
\begingroup%
  \makeatletter%
  \providecommand\color[2][]{%
    \errmessage{(Inkscape) Color is used for the text in Inkscape, but the package 'color.sty' is not loaded}%
    \renewcommand\color[2][]{}%
  }%
  \providecommand\transparent[1]{%
    \errmessage{(Inkscape) Transparency is used (non-zero) for the text in Inkscape, but the package 'transparent.sty' is not loaded}%
    \renewcommand\transparent[1]{}%
  }%
  \providecommand\rotatebox[2]{#2}%
  \newcommand*\fsize{\dimexpr\f@size pt\relax}%
  \newcommand*\lineheight[1]{\fontsize{\fsize}{#1\fsize}\selectfont}%
  \ifx\svgwidth\undefined%
    \setlength{\unitlength}{69.84283447bp}%
    \ifx\svgscale\undefined%
      \relax%
    \else%
      \setlength{\unitlength}{\unitlength * \real{\svgscale}}%
    \fi%
  \else%
    \setlength{\unitlength}{\svgwidth}%
  \fi%
  \global\let\svgwidth\undefined%
  \global\let\svgscale\undefined%
  \makeatother%
  \begin{picture}(1,2.00263308)%
    \lineheight{1}%
    \setlength\tabcolsep{0pt}%
    \put(0.53383463,0.07091786){\makebox(0,0)[lt]{\lineheight{1.25}\smash{\begin{tabular}[t]{l}$\Delta \tau_a(t_3)$\end{tabular}}}}%
    \put(0.31657818,0.06914916){\makebox(0,0)[lt]{\lineheight{1.25}\smash{\begin{tabular}[t]{l}1\end{tabular}}}}%
    \put(0,0){\includegraphics[width=\unitlength,page=1]{vickrey_example_3.pdf}}%
    \put(0.52651271,0.48363346){\makebox(0,0)[lt]{\lineheight{1.25}\smash{\begin{tabular}[t]{l}$\tau_a$\end{tabular}}}}%
    \put(0.31326744,0.479138){\makebox(0,0)[lt]{\lineheight{1.25}\smash{\begin{tabular}[t]{l}2\end{tabular}}}}%
    \put(0,0){\includegraphics[width=\unitlength,page=2]{vickrey_example_3.pdf}}%
    \put(0.52670408,1.08980169){\makebox(0,0)[lt]{\lineheight{1.25}\smash{\begin{tabular}[t]{l}$\tau_a$\end{tabular}}}}%
    \put(0.31519716,1.08530623){\makebox(0,0)[lt]{\lineheight{1.25}\smash{\begin{tabular}[t]{l}3\end{tabular}}}}%
    \put(0,0){\includegraphics[width=\unitlength,page=3]{vickrey_example_3.pdf}}%
    \put(0.04462939,1.48199388){\makebox(0,0)[lt]{\lineheight{1.25}\smash{\begin{tabular}[t]{l}$f(t_3) = 2$\end{tabular}}}}%
  \end{picture}%
\endgroup%

%% file: contents/3.Methodology.tex
\section{Methodology}\label{sec:methodology}
In this section, we present our methodology for solving large instances of Problem~\eqref{eq:problem}. Section~\ref{sec:alternative_routes} outlines the approach used to design alternative routes for balancing assignments. Section~\ref{sec:milp} introduces a mathematical formulation of the problem, including the notation and constraints that define feasible solutions. While this formulation serves as a valuable \emph{descriptive} tool, it is computationally impractical for large-scale instances. Therefore, in Section~\ref{sec:lns}, we propose a custom metaheuristic based on a \gls{lns} framework, designed to efficiently generate high-quality solutions at scale.

\subsection{Alternative route design}\label{sec:alternative_routes}  
Designing effective route alternatives is central to traffic balancing and congestion mitigation, yet it poses a critical trade-off between spatial diversity and route plausibility.

Solving a $k$-shortest path problem tends to produce highly overlapping routes, funneling demand onto a few shared arcs. This results in unrealistically high congestion and limits the \gls{amod} controller’s ability to distribute traffic effectively. In contrast, $k$-disjoint shortest paths enhance route diversity but often span significantly greater distances than the shortest route. Consequently, trips follow implausible detours that distort congestion estimates, and \gls{amod} controllers frequently reject them due to tight arrival time constraints.

To balance these limitations, we solve a \gls{kspwlo} problem for each trip \( \Trip \). Given a desired number of alternative routes \( k \) and a similarity threshold \( \theta\% \), the goal is to compute a set \( \SetPaths{\Trip} \) of \( k \) routes connecting the origin and destination of \( \Trip \), such that all route pairs are sufficiently dissimilar. Specifically, for all \( \{\Route{\Trip}{i}, \Route{\Trip}{j}\} \in \binom{\SetPaths{\Trip}}{2} \), we require:
\[
\Sim{\Route{\Trip}{i}}{\Route{\Trip}{j}} = \frac{\sum_{a \in \Route{\Trip}{i} \cap \Route{\Trip}{j}} \ell(a)}{\min\big\{ \ell(\Route{\Trip}{i}), \ell(\Route{\Trip}{j}) \big\}} \leq \theta,
\]
with \( \Route{\Trip}{i} \cap \Route{\Trip}{j} \) denoting shared arcs, and \( \ell(\cdot) \) the arc or path length.

To solve the \gls{kspwlo} problem efficiently, we adopt the \texttt{SVP}$^+$ algorithm \citep{ChondrogiannisBourosEtAl2020}, a performance-oriented heuristic designed for large-scale networks. \texttt{SVP}$^+$ restricts the search space to \emph{single-via paths} \citep{AbrahamDellingEtAl2013}, obtained by concatenating a shortest path from the source to an intermediate node with another shortest path from that node to the target. 
The algorithm first computes the shortest paths from the source to all nodes and from all nodes to the target. It then evaluates candidate paths formed by combining these partial paths and ranks them by length. As a final step, the algorithm incrementally adds paths to the solution set if they are sufficiently dissimilar from those already selected.
Finally, we note that the resulting path sets satisfy the following properties:
\begin{enumerate}
    \item[i)] they always include the shortest path;
    \item[ii)] all selected paths meet the similarity threshold $\theta$;
    \item[iii)] the selected paths are as short as possible within the restricted set.
\end{enumerate}

\subsection{Mathematical programming formulation}\label{sec:milp}
Our goal is to formulate Problem~\eqref{eq:problem} as a mathematical program. We begin by introducing the necessary variables, then formalize the constraints using an indicator function, and finally discuss the linearization of this function and explain why it renders the formulation intractable.

For every trip $\Trip \in \SetTrips$, we define the binary variable $\PathSelected{i}{\Trip} \in {0,1}$, which takes the value $1$ when the $i$-th route serves trip~$\Trip$ and remains $0$ otherwise, with $1 \leq i \leq k$. To account for increased travel time due to congestion, we introduce a delay variable $\DelayOnArc{\Arc}{\Trip}$ for each arc appearing in the alternative routes $\SetPaths{\Trip}$ of trip~$\Trip$. This delay satisfies $\DelayOnArc{\Arc}{\Trip} = \DelayFunction$, where $\DelayFunction$ is a convex, non-decreasing function of flow. We also introduce non-negative travel time variables $\TravelTime{\Arc}{\Trip} \geq 0$ for the same set of arcs and trips. These variables take positive values only on arcs that belong to the selected route and are zero otherwise. To enforce this behavior, together with the objective function minimizing the sum of the travel time variables, we apply the following big-M constraints:
\begin{align*}
\TravelTime{\Arc}{\Trip} 
&\geq \NominalTravelTime{\Arc} + \DelayOnArc{\Arc}{\Trip} 
- M \bigl( 1 - \PathSelected{i}{\Trip} \bigr), 
&\quad \forall \Trip \in \SetTrips,\; 
i \leq k,\; 
\Arc \in \Route{\Trip}{i}
\end{align*}

where $M$ denotes a sufficiently large constant. 

To model the temporal progress of each trip along its route, we introduce continuous variables $\Departure{\Arc}{\Trip}$ and $\Arrival{\Arc}{\Trip}$ for each arc in every alternative route, representing the start and completion times of the traversal of arc $\Arc$ by trip~$\Trip$. Finally, we denote $\FirstArcRoute{\Trip}$ and $\LastArcRoute{\Trip}$ as the first and last arcs of route $\Route{\Trip}{i}$, respectively. We formulate the problem as:
\begin{subequations}\label{eq:model}
\begin{align}
    &\min \;  \sum_{\Trip \in \SetTrips} \sum_{i \leq k} \sum_{\Arc \in \Route{\Trip}{i}} \TravelTime{\Arc}{\Trip}  \label{model:objective} \\
    &\text{s.t. } \nonumber \\
    & \textstyle \sum_{i  = 1}^{k} \PathSelected{i}{\Trip} = 1, 
    & \forall \Trip \in \SetTrips \label{model:route_selection} \\
    & \DelayOnArc{\Arc}{\Trip} = \DelayFunction, 
    & \forall \Trip \in \SetTrips, i \leq k, \Arc \in \Route{\Trip}{i} \label{model:delay_one} \\
    & \TravelTime{\Arc}{\Trip} \geq \NominalTravelTime{\Arc} + \DelayOnArc{\Arc}{\Trip} - M (1 - \PathSelected{i}{\Trip}), 
    & \forall \Trip \in \SetTrips, i \leq k, \Arc \in \Route{\Trip}{i} \label{model:delay_two} \\
    & \Flow{\Arc}{\Trip} = \textstyle \sum_{\SecondTrip \in \SetTrips \setminus\!\{\Trip\}} \mathbb{I}(\Departure{\Arc}{\SecondTrip} \!\leq\! \Departure{\Arc}{\Trip}  \!<\! \Arrival{\Arc}{\SecondTrip}), 
    & \forall \Trip \in \SetTrips, i \leq k, \Arc \in \Route{\Trip}{i} \label{model:indicator} \\
    & \Arrival{\Arc}{\Trip} = \Departure{\Arc}{\Trip} + \TravelTime{\Arc}{\Trip}, 
    & \forall \Trip \in \SetTrips, i \leq k, \Arc \in \Route{\Trip}{i} \label{model:continuity} \\
    & \TripEarliestDeparture{\Trip} \leq \Departure{\Arc}{\Trip} \leq \TripEarliestDeparture{\Trip} + \TripMaxStaggering{\Trip}, 
    & \forall \Trip \in \SetTrips, i \leq k, \Arc = \FirstArcRoute{\Trip} \label{model:earliest_departure} \\
    & \Arrival{\Arc}{\Trip} \leq \TripLatestArrival{\Trip}, 
    & \forall \Trip \in \SetTrips, i \leq k, \Arc = \LastArcRoute{\Trip} \label{model:latest_arrival} \\
    & \Departure{\Arc}{\Trip},\Arrival{\Arc}{\Trip}, \Flow{\Arc}{\Trip}, \DelayOnArc{\Arc}{\Trip}, \TravelTime{\Arc}{\Trip} \geq 0, 
    & \forall \Trip \in \SetTrips, i \leq k, \Arc \in \Route{\Trip}{i} \label{model:continuous_variables} \\
    & \PathSelected{i}{\Trip} \in \{0, 1\}, 
    & \forall \Trip \in \SetTrips, i \leq k \label{model:binary_variables}
\end{align}
\end{subequations}
The objective in \eqref{model:objective} minimizes the total travel time of the fleet. Constraint~\eqref {model:route_selection} ensures the selection of exactly one route for each trip. Constraints \eqref{model:delay_one}-\eqref{model:delay_two}, together with the definitions of the variables $\DelayOnArc{\Arc}{\Trip}$ and $\TravelTime{\Arc}{\Trip}$ in \eqref{model:continuous_variables}, calculate the travel time for trips on the arcs of their selected routes and ensure that travel time is zero for non-selected routes.
Constraint~\eqref{model:indicator} determines the number of trips encountered by a trip $\Trip$ upon entering arc~$\Arc$, using the indicator function $\mathbb{I}(\cdot)$. Constraint~\eqref{model:continuity} propagates the departure times of trips along their routes. Constraints \eqref{model:earliest_departure}-\eqref{model:latest_arrival} enforce time window restrictions for each trip. Finally, Constraints~\eqref{model:continuous_variables}-\eqref{model:binary_variables} define the domains of the decision variables.

Directly solving the formulation presented is impractical due to the presence of the indicator function in \eqref{model:indicator}. To address this, we replace it with an additional set of big-M constraints, obtaining a mixed-integer convex program. For a detailed description of this approach, we refer the reader to Section~3.1 of \citet{CoppolaHiermannEtAl2024}, which provides a procedure that applies almost directly to our extended setting. 
The key distinction in this case lies in the requirement to define a set of big-M constraints for all pairs of trips that may \emph{potentially} traverse a given arc, that is, for trips for which the arc belongs to at least one alternative route.
Even under a simple specification of $\DelayFunction$, such as a linear or piecewise linear function (resulting in a mixed-integer linear program), the number of big-M constraints required for linearization is significantly higher compared to a setting with fixed routes, rendering the formulation computationally intractable.
Therefore, rather than attempting to solve this \gls{milp} directly, we develop a metaheuristic approach capable of finding feasible solutions without the need to explicitly define each of these constraints.

\subsection{Algorithm}\label{sec:lns}

To solve large instances of Problem~\eqref{eq:problem}, we develop a custom metaheuristic algorithm. 
\begin{algorithm}[!b]
\footnotesize
\caption{Framework overview}\label{algo:hierarchical}
\hspace{\algorithmicindent}\textbf{Input:} Set of trips $\SetTrips$, parameters for LNS and routing
\vspace{-0em}
\begin{algorithmic}[1]
\itemsep-0em
    \State $\Solution_{\text{RDUO}} \gets \textsc{buildRDUO}(\SetTrips)$
    \State $\Solution_{\text{greedy}} \gets \textsc{greedyAssignment}(\SetTrips)$
    \State $\Solution \gets \arg\min\{\SolutionCost{\Solution_{\text{RDUO}}}, \SolutionCost{\Solution_{\text{greedy}}}\}$
    \State $\Solution \gets \textsc{LNS}(\Solution)$
    \State \textbf{return} $\Solution$
\end{algorithmic}
\end{algorithm}
\begin{figure}[!b]
\centering
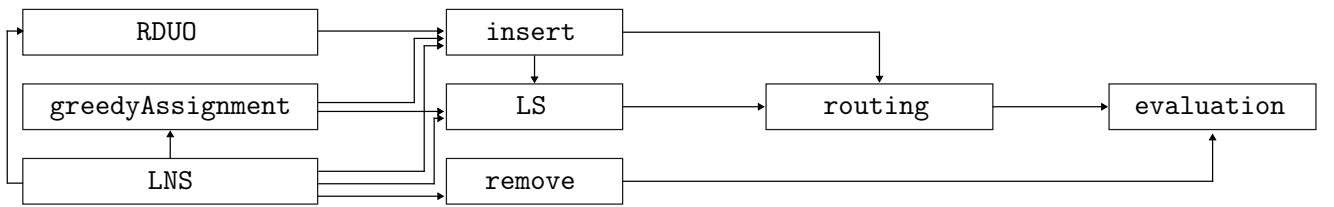
\caption{Algorithm operators in hierarchical order. Solution-finding modules (left) generate or improve solutions by invoking intermediate-level operators (center), which, in turn, rely on lower-level logic for routing trips and evaluating the resulting solutions (right).}
\label{fig:operators}
\end{figure}
Algorithm~\ref{algo:hierarchical} outlines the overall structure of our solution approach, while Figure~\ref{fig:operators} illustrates the hierarchical dependencies among the algorithm’s operators.
At a high level, the algorithm:
\begin{enumerate}
\item[i)] Builds a \gls{rduo}.
\item[ii)] Builds a solution that greedily routes trips to minimize the system’s total travel time.
\item[iii)] Chooses the better of i) and ii) to initialize a \gls{lns} for further improvement.
\end{enumerate}
These solution-finding modules rely on several intermediate-level operators:
\begin{enumerate}
    \item[i)] An \emph{insert} operator assigning trips to partial solutions.    
    \item[ii)] A \emph{local search} that adjusts routes and departure times for inserted trips.
    \item[iii)] A \emph{remove} operator that eliminates trips from a solution.
\end{enumerate}
The \emph{insert} and \emph{local search} operators share a common low-level logic for identifying route and departure time combinations. After any insertion or removal, the algorithm evaluates the updated solution. Since evaluation constitutes the most computationally expensive step, we develop an efficient procedure to perform it.
As we allow the algorithm to explore intermediate infeasible solutions, the algorithm ranks solutions \( \Solution \) using the following cost function:
\[
\SolutionCost{\Solution} = \TotalDelay{\Solution} + \alpha \cdot \Infeasibility{\Solution},
\]
where \( \TotalDelay{\Solution} \) is the \emph{total delay} and \( \Infeasibility{\Solution} \) is the \emph{total infeasibility} associated with \( \Solution \), and \( \alpha \) controls the penalty applied to infeasibility. The total delay \( \TotalDelay{\Solution} \) measures the excess travel time relative to the shortest possible free-flow travel times. It consists of two components: a \emph{detour delay}, the extra time from selecting routes longer in distance than the shortest routes, and a \emph{congestion delay}, the additional time caused by traffic congestion along the chosen routes. The infeasibility term \( \Infeasibility{\Solution} \) quantifies the total excess travel time of trips arriving later than their latest allowed arrival times.

The remainder of this section is organized as follows. First, we describe the solution-finding modules in Section~\ref{sec:high_level_operatators}, followed by the intermediate operators in Section~\ref{sec:intermediate_operators}. Finally, we present the lower-level operators in Section~\ref{sec:low_level_operators}.

\subsubsection{Solution-finding modules}\label{sec:high_level_operatators}
\subparagraph{Reactive dynamic user optimum}
The \gls{rduo} models users departing at their preferred times and independently selecting routes that minimize their travel time based on prevailing traffic conditions at the time of departure. As users do not account for future congestion, their chosen routes remain fixed once the journey commences. The \gls{rduo} computation starts with an empty solution $\Solution = (\SolutionPaths, \SolutionTimes)$, where $\SolutionPaths = \emptyset$ and $\SolutionTimes = \emptyset$ denote the initial route and departure time assignment sets.
The procedure inserts trips sequentially, in order of increasing departure times. For each trip, $\SolutionTimes \gets \SolutionTimes \cup \TripEarliestDeparture{\Trip}$ and $\SolutionPaths \gets \SolutionPaths \cup \Route{\Trip}{e}$, where $\Route{\Trip}{e}$ minimizes travel time at $\TripEarliestDeparture{\Trip}$ given previously assigned trips. {We provide the pseudocode in Appendix~\ref{appendix:baseline-pseudocode}}.

\subparagraph{Greedy assignment}
The greedy assignment sequentially adds trips to an initially empty solution, selecting staggering and routing options to minimize total system travel time. As in the \gls{rduo}, the algorithm starts from an empty solution and inserts trips one at a time. This process may yield infeasible solutions, as the greedy insertion of early trips can cause later trips to violate their latest arrival times. To address this, the algorithm prioritizes the insertion of trips with the most restrictive latest arrival times, thereby promoting feasibility while still minimizing total delay. If infeasibilities remain at the end of the insertion procedure, the algorithm invokes a \emph{local search}, which assigns new route and departure time combinations to late-arriving trips as a repair mechanism. {We include the corresponding pseudocode in Appendix~\ref{appendix:baseline-pseudocode}.}

\subparagraph{Large neighborhood search}
\begin{algorithm}[!t]
\footnotesize
\caption{\textsc{largeNeighborhoodSearch}}\label{algo:lns}
\hspace{\algorithmicindent}\textbf{Input:} Set of trips $\SetTrips$, $\Solution_{\text{greedy}}$, $\Solution_{\text{RDUO}}$, removal pool size $x$, destruction percentage $y$, max cycles $z$
\vspace{-0em}
\begin{algorithmic}[1]
\itemsep-0em
    \State $\Solution \gets \arg\min\{\SolutionCost{\Solution_{\text{RDUO}}}, \SolutionCost{\Solution_{\text{greedy}}}\}$ \label{line:init_solution}
    \While{True} \label{line:while_start}
        \State $\texttt{improved} \gets \texttt{False}$ \label{line:init_improved}
        \For{$\textit{operator}$ in \{$\texttt{COSTLY\_TRIPS, UNTOUCHED\_TRIPS}$\}} \label{line:loop_destroy_ops}
            \State $\tilde{\Solution} \gets \Solution$ 
            \Comment{Best candidate for this operator} \label{line:init_best_candidate}
            \State $\SetTrips_{\text{p}} \gets \textsc{selectRemovalPool}(\Solution, x, \textit{operator})$ \label{line:pool_per_operator}
            \For{$i$ in $\{1 \dots z\}$} \label{line:destroy_cycles}
                \State $\SetTrips_{\text{d}} \gets \textsc{randomSample}(\SetTrips_{\text{p}}, y)$ \label{line:destroy_sample}
                \For{$\textit{order}$ in \textsc{shuffle}(\texttt{EARLIEST, LATEST, DELAY})} \label{line:reinsertion_orders}
                    \State $\Solution' \gets \textsc{remove}(\Solution, \SetTrips_{\text{d}})$ \label{line:remove}
                    \State $\Solution' \gets \textsc{insert}(\Solution', \SetTrips_{\text{d}}, \textit{order})$ \label{line:insert}
                    \If{$\SolutionCost{\Solution'} < \SolutionCost{\Solution}$} \label{line:check_improvement}
                        \State $\Solution \gets \Solution'$ \label{line:update_solution}
                        \State $\texttt{improved} \gets \texttt{True}$ \label{line:set_improved}
                        \State \textbf{goto} line~\ref{line:pool_per_operator} \Comment{Select new removal pool with this operator} \label{line:goto_restart_operator}
                    \ElsIf{$\SolutionCost{\Solution'} < \SolutionCost{\tilde{\Solution}}$} \label{line:check_candidate}
                        \State $\tilde{\Solution} \gets \Solution'$ \label{line:update_candidate}
                    \EndIf
                \EndFor \label{line:end_reinsertion_orders}
            \EndFor \label{line:end_destroy_cycles}
            \State $\SetTrips_{\text{c}} \gets \textsc{selectChangedTrips}(\tilde{\Solution}, \Solution)$ \label{line:select_changed_trips}
            \State $\Solution' \gets \textsc{localSearch}(\tilde{\Solution}, \SetTrips_{\text{c}})$ \label{line:local_search}
            \If{$\SolutionCost{\Solution'} < \SolutionCost{\Solution}$} \label{line:check_local_improvement}
                \State $\Solution \gets \Solution'$ \label{line:update_after_local}
                \State $\texttt{improved} \gets \texttt{True}$ \label{line:set_improved_after_local}
            \EndIf
        \EndFor \label{line:end_operator_loop}
        \If{\textbf{not} \texttt{improved}} \label{line:check_any_improvement}
            \State \textbf{break} \Comment{No improvement in any operator} \label{line:terminate}
        \EndIf
    \EndWhile \label{line:end_while}
    \State \textbf{return} $\Solution$ \label{line:return_solution}
\end{algorithmic}
\end{algorithm}

The \gls{lns} iteratively improves a solution by progressively selecting destroy and repair operators applied to subsets of trips, thereby gradually exploring better combinations of trip routes and departure times. Algorithm~\ref{algo:lns} summarizes the steps of the \gls{lns}. The algorithm begins by selecting the lower-cost solution between the \gls{rduo} and the greedy assignment as the \emph{incumbent} (l.\ref{line:init_solution}). Each iteration begins by selecting the \emph{remove} operator, applying one of the following strategies:
\begin{enumerate}
\item[i)] removal of the trips with the highest individual cost;
\item[ii)] removal of the trips scheduled earliest along their shortest-distance route.
\end{enumerate}
These two operators target distinct classes of trips, yet both offer strong potential for improvement: the first removes the most costly trips that contribute significantly to the overall solution cost, while the second focuses on reassigning trips with untapped potential for improved balancing or staggering.
After the operator selection, the algorithm builds a removal pool (l.\ref{line:pool_per_operator}) consisting of $x\%$ of all trips, along with any infeasible ones. Up to $z$ attempts (l.\ref{line:destroy_cycles}), it randomly samples $y\%$ of the pool (l.\ref{line:destroy_sample}) and removes the selected trips from the current solution (l.\ref{line:remove}).
It then uses the \emph{insert} operator to reinsert the removed trips, which explores up to three insertion orders applied in random sequence (l.\ref{line:reinsertion_orders}--\ref{line:insert}):
%
If any reinsertion produces a new incumbent solution, the algorithm immediately restarts the destroy-repair cycle by selecting a new removal pool using the last remove operator invoked (l.\ref{line:update_solution}--\ref{line:goto_restart_operator}).
If no incumbent updates occur after $z$ cycles, the algorithm invokes the \emph{local search} procedure to improve the quality of the best solution found from this phase (l.\ref{line:select_changed_trips}--\ref{line:set_improved_after_local}).
While effective, the local search is computationally expensive, as it operates on full solutions. Therefore, during its first invocation, the algorithm applies the local search to all trips in $\SetTrips$, and restricts it in subsequent calls to trips whose delay or infeasibility has changed since the previous invocation. To support this, the algorithm maintains two vectors storing each trip’s total delay and infeasibility, updated after every \emph{local search} call. 

\noindent After completing the {local search}, the algorithm proceeds to the next destroy operator. If at least one incumbent update occurs after evaluating all destroy operators, the algorithm restarts from the first operator (l.\ref{line:end_operator_loop}).
Otherwise, it terminates (l.\ref{line:terminate}) and returns the best solution found (l.\ref{line:return_solution}).

To help the \gls{lns} balance feasibility and solution quality, we utilize an adaptive infeasibility weight~$\alpha$, which penalizes constraint violations in the solution cost. Starting with $\alpha$ equal to ten, we divide it by ten if ten consecutive \emph{insert} calls yield feasible solutions, thereby reducing the penalty and encouraging exploration near the boundary of the feasible region. Conversely, if ten consecutive calls result in infeasible solutions, we multiply $\alpha$ by ten to guide the search back toward feasibility. To prevent both excessive penalization and loss of feasibility, we constrain $\alpha$ to the interval $[10^{-2}, 10^{3}]$.

\subsubsection{Intermediate-level operators}\label{sec:intermediate_operators}

\subparagraph{Insert}
The \emph{insert} operator takes as input a partial solution, a set of trips to insert, and a rule specifying the insertion order. Operating on one trip at a time, it tests and evaluates a series of alternative route and departure time combinations $(\Route{\Trip}{}, \Departure{\Trip}{})$, identified using the low-level \emph{routing} logic. It selects the combination yielding the lowest cost and updates the solution by adding the chosen route and departure time, updating $\SolutionPaths \gets \SolutionPaths \cup \Route{\Trip}{}$ and $\SolutionTimes \gets \SolutionTimes \cup \Departure{\Trip}{}$.

\subparagraph{Local search}
The \emph{local search} operator takes as input a solution $\Solution$ and reassigns routes and departure times for a specified subset of its trips.
Unlike \emph{insert}, it does not alter the number of trips in the solution. Proceeding in order of earliest permissible departure times $\TripEarliestDeparture{\Trip}$, the operator iteratively generates a set of alternative route and departure time combinations $(\Route{\Trip}{}, \Departure{\Trip}{})$ for each trip using the low-level \emph{routing} logic. It evaluates all combinations, selects the one with the lowest cost, and updates the solution by replacing the trip’s current route and departure time accordingly.

\subparagraph{Remove}
The \emph{remove} operator takes as input a solution and a set of trips to remove. For each trip, the operator deletes the corresponding route and departure time from the solution, updating $\SolutionPaths \gets \SolutionPaths \setminus \Route{\Trip}{}$ and $\SolutionTimes \gets \SolutionTimes \setminus \Departure{\Trip}{}$. After removing all specified trips, it evaluates the resulting solution.

\subsubsection{Low-level operators}\label{sec:low_level_operators}

\subparagraph{Evaluation}\label{sec:schedule_evaluation}
Every modification to a solution — whether a trip insertion, removal, departure time shift, or trip rerouting — alters the solution cost, which the algorithm must quantify. Computing $\SolutionCost{\Solution}$ requires assessing both the congestion experienced by trips and any late arrivals at their destinations. To this end, the algorithm assigns each solution a \emph{schedule} $\Schedule$, which specifies, for each trip, the set of arc departure times $\Departure{\Arc}{\Trip}$ determined by the selected start times and routes in $\Solution$.
To determine the actual travel time of a trip $\Trip$ on arc $\Arc$, the algorithm first calculates the flow on the arc, $\Flow{\Arc}{\Trip}$, and then applies Equation~\eqref{eq:travel_time}. The algorithm calculates $\Flow{\Arc}{\Trip}$ by counting the number of \emph{conflicts} with other trips $\SecondTrip$ on the same arc.
A conflict occurs when $\Trip$ departs from $\Arc$ after $\SecondTrip$ has started but before $\SecondTrip$ has completed its traversal of the arc, including any delays; that is, when 
$\Departure{\Arc}{\SecondTrip} \leq \Departure{\Arc}{\Trip} < \Arrival{\Arc}{\SecondTrip}$.
Previous work on staggered routing proposes the \textsc{constructSchedule} procedure to compute a schedule from scratch; we provide its details in Appendix~\ref{appendix:construct_schedule}. While this function still finds application in our algorithm—for instance, after the \emph{remove} operator deletes a subset of trips from the current solution— it becomes inefficient to evaluate a large number of new assignments, as it rebuilds the entire schedule even when the change affects only the trip under consideration.
\begin{algorithm}[!b]
\footnotesize 
\caption{\textsc{updateSchedule}}\label{algo:update_schedule}
\hspace{\algorithmicindent}\textbf{Input:} Schedule $\Schedule$, solution $\Solution$
\begin{algorithmic}[1]
\itemsep-0em  
\State $\Queue \gets \textsc{getChanges}(\Schedule,\Solution)$ \label{line:stagger_departure}
\While{$\Queue$ is not empty} 
    \State $(\Trip,\Arc,\Departure{\Arc}{\Trip}) \gets \Queue.\textsc{pop}()$ \label{line:processing}
    \State $\Schedule, \SecondTrip \gets \textsc{processTrip}(\Trip,\Arc,\Departure{\Arc}{\Trip})$ \label{line:check_new_conflicts}
    \If{$\SecondTrip \neq \textsc{false}$}\label{line:trip_flag_is_false}
        \State $\Queue \gets \textsc{insertTrips}(\Trip,\SecondTrip,\Arc,\Departure{\Arc}{\Trip},\Departure{\Arc}{\SecondTrip})$ \label{line:reinsert_new_conflict}
    \Else \label{line:stag_dep_else}
        \While{$\SecondTrip = \textsc{false}$}
            \State $\Arc' \gets \textsc{getNextArc}(\Trip,\Arc, \Solution)$
            \State $\Schedule,\SecondTrip \gets \textsc{processTrip}(\Trip,\Arc',\Departure{\Arc'}{\Trip})$
            \State $\Arc \gets \Arc'$
            \If{$\SecondTrip \neq \textsc{false}$}
                \State $\Queue \gets \textsc{insertTrips}(\Trip,\SecondTrip,\Arc,\Departure{\Arc}{\Trip},\Departure{\Arc}{\SecondTrip})$ \label{line:insert_new_trips_2} 
            \EndIf
        \EndWhile
    \EndIf
\EndWhile
\State \textbf{return} $\Schedule$
\end{algorithmic}
\end{algorithm}
Algorithm~\ref{algo:update_schedule} outlines our proposed approach for efficiently updating the schedule during conflict resolution. The core challenge that the algorithm addresses lies in identifying how a modification to a trip’s route or departure time propagates through the schedule. In the best case, the change affects only the modified trip; in the worst case, it triggers a cascade of updates that shift the departure times of many other trips. Our method extends the \textsc{updateSchedule} procedure from previous work on staggered routing by supporting multiple route alternatives and introducing further refinements to improve computational efficiency.

To initiate the search, we begin with an existing schedule that serves as a benchmark for evaluating changes resulting from start time or route adjustments. At the core of the algorithm is a priority queue $\Queue$, initially empty. This queue holds tuples that include trip-arc-departure labels $(\Trip, \Arc, \Departure{\Arc}{\Trip})$. To manage $\Queue$, we sort the tuples in ascending order based on their start times~$\Departure{\Arc}{\Trip}$. We fill the queue with all trips potentially impacted by modifications in their scheduled departures~(l.\ref{line:stagger_departure}). The approach to populating this queue varies based on the specific operation being executed:
\begin{enumerate}
    \item[i) Trip insertion or staggering:] We add an initial tuple containing the trip, the first arc of its assigned route, and its scheduled departure time.
    \item[ii) Trip rerouting:] We add all trips affected by the change on the previous route to the queue, along with a tuple for the rerouted trip at its earliest entry time on the new route. Specifically, for each arc $\Arc$ along the previous route, we insert $(\SecondTrip, \Arc, \Departure{\Arc}{\SecondTrip})$ for every trip $\SecondTrip$ that conflicts with the rerouted trip in the original schedule.
\end{enumerate}
We start processing trips in $\Queue$ (l.\ref{line:processing}) by checking whether the trip's updated departure time $\Departure{\Arc}{\Trip}$ creates a change in the arc entry times of another trip $\SecondTrip$ (l.\ref{line:check_new_conflicts}). In this case, we insert $\SecondTrip$ into $\Queue$. If this insertion disrupts the queue order — specifically, if $\Departure{\Arc}{\SecondTrip} < \Departure{\Arc}{\Trip}$ — we flag the conflict by setting $\SecondTrip \neq \textsc{false}$, reinsert $\Departure{\Arc}{\Trip}$ into $\Queue$, and restart the evaluation from line~\ref{line:trip_flag_is_false} to~\ref{line:reinsert_new_conflict}. If $\Departure{\Arc}{\Trip}$ does not induce changes ($\SecondTrip = \textsc{false}$) in other trips, we recursively propagate the change in trip $\Trip$'s arc entry times on subsequent arcs until the trip reaches its destination~(l.\ref{line:stag_dep_else}-\ref{line:insert_new_trips_2}). 
Our overall update ends once $\Queue$ is empty and all arc entry times have been correctly adjusted.

\begin{figure}[!b]
\centering
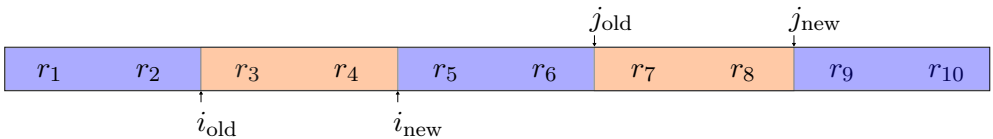
\caption{Identification of trips to activate on an arc, given a set of labels sorted by departure time. The indices \( i_\mathrm{old} \) and \( i_\mathrm{new} \) mark the previous and updated departure times of trip \( \Trip \), while \( j_\mathrm{old} \) and \( j_\mathrm{new} \) correspond to its previous and updated arrival times. Trips before (\( r_1 \), \( r_2 \)), after (\( r_9 \), \( r_{10} \)), or overlapping (\( r_5 \), \( r_6 \)) remain unaffected (highlighted in blue). The algorithm pushes the remaining trips (highlighted in orange) to \( \Queue \) because their conflict status with \( \Trip \) changes: either the update resolves a previous conflict (\( r_3 \), \( r_4 \)) or introduces a new one (\( r_7 \), \( r_8 \)).}
\label{fig:index_in_set}
\end{figure}

\begin{figure}[!b]
\centering
\def\svgwidth{0.22\linewidth}
\subfigure[Alternative routes]{%
    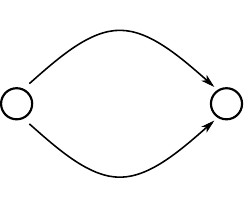}
\hfill
\def\svgwidth{0.22\linewidth}
\subfigure[Initial schedule]{%
    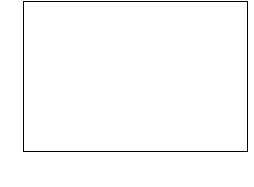}
\hfill
\def\svgwidth{0.22\linewidth}
\subfigure[Staggering]{%
    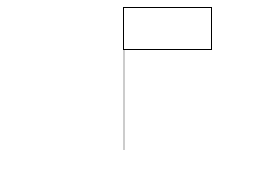}
\hfill
\def\svgwidth{0.22\linewidth}
\subfigure[Rerouting]{%
    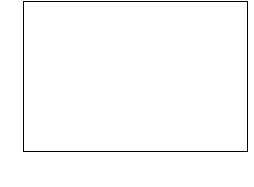}
\caption{Examples of new start times and route identification on alternative routes \( a \) and \( b \), with \( \tau_a < \tau_b \).}
\label{fig:ls_moves}
\end{figure}

To efficiently determine which trips to push to $\Queue$ and count the resulting conflicts when processing trip $\Trip$ on arc $\Arc$, we maintain auxiliary sets $\SortedTripsOnArc{\Arc}$, which store the labels of all trips currently assigned to arc $\Arc$, sorted by departure time. We partition this set into:
\begin{enumerate}
    \item[i)] $\SortedInactiveTripsOnArc{\Arc}$: labels of \emph{inactive} trips, i.e., those not yet enqueued in $\Queue$.
    \item[ii)] $\SortedActiveTripsOnArc{\Arc}$: labels of \emph{active} trips, i.e., those already present in $\Queue$.
\end{enumerate}
Initially, $\SortedInactiveTripsOnArc{\Arc} = \SortedTripsOnArc{\Arc}$ and $\SortedActiveTripsOnArc{\Arc} = \emptyset$. As the algorithm activates new trips, it removes their labels from each $\SortedInactiveTripsOnArc{\Arc}$ for all arcs in their route and adds them to the corresponding $\SortedActiveTripsOnArc{\Arc}$. We deliberately avoid maintaining $\SortedActiveTripsOnArc{\Arc}$ in sorted order, since active departure times may change multiple times due to propagation effects, making constant re-sorting computationally expensive. We then compute the flow on arc $\Flow{\Arc}{\Trip}$ as:
\[
\Flow{\Arc}{\Trip} = \InactiveFlow{\Arc}{\Trip} + \ActiveFlow{\Arc}{\Trip},
\]
where $\InactiveFlow{\Arc}{\Trip}$ and $\ActiveFlow{\Arc}{\Trip}$ represent the flows contributed by inactive and active trips, respectively. 
Since $\SortedActiveTripsOnArc{\Arc}$ is unsorted, we evaluate $\ActiveFlow{\Arc}{\Trip}$ by iterating over all active trips:
\[
{\ActiveFlow{\Arc}{\Trip}} = \sum_{\SecondTrip \in \SortedActiveTripsOnArc{\Arc} \setminus \{\Trip\}} \mathbb{I}({\Departure{\Arc}{\SecondTrip}} \leq {\Departure{\Arc}{\Trip}} < {\Arrival{\Arc}{\SecondTrip}}).
\]
{Conversely, we exploit the sorted structure of $\SortedInactiveTripsOnArc{\Arc}$ to compute $\InactiveFlow{\Arc}{\Trip}$ through an efficient counting procedure: we locate the latest trip $\SecondTrip$ in $\SortedInactiveTripsOnArc{\Arc}$ with $\Departure{\Arc}{\SecondTrip}\! < \Departure{\Arc}{\Trip}$ and traverse the set backward from that position, incrementing a counter at each step, until we find a trip with $\Arrival{\Arc}{\SecondTrip} \leq \Departure{\Arc}{\Trip}$. When the \gls{fifo} property does not hold on an arc, the procedure may yield occasional inaccuracies in delay calculations. We therefore repair the schedule upon completion of each operator by invoking \textsc{constructSchedule}.} Finally, the algorithm calculates the corresponding delay using the delay function $\DelayFunction$ defined in Equation~\eqref{eq:travel_time}.\\
The algorithm adds trip $\SecondTrip$ to $\Queue$ if an update to $\Trip$ either introduces a new conflict or resolves an existing one between $\Trip$ and $\SecondTrip$. To identify such trips efficiently, we leverage the sorted structure of $\SortedInactiveTripsOnArc{\Arc}$.
We define four indices in $\SortedInactiveTripsOnArc{\Arc}$, each representing the insertion position of $\Trip$ based on:
\begin{enumerate}
    \item[i)] its departure time \emph{before} the update ($i_\mathrm{old}$).
    \item[ii)] its departure time \emph{after} the update ($i_\mathrm{new}$).
    \item[iii)] its arrival time \emph{before} the update ($j_\mathrm{old}$).
    \item[iv)] its arrival time \emph{after} the update ($j_\mathrm{new}$).
\end{enumerate}
Trips located between $i_\mathrm{old}$ and $i_\mathrm{new}$ (departure range) or between $j_\mathrm{old}$ and $j_\mathrm{new}$ (arrival range) become active, as their conflict status with $\Trip$ has changed. Trips outside both ranges remain inactive. Those within the intersection of the two ranges also stay inactive, since their conflict status with $\Trip$ is unaffected. Figure~\ref{fig:index_in_set} illustrates these cases with an example.
\subparagraph{Routing}  
This operator identifies new combinations of start times and routes, either to insert a trip with minimal impact on the solution cost (\emph{insert}) or to improve an existing assignment (\emph{local search}). First, it iteratively refines the departure time on the trip's initial route to maximize the number of resolved conflicts. As an additional refinement step, the method removes portions of the applied departure shifts that do not contribute to reducing the solution cost, thereby promoting parsimonious staggering. Then, it reroutes the trip as early as possible on each alternative route and repeats the iterative start time identification procedure. After evaluating all routes, it selects the route--departure time combination that yields the lowest cost. Figure~\ref{fig:ls_moves} illustrates this process using a small example.
Given a trip $\Trip$, the algorithm analyzes conflicts along its current route $\TripPath{\Trip}$ on an arc-by-arc basis, based on its current departure time. On each arc $\Arc$ of the route, the algorithm evaluates whether delaying $\Trip$ can resolve conflicts without introducing new ones. These evaluations are performed independently for each arc and ignore how staggering affects delays on preceding arcs, which may, in turn, influence downstream conflicts. As a result, the conflict-based evaluation serves only as a heuristic proxy for the true impact on overall solution quality.

To detect conflicts efficiently, the algorithm uses the sorted set $\SortedTripsOnArc{\Arc}$, which orders trip labels on arc $\Arc$ by their departure times. It locates $\Trip$ in $\SortedTripsOnArc{\Arc}$ and identifies the immediately preceding trip $\SecondTrip$. If $\Arrival{\Arc}{\SecondTrip} < \Departure{\Arc}{\Trip}$, no overlap occurs, and the algorithm proceeds to the next arc. Otherwise, it delays $\Trip$ just enough to eliminate the conflict, as long as the new departure time remains within its feasible window:
\[
\StartTime{\Trip} \gets \min\left(
\StartTime{\Trip} + \left(\Arrival{\Arc}{\SecondTrip} - \Departure{\Arc}{\Trip}\right)\hspace{-.2em}, \TripEarliestDeparture{\Trip} + \TripMaxStaggering{\Trip}
\right).
\]
If trips on the arc follow a \gls{fifo} discipline, resolving the conflict with $\SecondTrip$ also eliminates all earlier conflicts on the arc. While this assumption does not always hold, the algorithm adopts it as a computational shortcut.
After updating the departure time of $\Trip$, the algorithm first estimates its new arrival time, incorporating both the applied staggering and any delay reductions due to resolved conflicts with preceding trips. It then scans forward through $\SortedTripsOnArc{\Arc}$ to detect changes in conflict status with subsequent trips $\Trip'$. A new conflict arises if $\Arrival{\Arc}{\Trip} > \Departure{\Arc}{\Trip'}$, whereas no such conflict existed before the update. Conversely, a conflict is resolved if the reduced delay leads to $\Arrival{\Arc}{\Trip} \leq \Departure{\Arc}{\Trip'}$.
The algorithm aggregates all resolved conflicts as $f_{\Arc}^{-}$, and all newly introduced conflicts as $f_{\Arc}^{+}$, and computes the net conflict gain as:
\[
\Delta f_{\Arc} = f_{\Arc}^{-} - f_{\Arc}^{+}.
\]
This value heuristically indicates the benefit of setting the trip’s start time to $\StartTime{\Trip}$ on arc $\Arc$. After evaluating all arcs, the algorithm selects the departure time candidate with the highest net conflict gain $\Delta f_{\Arc}$. If the adjustment improves the solution, the algorithm accepts it and continues refining the start time on the same route until it finds no additional improvement. 
Before concluding the analysis of the current route, the method trims any part of the applied departure shift that does not reduce the solution cost. This refinement brings two main advantages: it unlocks further staggering opportunities for the same trip in later iterations and minimizes overlap with downstream trips on shared arcs.
To compute the maximum feasible forward shift for a trip, we examine each arc along its path and determine how much earlier it can depart without introducing new conflicts. For each arc~$a$, we define $R'_a$ as the set of trips $r' \in \SortedTripsOnArc{a}$ such that $\Departure{a}{r'} < \Departure{a}{\Trip}$. If $R'_a \neq \emptyset$, we compute the shift to apply to $\Trip$'s departure time to match it with the latest departure or arrival time among these preceding trips as follows:
\[
\Delta s(a) = 
\Departure{a}{\Trip} - 
\max_{r' \in R'_a} 
\left\{ 
\Departure{a}{r'},\ 
\Arrival{a}{r'} 
\right\}.
\]
We shift the trip's departure time forward by the smallest $\Delta s(a)$ across all arcs, as long as the resulting time respects feasibility constraints. Formally:
\[
\StartTime{\Trip} \gets \max\left(
\TripEarliestDeparture{\Trip},\ 
\StartTime{\Trip} - 
\min_{a \in \TripPath{\Trip}} \Delta s(a)
\right).
\]
Once the procedure converges on the currently selected route, the algorithm resets $\Trip$'s departure time to $\TripEarliestDeparture{\Trip}$ and proceeds to evaluate the next available route. For each alternative route, it applies the same start time selection strategy. After exploring all route options, the algorithm selects the route--departure time combination that yields the lowest overall cost.

To conclude this section, we note that by slightly modifying the logic of the \textit{routing} operator, we can disable either the departure time assignment or the route assignment, resulting in variants that solve only the staggering or balancing problem, respectively. In our case study, we use these modes to isolate and compare the individual contribution of each component to the algorithm’s overall performance.

%% file: figures/methodology/operators_flowchart/operators_flowchart.pdf_tex
\begingroup%
  \makeatletter%
  \providecommand\color[2][]{%
    \errmessage{(Inkscape) Color is used for the text in Inkscape, but the package 'color.sty' is not loaded}%
    \renewcommand\color[2][]{}%
  }%
  \providecommand\transparent[1]{%
    \errmessage{(Inkscape) Transparency is used (non-zero) for the text in Inkscape, but the package 'transparent.sty' is not loaded}%
    \renewcommand\transparent[1]{}%
  }%
  \providecommand\rotatebox[2]{#2}%
  \newcommand*\fsize{\dimexpr\f@size pt\relax}%
  \newcommand*\lineheight[1]{\fontsize{\fsize}{#1\fsize}\selectfont}%
  \ifx\svgwidth\undefined%
    \setlength{\unitlength}{488.56616403bp}%
    \ifx\svgscale\undefined%
      \relax%
    \else%
      \setlength{\unitlength}{\unitlength * \real{\svgscale}}%
    \fi%
  \else%
    \setlength{\unitlength}{\svgwidth}%
  \fi%
  \global\let\svgwidth\undefined%
  \global\let\svgscale\undefined%
  \makeatother%
  \begin{picture}(1,0.15071217)%
    \lineheight{1}%
    \setlength\tabcolsep{0pt}%
    \put(0.36549768,0.1268111){\makebox(0,0)[lt]{\lineheight{1.25}\smash{\begin{tabular}[t]{l}\texttt{insert}\end{tabular}}}}%
    \put(0.38907593,0.06947952){\makebox(0,0)[lt]{\lineheight{1.25}\smash{\begin{tabular}[t]{l}\texttt{LS}\end{tabular}}}}%
    \put(0,0){\includegraphics[width=\unitlength,page=1]{operators_flowchart.pdf}}%
    \put(0.36527582,0.01214793){\makebox(0,0)[lt]{\lineheight{1.25}\smash{\begin{tabular}[t]{l}\texttt{remove}\end{tabular}}}}%
    \put(0.62367175,0.06947952){\makebox(0,0)[lt]{\lineheight{1.25}\smash{\begin{tabular}[t]{l}\texttt{routing}\end{tabular}}}}%
    \put(0,0){\includegraphics[width=\unitlength,page=2]{operators_flowchart.pdf}}%
    \put(0.86101412,0.06947952){\makebox(0,0)[lt]{\lineheight{1.25}\smash{\begin{tabular}[t]{l}\texttt{evaluation}\end{tabular}}}}%
    \put(0,0){\includegraphics[width=\unitlength,page=3]{operators_flowchart.pdf}}%
    \put(0.09919609,0.1268111){\makebox(0,0)[lt]{\lineheight{1.25}\smash{\begin{tabular}[t]{l}\texttt{RDUO}\end{tabular}}}}%
    \put(0.033,0.06947952){\makebox(0,0)[lt]{\lineheight{1.25}\smash{\begin{tabular}[t]{l}\texttt{greedyAssignment}\end{tabular}}}}%
    \put(0,0){\includegraphics[width=\unitlength,page=4]{operators_flowchart.pdf}}%
    \put(0.10787301,0.01214793){\makebox(0,0)[lt]{\lineheight{1.25}\smash{\begin{tabular}[t]{l}\texttt{LNS}\end{tabular}}}}%
    \put(0,0){\includegraphics[width=\unitlength,page=5]{operators_flowchart.pdf}}%
  \end{picture}%
\endgroup%

%% file: figures/methodology/index_in_set/index_in_set.pdf_tex
\begingroup%
  \makeatletter%
  \providecommand\color[2][]{%
    \errmessage{(Inkscape) Color is used for the text in Inkscape, but the package 'color.sty' is not loaded}%
    \renewcommand\color[2][]{}%
  }%
  \providecommand\transparent[1]{%
    \errmessage{(Inkscape) Transparency is used (non-zero) for the text in Inkscape, but the package 'transparent.sty' is not loaded}%
    \renewcommand\transparent[1]{}%
  }%
  \providecommand\rotatebox[2]{#2}%
  \newcommand*\fsize{\dimexpr\f@size pt\relax}%
  \newcommand*\lineheight[1]{\fontsize{\fsize}{#1\fsize}\selectfont}%
  \ifx\svgwidth\undefined%
    \setlength{\unitlength}{377.69611989bp}%
    \ifx\svgscale\undefined%
      \relax%
    \else%
      \setlength{\unitlength}{\unitlength * \real{\svgscale}}%
    \fi%
  \else%
    \setlength{\unitlength}{\svgwidth}%
  \fi%
  \global\let\svgwidth\undefined%
  \global\let\svgscale\undefined%
  \makeatother%
  \begin{picture}(1,0.13456773)%
    \lineheight{1}%
    \setlength\tabcolsep{0pt}%
    \put(0,0){\includegraphics[width=\unitlength,page=1]{index_in_set.pdf}}%
    \put(0.03165756,0.0611599){\makebox(0,0)[lt]{\lineheight{1.25}\smash{\begin{tabular}[t]{l}$r_1$\end{tabular}}}}%
    \put(0.12946221,0.0611599){\makebox(0,0)[lt]{\lineheight{1.25}\smash{\begin{tabular}[t]{l}$r_2$\end{tabular}}}}%
    \put(0.22726683,0.0611599){\color[rgb]{0.10196078,0.10196078,0.10196078}\makebox(0,0)[lt]{\lineheight{1.25}\smash{\begin{tabular}[t]{l}$r_3$\end{tabular}}}}%
    \put(0.42287607,0.0611599){\makebox(0,0)[lt]{\lineheight{1.25}\smash{\begin{tabular}[t]{l}$r_5$\end{tabular}}}}%
    \put(0.32507145,0.0611599){\makebox(0,0)[lt]{\lineheight{1.25}\smash{\begin{tabular}[t]{l}$r_4$\end{tabular}}}}%
    \put(0.52068069,0.0611599){\makebox(0,0)[lt]{\lineheight{1.25}\smash{\begin{tabular}[t]{l}$r_6$\end{tabular}}}}%
    \put(0.61848537,0.0611599){\makebox(0,0)[lt]{\lineheight{1.25}\smash{\begin{tabular}[t]{l}$r_7$\end{tabular}}}}%
    \put(0.71628999,0.0611599){\makebox(0,0)[lt]{\lineheight{1.25}\smash{\begin{tabular}[t]{l}$r_8$\end{tabular}}}}%
    \put(0.81409455,0.0611599){\makebox(0,0)[lt]{\lineheight{1.25}\smash{\begin{tabular}[t]{l}$r_9$\end{tabular}}}}%
    \put(0.91189923,0.0611599){\makebox(0,0)[lt]{\lineheight{1.25}\smash{\begin{tabular}[t]{l}$r_{10}$\end{tabular}}}}%
    \put(0,0){\includegraphics[width=\unitlength,page=2]{index_in_set.pdf}}%
    \put(0.18940035,0.00685694){\makebox(0,0)[lt]{\lineheight{1.25}\smash{\begin{tabular}[t]{l}$i_{\text{old}}$\end{tabular}}}}%
    \put(0,0){\includegraphics[width=\unitlength,page=3]{index_in_set.pdf}}%
    \put(0.38369235,0.00685694){\makebox(0,0)[lt]{\lineheight{1.25}\smash{\begin{tabular}[t]{l}$i_{\text{new}}$\end{tabular}}}}%
    \put(0,0){\includegraphics[width=\unitlength,page=4]{index_in_set.pdf}}%
    \put(0.57770939,0.11176296){\makebox(0,0)[lt]{\lineheight{1.25}\smash{\begin{tabular}[t]{l}$j_{\text{old}}$\\\\\end{tabular}}}}%
    \put(0,0){\includegraphics[width=\unitlength,page=5]{index_in_set.pdf}}%
    \put(0.77471367,0.11176296){\makebox(0,0)[lt]{\lineheight{1.25}\smash{\begin{tabular}[t]{l}$j_{\text{new}}$\\\end{tabular}}}}%
    \put(0,0){\includegraphics[width=\unitlength,page=6]{index_in_set.pdf}}%
  \end{picture}%
\endgroup%

%% file: figures/methodology/ls_move/move_a.pdf_tex
\begingroup%
  \makeatletter%
  \providecommand\color[2][]{%
    \errmessage{(Inkscape) Color is used for the text in Inkscape, but the package 'color.sty' is not loaded}%
    \renewcommand\color[2][]{}%
  }%
  \providecommand\transparent[1]{%
    \errmessage{(Inkscape) Transparency is used (non-zero) for the text in Inkscape, but the package 'transparent.sty' is not loaded}%
    \renewcommand\transparent[1]{}%
  }%
  \providecommand\rotatebox[2]{#2}%
  \newcommand*\fsize{\dimexpr\f@size pt\relax}%
  \newcommand*\lineheight[1]{\fontsize{\fsize}{#1\fsize}\selectfont}%
  \ifx\svgwidth\undefined%
    \setlength{\unitlength}{116.704617bp}%
    \ifx\svgscale\undefined%
      \relax%
    \else%
      \setlength{\unitlength}{\unitlength * \real{\svgscale}}%
    \fi%
  \else%
    \setlength{\unitlength}{\svgwidth}%
  \fi%
  \global\let\svgwidth\undefined%
  \global\let\svgscale\undefined%
  \makeatother%
  \begin{picture}(1,0.8497888)%
    \lineheight{1}%
    \setlength\tabcolsep{0pt}%
    \put(0,0){\includegraphics[width=\unitlength,page=1]{move_a.pdf}}%
    \put(0.46249548,0.77598475){\makebox(0,0)[lt]{\lineheight{1.25}\smash{\begin{tabular}[t]{l}$\tau_a$\end{tabular}}}}%
    \put(0.4612905,0.02219141){\makebox(0,0)[lt]{\lineheight{1.25}\smash{\begin{tabular}[t]{l}$\tau_b$\end{tabular}}}}%
  \end{picture}%
\endgroup%

%% file: figures/methodology/ls_move/move_b.pdf_tex
\begingroup%
  \makeatletter%
  \providecommand\color[2][]{%
    \errmessage{(Inkscape) Color is used for the text in Inkscape, but the package 'color.sty' is not loaded}%
    \renewcommand\color[2][]{}%
  }%
  \providecommand\transparent[1]{%
    \errmessage{(Inkscape) Transparency is used (non-zero) for the text in Inkscape, but the package 'transparent.sty' is not loaded}%
    \renewcommand\transparent[1]{}%
  }%
  \providecommand\rotatebox[2]{#2}%
  \newcommand*\fsize{\dimexpr\f@size pt\relax}%
  \newcommand*\lineheight[1]{\fontsize{\fsize}{#1\fsize}\selectfont}%
  \ifx\svgwidth\undefined%
    \setlength{\unitlength}{121.51357936bp}%
    \ifx\svgscale\undefined%
      \relax%
    \else%
      \setlength{\unitlength}{\unitlength * \real{\svgscale}}%
    \fi%
  \else%
    \setlength{\unitlength}{\svgwidth}%
  \fi%
  \global\let\svgwidth\undefined%
  \global\let\svgscale\undefined%
  \makeatother%
  \begin{picture}(1,0.73797985)%
    \lineheight{1}%
    \setlength\tabcolsep{0pt}%
    \put(0,0){\includegraphics[width=\unitlength,page=1]{move_b.pdf}}%
    \put(-0.00520775,0.60021417){\makebox(0,0)[lt]{\lineheight{1.25}\smash{\begin{tabular}[t]{l}$r$\end{tabular}}}}%
    \put(-0.00520775,0.24161143){\makebox(0,0)[lt]{\lineheight{1.25}\smash{\begin{tabular}[t]{l}$r'$\end{tabular}}}}%
    \put(0,0){\includegraphics[width=\unitlength,page=2]{move_b.pdf}}%
    \put(0.34393751,0.0213132){\makebox(0,0)[lt]{\lineheight{1.25}\smash{\begin{tabular}[t]{l}$e^r$\end{tabular}}}}%
    \put(0.51727418,0.60426462){\makebox(0,0)[lt]{\lineheight{1.25}\smash{\begin{tabular}[t]{l}$\tau_a$\end{tabular}}}}%
    \put(0.76362612,0.60426462){\makebox(0,0)[lt]{\lineheight{1.25}\smash{\begin{tabular}[t]{l}$d_a$\end{tabular}}}}%
    \put(0,0){\includegraphics[width=\unitlength,page=3]{move_b.pdf}}%
    \put(0.10696942,0.0213132){\makebox(0,0)[lt]{\lineheight{1.25}\smash{\begin{tabular}[t]{l}$e^{r'}$\end{tabular}}}}%
    \put(0,0){\includegraphics[width=\unitlength,page=4]{move_b.pdf}}%
    \put(0.26843324,0.2563668){\makebox(0,0)[lt]{\lineheight{1.25}\smash{\begin{tabular}[t]{l}$\tau_a$\end{tabular}}}}%
    \put(0,0){\includegraphics[width=\unitlength,page=5]{move_b.pdf}}%
  \end{picture}%
\endgroup%

%% file: figures/methodology/ls_move/move_c.pdf_tex
\begingroup%
  \makeatletter%
  \providecommand\color[2][]{%
    \errmessage{(Inkscape) Color is used for the text in Inkscape, but the package 'color.sty' is not loaded}%
    \renewcommand\color[2][]{}%
  }%
  \providecommand\transparent[1]{%
    \errmessage{(Inkscape) Transparency is used (non-zero) for the text in Inkscape, but the package 'transparent.sty' is not loaded}%
    \renewcommand\transparent[1]{}%
  }%
  \providecommand\rotatebox[2]{#2}%
  \newcommand*\fsize{\dimexpr\f@size pt\relax}%
  \newcommand*\lineheight[1]{\fontsize{\fsize}{#1\fsize}\selectfont}%
  \ifx\svgwidth\undefined%
    \setlength{\unitlength}{121.33241308bp}%
    \ifx\svgscale\undefined%
      \relax%
    \else%
      \setlength{\unitlength}{\unitlength * \real{\svgscale}}%
    \fi%
  \else%
    \setlength{\unitlength}{\svgwidth}%
  \fi%
  \global\let\svgwidth\undefined%
  \global\let\svgscale\undefined%
  \makeatother%
  \begin{picture}(1,0.73548072)%
    \lineheight{1}%
    \setlength\tabcolsep{0pt}%
    \put(-0.00175082,0.60206335){\makebox(0,0)[lt]{\lineheight{1.25}\smash{\begin{tabular}[t]{l}$r$\end{tabular}}}}%
    \put(-0.00175082,0.24292517){\makebox(0,0)[lt]{\lineheight{1.25}\smash{\begin{tabular}[t]{l}$r'$\end{tabular}}}}%
    \put(0,0){\includegraphics[width=\unitlength,page=1]{move_c.pdf}}%
    \put(0.33204271,0.02134503){\makebox(0,0)[lt]{\lineheight{1.25}\smash{\begin{tabular}[t]{l}$\rightarrow s^r$\end{tabular}}}}%
    \put(0.62636151,0.60611985){\makebox(0,0)[lt]{\lineheight{1.25}\smash{\begin{tabular}[t]{l}$\tau_a$\end{tabular}}}}%
    \put(0,0){\includegraphics[width=\unitlength,page=2]{move_c.pdf}}%
    \put(0.11059384,0.02134503){\makebox(0,0)[lt]{\lineheight{1.25}\smash{\begin{tabular}[t]{l}$e^{r'}$\end{tabular}}}}%
    \put(0,0){\includegraphics[width=\unitlength,page=3]{move_c.pdf}}%
    \put(0.27229875,0.25770257){\makebox(0,0)[lt]{\lineheight{1.25}\smash{\begin{tabular}[t]{l}$\tau_a$\end{tabular}}}}%
    \put(0,0){\includegraphics[width=\unitlength,page=4]{move_c.pdf}}%
  \end{picture}%
\endgroup%

%% file: figures/methodology/ls_move/move_d.pdf_tex
\begingroup%
  \makeatletter%
  \providecommand\color[2][]{%
    \errmessage{(Inkscape) Color is used for the text in Inkscape, but the package 'color.sty' is not loaded}%
    \renewcommand\color[2][]{}%
  }%
  \providecommand\transparent[1]{%
    \errmessage{(Inkscape) Transparency is used (non-zero) for the text in Inkscape, but the package 'transparent.sty' is not loaded}%
    \renewcommand\transparent[1]{}%
  }%
  \providecommand\rotatebox[2]{#2}%
  \newcommand*\fsize{\dimexpr\f@size pt\relax}%
  \newcommand*\lineheight[1]{\fontsize{\fsize}{#1\fsize}\selectfont}%
  \ifx\svgwidth\undefined%
    \setlength{\unitlength}{121.51357936bp}%
    \ifx\svgscale\undefined%
      \relax%
    \else%
      \setlength{\unitlength}{\unitlength * \real{\svgscale}}%
    \fi%
  \else%
    \setlength{\unitlength}{\svgwidth}%
  \fi%
  \global\let\svgwidth\undefined%
  \global\let\svgscale\undefined%
  \makeatother%
  \begin{picture}(1,0.73797985)%
    \lineheight{1}%
    \setlength\tabcolsep{0pt}%
    \put(0,0){\includegraphics[width=\unitlength,page=1]{move_d.pdf}}%
    \put(-0.00520782,0.60021417){\makebox(0,0)[lt]{\lineheight{1.25}\smash{\begin{tabular}[t]{l}$r$\end{tabular}}}}%
    \put(-0.00520782,0.24161143){\makebox(0,0)[lt]{\lineheight{1.25}\smash{\begin{tabular}[t]{l}$r'$\end{tabular}}}}%
    \put(0,0){\includegraphics[width=\unitlength,page=2]{move_d.pdf}}%
    \put(0.34393744,0.0213132){\makebox(0,0)[lt]{\lineheight{1.25}\smash{\begin{tabular}[t]{l}$e^r$\end{tabular}}}}%
    \put(0.51727411,0.60426462){\makebox(0,0)[lt]{\lineheight{1.25}\smash{\begin{tabular}[t]{l}$\tau_{a \rightarrow b}$\end{tabular}}}}%
    \put(0,0){\includegraphics[width=\unitlength,page=3]{move_d.pdf}}%
    \put(0.10696935,0.0213132){\makebox(0,0)[lt]{\lineheight{1.25}\smash{\begin{tabular}[t]{l}$e^{r'}$\end{tabular}}}}%
    \put(0,0){\includegraphics[width=\unitlength,page=4]{move_d.pdf}}%
    \put(0.26843317,0.2563668){\makebox(0,0)[lt]{\lineheight{1.25}\smash{\begin{tabular}[t]{l}$\tau_a$\end{tabular}}}}%
    \put(0,0){\includegraphics[width=\unitlength,page=5]{move_d.pdf}}%
  \end{picture}%
\endgroup%

%% file: contents/4.CaseStudy.tex
\section{Design of experiments}\label{sec:case_study}
In this section, we present the setup for our case study, including the preparation of the street network, the parameterization of the travel time function, and the trip design. We conclude with a description of the algorithm specifics.
 
\subparagraph{Network}
We obtain the street network for Manhattan, New York City, from OpenStreetMap \citep{OpenStreetMap}. As part of preprocessing, we remove motorway roads, retain only the shortest arc among parallel connections between the same node pair, and consolidate intersections located within a 20-meter radius. Rather than using the entire street network, we extract the southernmost 10\% of nodes and their connecting arcs, then retain the largest strongly connected component. The resulting network consists of 302 nodes and 734 arcs. We compute the nominal travel time on each arc by assuming a constant speed of 20 kilometers per hour.
\subparagraph{Trip design} 
To model realistic traffic demand, we estimate the number of trips within our study area based on empirical data. During the evening peak on a typical fall weekday, approximately 70000 motorized vehicles  — including both base load and taxi demand — enter and leave Manhattan’s Central Business District \citep{NYMTC2023}. Assuming this approximates the total number of trips within Manhattan, and given that our study focuses on a subnetwork covering roughly 10\% of the area, we scale this value accordingly and consider 7000 trips per instance to be a reasonable estimate.

We consider 31 instances, one for each day of a month. For each instance, we construct the trip set $\SetTrips$ using the 2009 NYC Taxi \& Limousine Commission dataset \citep{NTLC2009}, which provides origin, destination, and departure time information for individual taxi trips. We include all trips with both origin and destination located within the subnetwork over the full 24-hour period of each day, using data from the months of January, March, May, and July (i.e., all 31-day months). To simulate realistic congestion, we compress the departure times into a one-hour window while preserving the original minute and second resolution. This procedure yields instances with an average of approximately 6000 trips, which aligns well with our target scale.

The design of alternative routes is central to the characterization of trip data in our case study. It affects both traffic balancing in our algorithmic solutions and congestion levels observed in the \gls{rduo} baseline, as trips in both settings select routes from the same predefined set. 
We solve the \gls{kspwlo} problem with a limited overlap setting the number of alternative paths $k$ to five and the similarity threshold $\SimilarityThs$ to 60\%. This choice aligns with the goal of ensuring good network coverage through alternative routes, both to reproduce realistic congestion levels in \gls{rduo} and to facilitate traffic balancing, while keeping the number of alternatives small enough to maintain a manageable search space. We discuss this design choice in further detail in Appendix~\ref{appendix:route_design}.

To complete the trip's characterization, we define each trip’s latest arrival time $\TripLatestArrival{\Trip}$ based on the travel time observed in the \gls{rduo} solution, extended by 25\% to account for congestion variability. We also set the maximum staggering time $\TripMaxStaggering{\Trip}$ to 20\% of the shortest-route free-flow travel time. We provide a detailed overview of all instances in Appendix~\ref{appendix:instances_description}.

\subparagraph{Delay function parameterization}
We characterize the delay function $\DelayFunction$ as a polynomial expression of the form:
\[
\DelayFunction = \NominalTravelTime{\Arc} \cdot \alpha \left[\left(\frac{\Flow{\Arc}{\Trip} + \beta}{\NominalTravelTime{\Arc}} \right)^\gamma - \left(\frac{\beta}{\NominalTravelTime{\Arc}} \right)^\gamma \right],
\]
where $\alpha$ controls the overall scaling of the delay, $\beta$ introduces a baseline shift to regulate delays under low congestion, and $\gamma$ is the polynomial exponent. The subtraction term ensures that the delay function evaluates to zero in the absence of other trips.

To parameterize this function, we used data from the TLC dataset, which includes pickup times, drop-off times, and trip distances. Assuming a constant free-flow speed of 20 km/h, we computed the free-flow travel time and compared it with the trip-recorded durations. In the same network area considered for our experiments, the data indicate a median delay of approximately one and a half minutes during peak hours. Furthermore, commercial mapping services report that travel between the Financial District and Lower Manhattan experiences a delay of approximately 25\% relative to the free-flow travel time during the peak hour.

Given these observations, we set the parameters of the delay function to $\alpha = 0.1$, $\beta = 35$, and $\gamma = 3$. This configuration yields instances with a \gls{rduo} characterized by an average total delay representing 12\% of the total travel time, and individual trip delays with a median of approximately 0.7 minutes. Although these figures indicate that our instances may underestimate real-world congestion, they effectively serve to demonstrate the core concepts of our algorithm. 

\subparagraph{Algorithm specifics}  
We configure the \gls{lns} algorithm with the following parameters: an initial pool size of $x = 40\%$, a sample size of $y = 10\%$ drawn from this pool, and a maximum of $z = 2$ destroy-repair cycles before invoking the local search for intensification. These values reflect the outcome of a structured parametric analysis, which we report in detail in Appendix~\ref{appendix:parameters_search}.

\section{Results}\label{sec:results}
In this section, we present the results of our numerical experiments. 
As a first step, we compare the staggering-only variant of our algorithm with the matheuristic proposed in \citet{CoppolaHiermannEtAl2024} to benchmark performance and demonstrate that the new algorithm consistently outperforms it across different levels of flexibility in trip start times.

We then move to the analysis of an idealized scenario in which we assume full control over all vehicles in the system. This setup allows us to establish an upper bound on the performance achievable through centralized routing. Within this framework, we isolate the individual effects of staggering and balancing on fleet and network performance and compare them to the performance of the fully integrated algorithm.

We then examine how to apply staggering and balancing in practice, assess their impact on individual travel and arrival times, and analyze how centralized control reshapes the spatial distribution of congestion. To complement this analysis, we consider more realistic scenarios in which only a portion of the traffic is under centralized control. We explore this from two perspectives: first, a municipality aiming to minimize total system-wide travel time; second, a profit-oriented \gls{amod} service provider focused solely on minimizing the travel time of its own fleet.

We implemented our metaheuristic in C\texttt{++} and compiled it using GCC version 11.1.0 with \texttt{-O3} optimization flags. We ran all experiments on an Intel(R) i9-9900 CPU at 3.1 GHz with 10 GB of RAM. {The full implementation and experimental setup are available at \url{https://github.com/tumBAIS/integ_bal_stag}}.

\subsection{Comparison with matheuristic baseline}
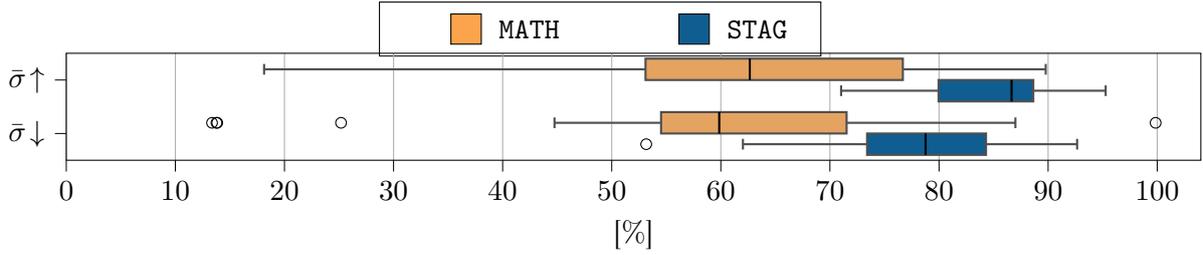
\begin{figure}[!t]
    \centering
    \subfigcapskip = -10pt
    \begin{minipage}{\textwidth}
        \centering
        \import{figures/experiments/math_comparison/}{legend.tex}

        \vspace{-0.5em} 
        \import{figures/experiments/math_comparison/}{relative_total_delay_reduction_congestion_HC.tex}
    \end{minipage}
    \caption{\raggedright Relative delay reductions achieved by \texttt{STAG} and \texttt{MATH} under reduced (\( \SigmaDown \)) and increased (\( \SigmaUp \)) levels of flexibility in trip start times.}
    \label{fig:math_comparison}
\end{figure}
Figure~\ref{fig:math_comparison} compares the performance of our staggering-only algorithm, \texttt{STAG}, with the matheuristic \texttt{MATH} from the literature, using relative delay reduction as the performance metric. {We evaluate both methods on a set of 31 instances under two levels of departure-time flexibility: $\SigmaDown$, where the maximum allowable start time shift $\bar{\sigma}_r$ is set to 10\% of the nominal travel time along each trip’s assigned route, and $\SigmaUp$, where this shift is increased to 20\%.}
For a detailed description of the experimental instances, we refer to Appendix~\ref{appendix:comparison_with_math}.

\texttt{STAG} outperforms \texttt{MATH} in terms of median delay reduction by approximately 20 and 25 percentage points under $\SigmaDown$ and $\SigmaUp$, respectively.
Interestingly, increasing flexibility from $\SigmaDown$ to $\SigmaUp$ brings only modest gains for \texttt{MATH} (around 3 percentage points in the median) and coincides with the loss of the sole optimal solution found. This highlights how \texttt{MATH} benefits from tightly bounded departure times, which enable leaner \gls{milp} models with fewer {${\text{big-M}}$} constraints. However, as flexibility grows, these structural advantages vanish: the search space becomes more symmetric, relaxations weaken, and performance declines.
{In contrast, \texttt{STAG} benefits from increased flexibility: its computational complexity and efficiency remain stable even with wider time windows, allowing it to improve median delay reduction by around ten percentage points and enhance robustness through stronger worst-case performance and reduced variance.}

\begin{result}
{
Across different levels of flexibility in trip start times, \texttt{STAG} consistently outperforms \texttt{MATH}, delivering higher delay reductions with reduced variance.}
\end{result}

Notably, the instances considered in this comparison already push \texttt{MATH} to its computational limits, since tracking pairwise trip interactions requires an exponential number of constraints, which rapidly increase memory usage.
In contrast, \texttt{STAG} circumvents this exponential growth entirely, making it well-suited for large-scale applications.
The combination of superior performance and scalability makes \texttt{STAG} the natural choice for solving the staggering-only variant of the problem in the remainder of this study, where the number of pairwise trip interactions is significantly higher than in the present benchmark.

\subsection{Algorithmic performance}
\begin{figure}[!t]
\centering
\subfigcapskip = -10pt
\begin{minipage}{\textwidth}
    \centering
    \import{figures/experiments/basic_setting_performance/}{legend.tex}
    \vspace{0.2cm}

    \subfigure[Relative delay variations\label{fig:performance_a}]{
        \import{figures/experiments/basic_setting_performance/}{relative_delay_reduction_boxplot.tex}
    }
    \hfill
    \subfigure[Absolute delay variations\label{fig:performance_b}]{
        \import{figures/experiments/basic_setting_performance/}{absolute_delay_reduction_boxplot.tex}
    }
\end{minipage}
\caption{Figure~\ref{fig:performance_a} shows the relative reductions in total and congestion delay compared to \texttt{RDUO} under full traffic control. Figure~\ref{fig:performance_b} shows the corresponding absolute values, including detour delay increases induced by balancing.}
\label{fig:performance}
\end{figure}
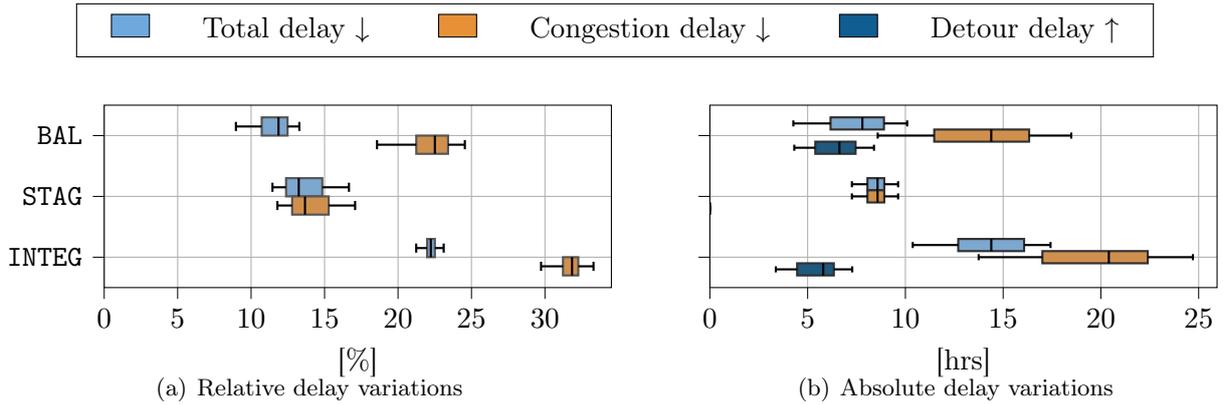

{Figure~\ref{fig:performance} illustrates the performance of three variants of our algorithm on the instances described in Section~\ref{sec:case_study}}: \texttt{BAL}, which applies only balancing operations; \texttt{STAG}, which only staggers trip departures; and \texttt{INTEG}, which combines both strategies. The figure reports total and congestion-induced delay reductions relative to \gls{rduo} under full traffic control, expressed in both percentage and absolute terms, along with the absolute detour delay introduced by balancing. All experiments terminated within a four-hour time limit.
\subparagraph{Total delay reductions.}
When isolating individual mechanisms, we observe that balancing alone (\texttt{BAL}) yields a median delay reduction of 12\% (approximately 8 hours). This result suggests that self-interested users, unaware of future congestion, often select suboptimal routes, while the balancing operator distributes traffic more efficiently across available paths. Staggering alone (\texttt{STAG}) achieves a comparable median delay reduction of 13\% (approximately 9 hours), demonstrating that our greedy staggering strategy, which adjusts departure times to resolve arc-level conflicts without accounting for upstream congestion effects, serves as an effective heuristic proxy for identifying promising trip start times.

The fully integrated algorithm (\texttt{INTEG}) consistently delivers strong performance across all real-world instances, achieving total delay reductions consistently around 23\%. In absolute terms, this corresponds to 10 to 17 hours of reduced travel time.
Crucially, the integrated approach appears to combine the benefits of both rerouting and rescheduling in a near-additive manner, highlighting their complementary roles in mitigating delay.

\begin{result}
The integrated algorithm achieves total delay reductions of 23\%, translating into up to 17 hours of saved travel time across the entire system.
\end{result}

\subparagraph{Trade-offs between congestion and detour.}
From the operator’s perspective, total delay reduction translates directly into shorter travel times and economic savings. Yet, it is equally important to assess network efficiency. Specifically, we assess the reduction in congestion delay, potentially at the cost of increased detour delay.

In \texttt{STAG}, all delay reduction corresponds to congestion relief, since routes remain fixed. Interestingly, balancing alone well mitigates congestion: \texttt{BAL} removes over 23\% of congestion-induced delay in median, which translates to approximately 15 hours. However, this improvement comes with additional detour delay — about six hours on average. 
Remarkably, the integrated approach delivers an even more favorable trade-off. By leveraging staggering to reduce congestion without modifying routes, \texttt{INTEG} removes more congestion delay compared to balancing alone (up to 25 hours) while inducing less additional detour delays.

\begin{result}
The integrated algorithm removes more congestion delay (up to 25 hours) than balancing alone (up to 18 hours), while introducing less detour delays (under 6 hours).
\end{result}
\begin{figure}[!b]
\centering
\begin{minipage}{\textwidth}
    \centering
    \import{figures/experiments/solution_structure/}{legend.tex}
    \vspace{0.2cm}
    \subfigcapskip = -10pt

    \subfigure[Selected routes\label{fig:solution_structure_a}]{
        \import{figures/experiments/solution_structure/}{user_equilibrium_route_distribution.tex}
    }
    \subfigure[Staggering applied\label{fig:solution_structure_b}]{
        \import{figures/experiments/solution_structure/}{staggering_frequency_barplot.tex}
    }
    \subfigure[Travel time difference\label{fig:solution_structure_c}]{
        \import{figures/experiments/solution_structure/}{travel_time_diff_frequency_barplot.tex}
    }
    \subfigure[Arrival time difference\label{fig:solution_structure_d}]{
        \import{figures/experiments/solution_structure/}{arrival_diff_frequency_barplot.tex}
    }
\end{minipage}
\caption{Frequency plots comparing the structure of the baseline and optimized solutions: 
(a) route assignments, with route IDs sorted from shortest to longest distance; 
(b) staggering operations applied, where zero indicates no staggering; 
(c) travel time differences, where negative values indicate shorter travel times in the solution relative to the baseline, and positive values indicate increases; and 
(d) arrival time differences, where negative values indicate earlier arrivals and positive values indicate delays relative to the baseline.}
\label{fig:solution_structure}
\end{figure}
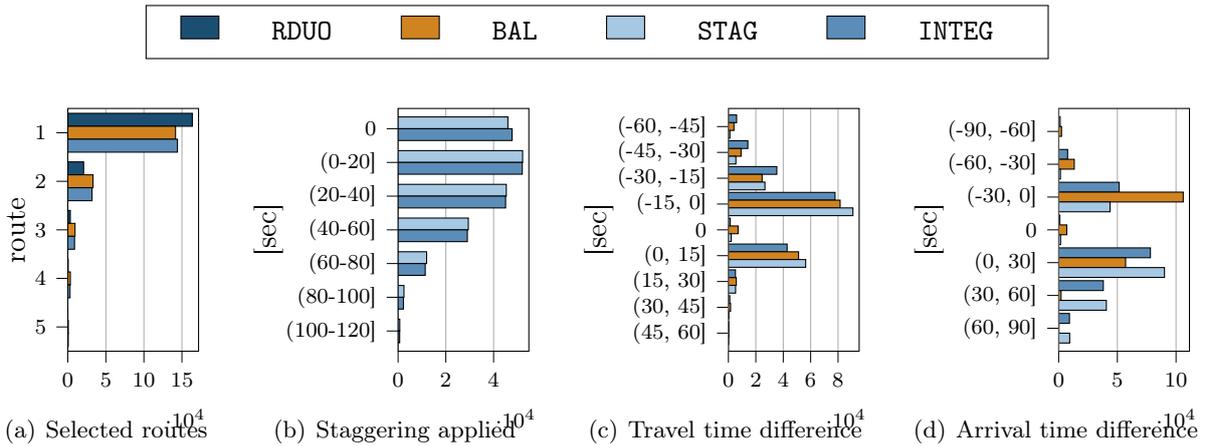
\subsection{Structural analysis}

In this section, we analyze the structural characteristics of the solutions generated by different control strategies. Figure~\ref{fig:solution_structure} presents: (a) the route assignment distribution for the baseline \gls{rduo} and for the rebalancing variants \texttt{BAL} and \texttt{INTEG}, (b) the frequency of staggering in \texttt{STAG} and \texttt{INTEG}, (c) travel time differences relative to the baseline, and (d) differences in arrival times.

We first observe that the \gls{rduo} configuration tends to overload the shortest available routes. In contrast, both \texttt{BAL} and \texttt{INTEG} redistribute a portion of the demand to longer distance alternatives to mitigate congestion. Expectedly, \texttt{INTEG} assigns a slightly larger fraction of trips to the shortest route than \texttt{BAL}: by also applying staggering, \texttt{INTEG} reduces reliance on detouring to achieve higher performance gains. 
Regarding staggering behavior, both \texttt{INTEG} and \texttt{STAG} most frequently apply small shifts, as shown by the higher concentration of trips in the zero to twenty-second departure shift bin.
Generally, \texttt{INTEG} applies staggering to a larger share of trips and with greater intensity than \texttt{STAG}, highlighting how the additional flexibility offered by route rebalancing unlocks further staggering opportunities. The resulting staggering distribution exhibits a clear decay: most adjustments fall within the first 40 seconds, with diminishing frequency thereafter.

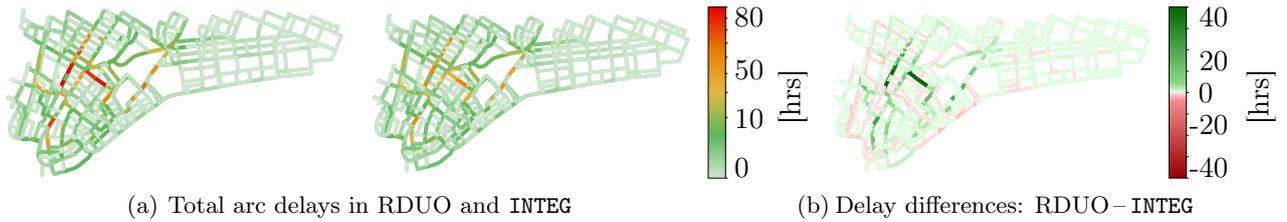
\begin{figure}[!t]
\centering
\newcommand{\colorbaroffset}{2.4}  
\newcommand{\colorbaroffsetsecond}{3.1}  
\hspace{-2cm}
\subfigure[Total arc delays in \gls{rduo} and \texttt{INTEG} \label{fig:heatmap_a}]{
    \begin{minipage}{0.3\textwidth}
        \includegraphics[width=\linewidth]{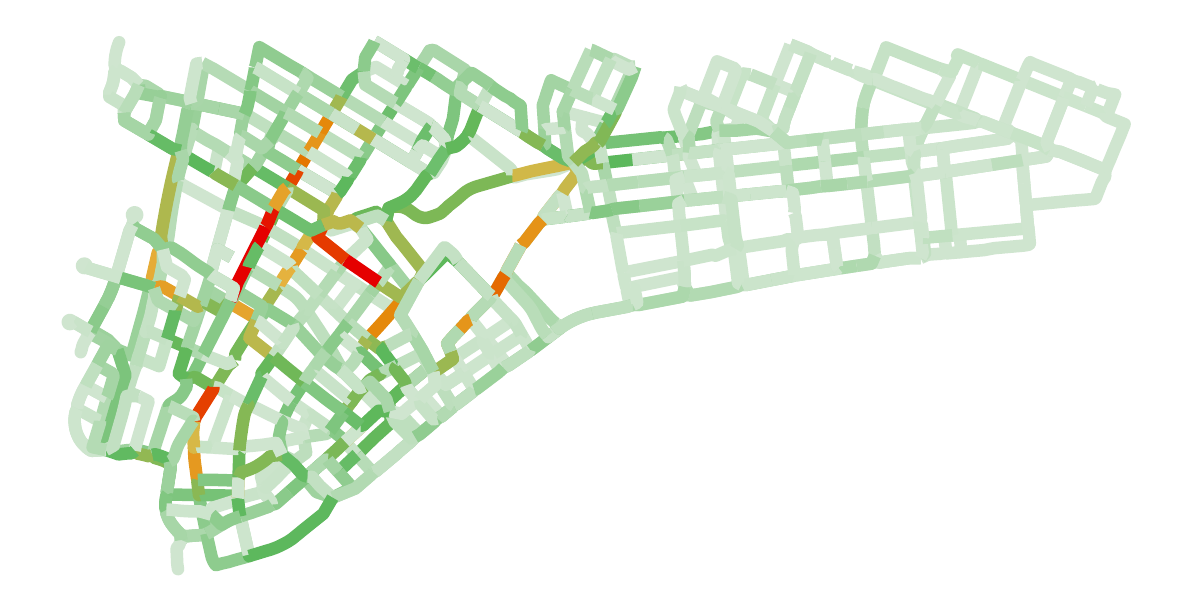}
    \end{minipage}
    \hspace{-.5cm}
    \begin{minipage}{0.3\textwidth}
        \begin{tikzpicture}
            \node at (0,0) {\includegraphics[width=\linewidth]{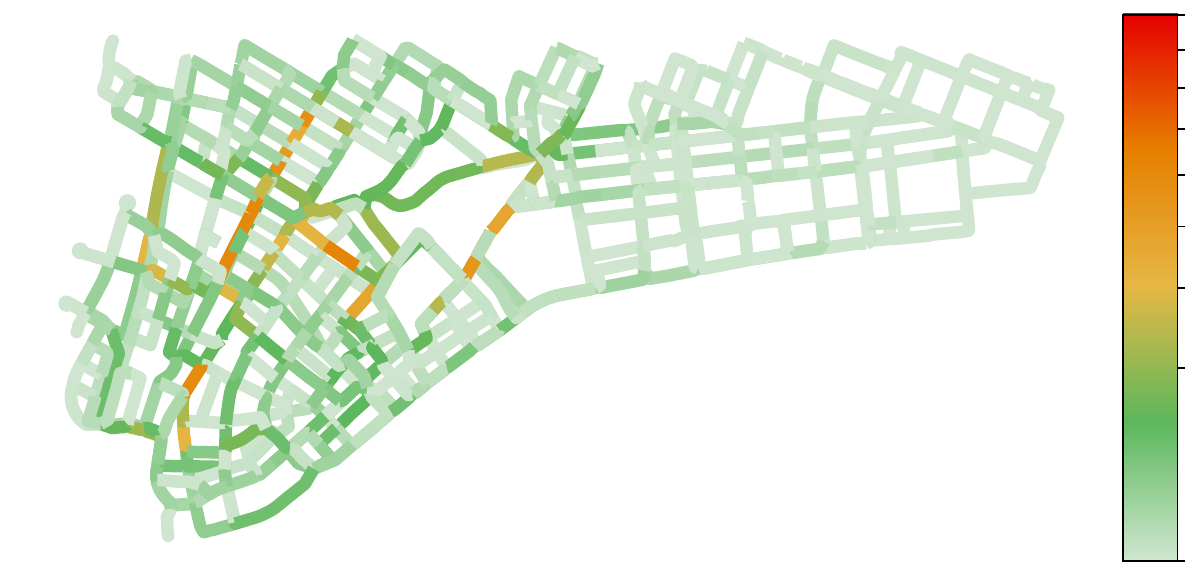}};
\node[anchor=west] at (\colorbaroffset, 1) {80};
\node[anchor=west] at (\colorbaroffset, 0.3) {50};
\node[anchor=west] at (\colorbaroffset, -0.35) {10};
\node[anchor=west] at (\colorbaroffset, -1) {0}; 
\node[anchor=west, rotate=90] at ({\colorbaroffset + .9}, -0.55) {[hrs]};
        \end{tikzpicture}
    \end{minipage}
}
\hspace{.6cm}
\subfigure[\!Delay differences: \gls{rduo}\! --\! \texttt{INTEG}\!\label{fig:heatmap_b}]{
    \begin{minipage}{0.3\textwidth}
        \begin{tikzpicture}
            \node at (0,0) {\includegraphics[width=\linewidth]{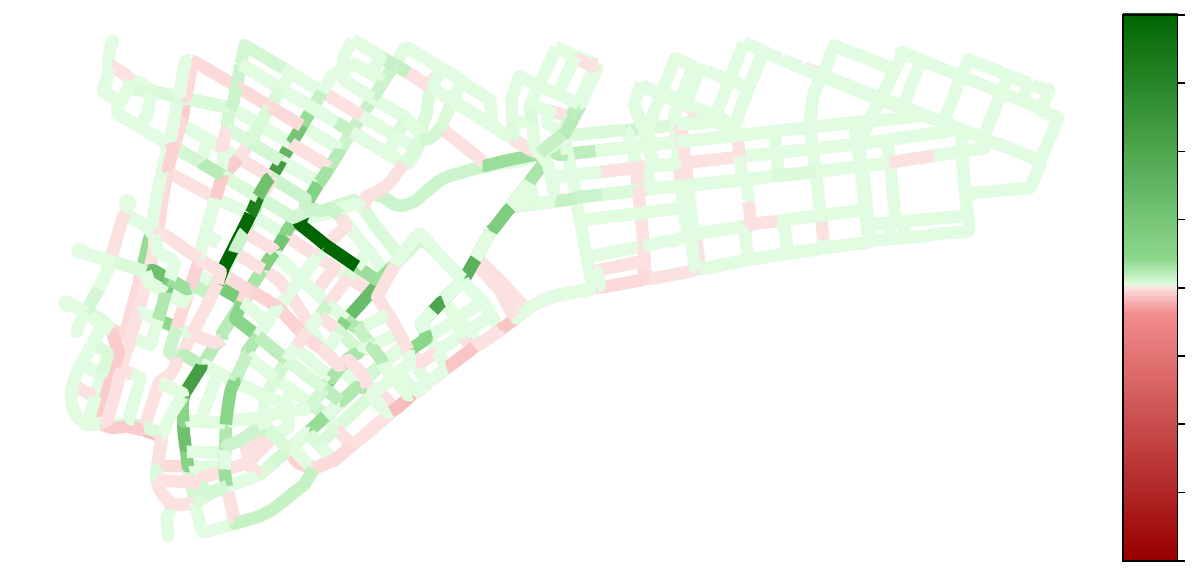}};
            \node[anchor=east] at (\colorbaroffsetsecond, 1) {40};
            \node[anchor=east] at (\colorbaroffsetsecond, 0.425) {20};
            \node[anchor=east] at (\colorbaroffsetsecond - .2, 0.00) {0};
            \node[anchor=east] at (\colorbaroffsetsecond, -0.425) {-20};
            \node[anchor=east] at (\colorbaroffsetsecond, -1) {-40};
            \node[anchor=west, rotate=90] at (\colorbaroffsetsecond + .3, -0.55) {[hrs]};
        \end{tikzpicture}
    \end{minipage}
}

\caption{Spatial distribution of arc-level delays over 31-day instances in \gls{rduo} and \texttt{INTEG} (a), and their delay differences (b).}
\label{fig:heatmap}
\end{figure}

The distribution of travel time differences reveals that achieving system-wide efficiency requires some trips to incur longer travel times. These increases are typically modest, mostly within 15 seconds, and only a small fraction exceeds 30 seconds. This trade-off is counterbalanced by a greater number of trips benefiting from travel time reductions, with improvements reaching up to 60 seconds.
Finally, arrival time differences are distributed almost symmetrically around zero, with only a small fraction of trips arriving at exactly the same time as in the \gls{rduo} configuration. In the \texttt{BAL} variant, the distribution is skewed toward earlier arrivals, whereas configurations involving staggering exhibit a slight skew toward positive differences, consistent with their design to delay trips strategically.
\begin{result}
Centralized coordination yields solutions that are structurally distinct from \gls{rduo}: (i) it alleviates congestion by rebalancing demand away from the shortest routes, (ii) it integrates balancing and staggering to identify larger pools of promising departure times, and (iii) it enhances overall efficiency by allowing slight travel time increases for some trips, while maintaining tight control over arrival times.
\end{result}
To complement the structural analysis and illustrate how our algorithm reshapes congestion patterns, Figure~\ref{fig:heatmap} compares cumulative arc-level delays across all 31-day instances between the \gls{rduo} and \texttt{INTEG} solutions.
In the \gls{rduo} configuration, congestion remains concentrated in the central region of the network. By contrast, the \texttt{INTEG} solution achieves a more balanced traffic distribution, effectively alleviating these bottlenecks, particularly in areas where alternative routes support rerouting. Importantly, this improvement does not cause delay displacement or introduce new congestion hotspots. As illustrated in the differential map, central arcs benefit from substantial delay reductions of up to 40 hours, while newly utilized arcs see only minor increases, typically below 5 hours.

\begin{result}
Our algorithm alleviates major bottlenecks, reducing delays by up to 40 hours, without shifting congestion elsewhere. Delay increases on newly utilized arcs remain minimal, generally under 5 hours, effectively preventing the emergence of new congestion points.
\end{result}

\subsection{Flow control analysis}

This section analyzes how total delay evolves as the fraction of controlled traffic $\FractionControlled$ increases.
In each instance, trips are randomly designated as \gls{amod} trips with probability $\FractionControlled$.
Baseload trips follow fixed behavior, consistent with the \gls{rduo} solution; we leave the study of responsive baseload dynamics to \gls{amod} behavior to future work.
We consider two control paradigms. In the first, a welfare-oriented authority aims to minimize total system delay by optimizing the travel times of all trips.
In the second, a self-interested operator controls only a fleet of \gls{amod} vehicles and minimizes the delay experienced by its own fleet, without accounting for delays among baseload trips.

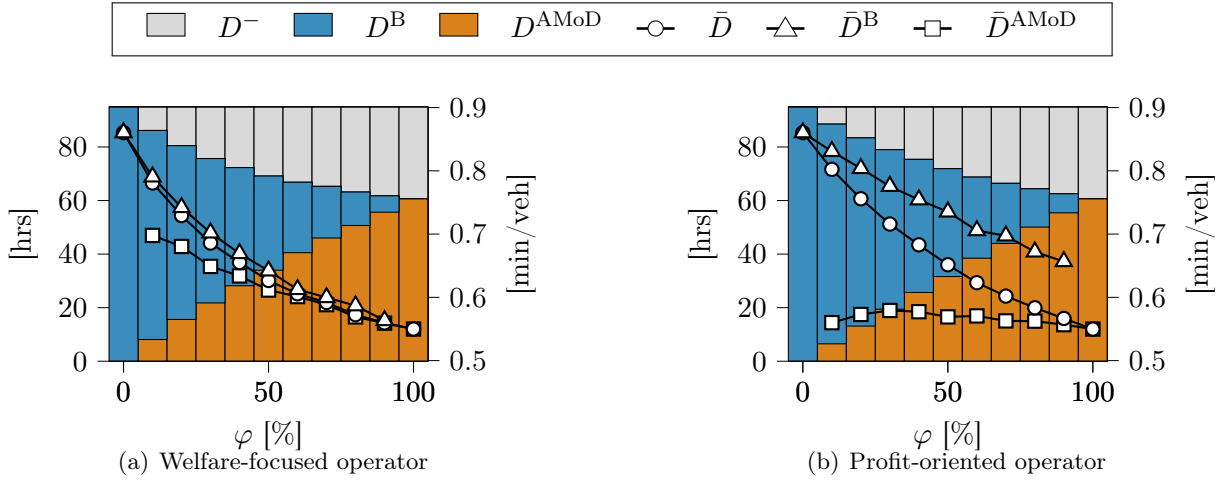
\begin{figure}[!t]
\centering
\begin{minipage}{\textwidth}
    \centering
    \subfigcapskip = -10pt
    \import{figures/experiments/flow_control/}{legend_plot.tex}

    \subfigure[Welfare-focused operator \label{fig:flow_control_a}]{
        \import{figures/experiments/flow_control/}{welfare_goal_plot.tex}
    }
    \hfill
    \subfigure[Profit-oriented operator\label{fig:flow_control_b}]{
        \import{figures/experiments/flow_control/}{profit_goal_plot.tex}
    }
\end{minipage}
\caption{Impact of the controlled traffic fraction \( \FractionControlled \) on total and average delay metrics. 
Figure~\ref{fig:flow_control_a} presents results for a welfare-focused operator aiming to minimize total system travel time. 
Figure~\ref{fig:flow_control_b} shows results for a profit-oriented operator focused on minimizing the travel time of its own fleet. 
Metrics include total delay removed \( \TotalDelayRemoved \), baseload traffic delay \( \TotalDelayExogenous \), fleet delay \( \TotalDelayFleet \), 
and average delays for all trips \( \AvgTotalDelay \), baseload trips \( \AvgTotalDelayExogenous \), and fleet trips \( \AvgTotalDelayFleet \).}
\label{fig:flow_control}
\end{figure}

Figure~\ref{fig:flow_control} presents results for a representative instance (day 21), illustrating how system-wide delay evolves as $\FractionControlled$ increases. In the welfare-oriented case, we observe diminishing returns, with the most substantial marginal delay reductions occurring when fewer than 50\% of trips are under control. In contrast, the profit-driven scenario exhibits a more linear trend, as improvements are concentrated on the operator's own fleet.

In the welfare-driven scenario, average delays for \gls{amod} and baseload trips remain similar, especially for $\FractionControlled > 30\%$.
In contrast, the profit-driven controller maintains minimal delay for its fleet, while baseload traffic delays remain relatively high.
Despite their differing objectives, both strategies lead to clear reductions in total and average delays for both trip classes, suggesting that even partial control yields substantial system-wide benefits.

\begin{result}
Controlling only a fraction of the traffic with staggered and balanced routing yields a win–win outcome: AMoD fleets achieve lower delays, while baseload traffic also benefits from improved flow. These system-wide gains arise regardless of whether the operator optimizes for overall welfare or fleet-specific efficiency.
\end{result}

%% file: figures/experiments/math_comparison/legend.tex
\makebox[\textwidth]{%
  \fbox{%
    \begin{tikzpicture}
    \centering
        \definecolor{MATHColor}{HTML}{FBA44C}   
        \definecolor{LNSColor}{HTML}{0F5D91}    

        \node at (-0.7, 0.2) {}; 

        \draw [fill=MATHColor, draw=black] (0,0) rectangle (0.4,0.4);
        \node [right] at (0.5,0.2) {\texttt{MATH}};

        \draw [fill=LNSColor, draw=black] (3,0) rectangle (3.4,0.4);
        \node [right] at (3.5,0.2) {\texttt{STAG}};

        \node at (4.5, 0.2) {}; 
    \end{tikzpicture}
  }
}

\vspace{3pt}  

%% file: figures/experiments/math_comparison/relative_total_delay_reduction_congestion_HC.tex
\begin{tikzpicture}

\definecolor{cadetblue84178182}{RGB}{15,93,145}
\definecolor{darkgray176}{RGB}{176,176,176}
\definecolor{darkslategray80}{RGB}{80,80,80}
\definecolor{sandybrown22916397}{RGB}{229,163,97}

\begin{axis}[
width=\MathComparisonWidth,
height=\MathComparisonHeight,
tick align=outside,
tick pos=left,
x grid style={darkgray176},
xlabel={[\%]},
xmajorgrids,
xmin=0, xmax=104,
xtick style={color=black},
y dir=reverse,
y grid style={darkgray176},
ymin=-0.5, ymax=1.5,
ytick style={color=black},
ytick={0,1},
yticklabels={$\mathrm{\SigmaUp}$,$\mathrm{\SigmaDown}$}
]

\path [draw=darkslategray80, fill=sandybrown22916397, thick]
(axis cs:53.0966674317165,-0.4)
--(axis cs:53.0966674317165,0)
--(axis cs:76.6850969824357,0)
--(axis cs:76.6850969824357,-0.4)
--(axis cs:53.0966674317165,-0.4)
--cycle;
\addplot [thick, darkslategray80]
table {%
53.0966674317165 -0.2
18.1322972724184 -0.2
};
\addplot [thick, darkslategray80]
table {%
76.6850969824357 -0.2
89.7809111665916 -0.2
};
\addplot [thick, darkslategray80]
table {%
18.1322972724184 -0.3
18.1322972724184 -0.1
};
\addplot [thick, darkslategray80]
table {%
89.7809111665916 -0.3
89.7809111665916 -0.1
};
\path [draw=darkslategray80, fill=sandybrown22916397, thick]
(axis cs:54.5214301285979,0.6)
--(axis cs:54.5214301285979,1)
--(axis cs:71.52898338185,1)
--(axis cs:71.52898338185,0.6)
--(axis cs:54.5214301285979,0.6)
--cycle;
\addplot [thick, darkslategray80]
table {%
54.5214301285979 0.8
44.7529905671361 0.8
};
\addplot [thick, darkslategray80]
table {%
71.52898338185 0.8
87.0005738989196 0.8
};
\addplot [thick, darkslategray80]
table {%
44.7529905671361 0.7
44.7529905671361 0.9
};
\addplot [thick, darkslategray80]
table {%
87.0005738989196 0.7
87.0005738989196 0.9
};
\addplot [black, opacity=1, mark=o, mark size=2, mark options={solid,fill opacity=0}, only marks]
table {%
13.8372158263184 0.8
13.7698206411314 0.8
13.3614985127846 0.8
25.1863235021337 0.8
99.8554546890652 0.8
};
\path [draw=darkslategray80, fill=cadetblue84178182, thick]
(axis cs:79.9573608265392,5.55111512312578e-17)
--(axis cs:79.9573608265392,0.4)
--(axis cs:88.6402137191415,0.4)
--(axis cs:88.6402137191415,5.55111512312578e-17)
--(axis cs:79.9573608265392,5.55111512312578e-17)
--cycle;
\addplot [thick, darkslategray80]
table {%
79.9573608265392 0.2
71.0243421751389 0.2
};
\addplot [thick, darkslategray80]
table {%
88.6402137191415 0.2
95.2765770369853 0.2
};
\addplot [thick, darkslategray80]
table {%
71.0243421751389 0.1
71.0243421751389 0.3
};
\addplot [thick, darkslategray80]
table {%
95.2765770369853 0.1
95.2765770369853 0.3
};
\path [draw=darkslategray80, fill=cadetblue84178182, thick]
(axis cs:73.4064217577594,1)
--(axis cs:73.4064217577594,1.4)
--(axis cs:84.3061724140389,1.4)
--(axis cs:84.3061724140389,1)
--(axis cs:73.4064217577594,1)
--cycle;
\addplot [thick, darkslategray80]
table {%
73.4064217577594 1.2
62.0257626873951 1.2
};
\addplot [thick, darkslategray80]
table {%
84.3061724140389 1.2
92.657864781531 1.2
};
\addplot [thick, darkslategray80]
table {%
62.0257626873951 1.1
62.0257626873951 1.3
};
\addplot [thick, darkslategray80]
table {%
92.657864781531 1.1
92.657864781531 1.3
};
\addplot [black, opacity=1, mark=o, mark size=2, mark options={solid,fill opacity=0}, only marks]
table {%
53.1602124284498 1.2
};
\addplot [thick, black]
table {%
62.6522481458051 -0.4
62.6522481458051 0
};
\addplot [thick, black]
table {%
59.8550734652765 0.6
59.8550734652765 1
};
\addplot [thick, black]
table {%
86.6510991330175 5.55111512312578e-17
86.6510991330175 0.4
};
\addplot [thick, black]
table {%
78.7680914845796 1
78.7680914845796 1.4
};
\end{axis}

\end{tikzpicture}

%% file: figures/experiments/basic_setting_performance/legend.tex
\hspace{20pt}
\begin{minipage}{\textwidth}
\definecolor{BlueTotal}{RGB}{111,168,220}      
\definecolor{OrangeCongestion}{RGB}{230,145,56} 
\definecolor{BlueDetour}{RGB}{15,93,145}       

\centering
\fbox{
\begin{tabular}{cccccc}
\tikz\draw[draw=black, fill=BlueTotal] (0,0) rectangle (0.5,0.3); & \hspace{4pt} Total delay ↓ \hspace{10pt} &
\tikz\draw[draw=black, fill=OrangeCongestion] (0,0) rectangle (0.5,0.3); & \hspace{4pt} Congestion delay ↓ \hspace{10pt} &
\tikz\draw[draw=black, fill=BlueDetour] (0,0) rectangle (0.5,0.3); & \hspace{4pt} Detour delay ↑ \\
\end{tabular}
}
\end{minipage}

%% file: figures/experiments/basic_setting_performance/relative_delay_reduction_boxplot.tex
\begin{tikzpicture}

\definecolor{cornflowerblue124167206}{RGB}{124,167,206}
\definecolor{darkgray176}{RGB}{176,176,176}
\definecolor{dimgray85}{RGB}{85,85,85}
\definecolor{peru20814477}{RGB}{208,144,77}

\begin{axis}[
height=\DelayReductionHeight,
tick align=outside,
tick pos=left,
width=\DelayReductionWidth,
x grid style={darkgray176},
xlabel={[\%]},
xmajorgrids,
xmin=0, xmax=34.524288078096,
xtick style={color=black},
y dir=reverse,
y grid style={darkgray176},
ymajorgrids,
ymin=-0.5, ymax=2.5,
ytick style={color=black},
ytick={0,1,2},
yticklabels={\texttt{BAL},\texttt{STAG},\texttt{INTEG}}
]
\draw[draw=dimgray85,fill=cornflowerblue124167206,thick] (axis cs:0,0) rectangle (axis cs:0,0);
\addlegendimage{ybar,ybar legend,draw=dimgray85,fill=cornflowerblue124167206,thick}

\draw[draw=dimgray85,fill=peru20814477,thick] (axis cs:0,0) rectangle (axis cs:0,0);
\addlegendimage{ybar,ybar legend,draw=dimgray85,fill=peru20814477,thick}

\path [draw=dimgray85, fill=cornflowerblue124167206, thick]
(axis cs:10.7220335638665,-0.3)
--(axis cs:10.7220335638665,0)
--(axis cs:12.481578306697,0)
--(axis cs:12.481578306697,-0.3)
--(axis cs:10.7220335638665,-0.3)
--cycle;
\addplot [thick, black]
table {%
10.7220335638665 -0.15
8.96515612261068 -0.15
};
\addplot [thick, black]
table {%
12.481578306697 -0.15
13.284441606761 -0.15
};
\addplot [thick, black]
table {%
8.96515612261068 -0.225
8.96515612261068 -0.075
};
\addplot [thick, black]
table {%
13.284441606761 -0.225
13.284441606761 -0.075
};
\path [draw=dimgray85, fill=cornflowerblue124167206, thick]
(axis cs:12.3794373367658,0.7)
--(axis cs:12.3794373367658,1)
--(axis cs:14.844545405403,1)
--(axis cs:14.844545405403,0.7)
--(axis cs:12.3794373367658,0.7)
--cycle;
\addplot [thick, black]
table {%
12.3794373367658 0.85
11.4549148236336 0.85
};
\addplot [thick, black]
table {%
14.844545405403 0.85
16.6607577267015 0.85
};
\addplot [thick, black]
table {%
11.4549148236336 0.775
11.4549148236336 0.925
};
\addplot [thick, black]
table {%
16.6607577267015 0.775
16.6607577267015 0.925
};
\path [draw=dimgray85, fill=cornflowerblue124167206, thick]
(axis cs:21.9815849889588,1.7)
--(axis cs:21.9815849889588,2)
--(axis cs:22.5096624671833,2)
--(axis cs:22.5096624671833,1.7)
--(axis cs:21.9815849889588,1.7)
--cycle;
\addplot [thick, black]
table {%
21.9815849889588 1.85
21.2377758614867 1.85
};
\addplot [thick, black]
table {%
22.5096624671833 1.85
23.1110166008832 1.85
};
\addplot [thick, black]
table {%
21.2377758614867 1.775
21.2377758614867 1.925
};
\addplot [thick, black]
table {%
23.1110166008832 1.775
23.1110166008832 1.925
};
\path [draw=dimgray85, fill=peru20814477, thick]
(axis cs:21.2422446805089,0)
--(axis cs:21.2422446805089,0.3)
--(axis cs:23.40553082202,0.3)
--(axis cs:23.40553082202,0)
--(axis cs:21.2422446805089,0)
--cycle;
\addplot [thick, black]
table {%
21.2422446805089 0.15
18.5709236411141 0.15
};
\addplot [thick, black]
table {%
23.40553082202 0.15
24.5501162130277 0.15
};
\addplot [thick, black]
table {%
18.5709236411141 0.075
18.5709236411141 0.225
};
\addplot [thick, black]
table {%
24.5501162130277 0.075
24.5501162130277 0.225
};
\path [draw=dimgray85, fill=peru20814477, thick]
(axis cs:12.792859994369,1)
--(axis cs:12.792859994369,1.3)
--(axis cs:15.2816979199495,1.3)
--(axis cs:15.2816979199495,1)
--(axis cs:12.792859994369,1)
--cycle;
\addplot [thick, black]
table {%
12.792859994369 1.15
11.7913065046907 1.15
};
\addplot [thick, black]
table {%
15.2816979199495 1.15
17.0788957888808 1.15
};
\addplot [thick, black]
table {%
11.7913065046907 1.075
11.7913065046907 1.225
};
\addplot [thick, black]
table {%
17.0788957888808 1.075
17.0788957888808 1.225
};
\path [draw=dimgray85, fill=peru20814477, thick]
(axis cs:31.2297987070279,2)
--(axis cs:31.2297987070279,2.3)
--(axis cs:32.2706749213297,2.3)
--(axis cs:32.2706749213297,2)
--(axis cs:31.2297987070279,2)
--cycle;
\addplot [thick, black]
table {%
31.2297987070279 2.15
29.736550962897 2.15
};
\addplot [thick, black]
table {%
32.2706749213297 2.15
33.3071865564062 2.15
};
\addplot [thick, black]
table {%
29.736550962897 2.075
29.736550962897 2.225
};
\addplot [thick, black]
table {%
33.3071865564062 2.075
33.3071865564062 2.225
};
\addplot [thick, black]
table {%
11.8647794931242 -0.3
11.8647794931242 0
};
\addplot [thick, black]
table {%
13.2436328767225 0.7
13.2436328767225 1
};
\addplot [thick, black]
table {%
22.2327243140107 1.7
22.2327243140107 2
};
\addplot [thick, black]
table {%
22.5085214080979 0
22.5085214080979 0.3
};
\addplot [thick, black]
table {%
13.6671784581947 1
13.6671784581947 1.3
};
\addplot [thick, black]
table {%
31.8466188049547 2
31.8466188049547 2.3
};
\end{axis}

\end{tikzpicture}

%% file: figures/experiments/basic_setting_performance/absolute_delay_reduction_boxplot.tex
\begin{tikzpicture}

\definecolor{cornflowerblue124167206}{RGB}{124,167,206}
\definecolor{darkgray176}{RGB}{176,176,176}
\definecolor{darkslategray48}{RGB}{48,48,48}
\definecolor{peru20814477}{RGB}{208,144,77}
\definecolor{teal3189128}{RGB}{31,89,128}

\begin{axis}[
height=\DelayReductionHeight,
tick align=outside,
tick pos=left,
width=\DelayReductionWidth,
x grid style={darkgray176},
xlabel={[hrs]},
xmajorgrids,
xmin=0, xmax=25.9411849375,
xtick style={color=black},
y dir=reverse,
y grid style={darkgray176},
ymajorgrids,
ymin=-0.5, ymax=2.5,
ytick style={color=black},
ytick={0,1,2},
yticklabels={}
]
\draw[draw=darkslategray48,fill=cornflowerblue124167206,thick] (axis cs:0,0) rectangle (axis cs:0,0);
\addlegendimage{ybar,ybar legend,draw=darkslategray48,fill=cornflowerblue124167206,thick}

\draw[draw=darkslategray48,fill=peru20814477,thick] (axis cs:0,0) rectangle (axis cs:0,0);
\addlegendimage{ybar,ybar legend,draw=darkslategray48,fill=peru20814477,thick}

\draw[draw=darkslategray48,fill=teal3189128,thick] (axis cs:0,0) rectangle (axis cs:0,0);
\addlegendimage{ybar,ybar legend,draw=darkslategray48,fill=teal3189128,thick}

\path [draw=darkslategray48, fill=cornflowerblue124167206, thick]
(axis cs:6.17108477777778,-0.3)
--(axis cs:6.17108477777778,-0.1)
--(axis cs:8.90352759722222,-0.1)
--(axis cs:8.90352759722222,-0.3)
--(axis cs:6.17108477777778,-0.3)
--cycle;
\addplot [thick, black]
table {%
6.17108477777778 -0.2
4.26844147222223 -0.2
};
\addplot [thick, black]
table {%
8.90352759722222 -0.2
10.0914305277778 -0.2
};
\addplot [thick, black]
table {%
4.26844147222223 -0.25
4.26844147222223 -0.15
};
\addplot [thick, black]
table {%
10.0914305277778 -0.25
10.0914305277778 -0.15
};
\path [draw=darkslategray48, fill=cornflowerblue124167206, thick]
(axis cs:8.05974970833333,0.7)
--(axis cs:8.05974970833333,0.9)
--(axis cs:8.91859658333333,0.9)
--(axis cs:8.91859658333333,0.7)
--(axis cs:8.05974970833333,0.7)
--cycle;
\addplot [thick, black]
table {%
8.05974970833333 0.8
7.26423222222223 0.8
};
\addplot [thick, black]
table {%
8.91859658333333 0.8
9.62724525 0.8
};
\addplot [thick, black]
table {%
7.26423222222223 0.75
7.26423222222223 0.85
};
\addplot [thick, black]
table {%
9.62724525 0.75
9.62724525 0.85
};
\path [draw=darkslategray48, fill=cornflowerblue124167206, thick]
(axis cs:12.7025911666667,1.7)
--(axis cs:12.7025911666667,1.9)
--(axis cs:16.0709402638889,1.9)
--(axis cs:16.0709402638889,1.7)
--(axis cs:12.7025911666667,1.7)
--cycle;
\addplot [thick, black]
table {%
12.7025911666667 1.8
10.3750186388889 1.8
};
\addplot [thick, black]
table {%
16.0709402638889 1.8
17.4237099722222 1.8
};
\addplot [thick, black]
table {%
10.3750186388889 1.75
10.3750186388889 1.85
};
\addplot [thick, black]
table {%
17.4237099722222 1.75
17.4237099722222 1.85
};
\path [draw=darkslategray48, fill=peru20814477, thick]
(axis cs:11.4746338194444,-0.1)
--(axis cs:11.4746338194444,0.1)
--(axis cs:16.3284115694444,0.1)
--(axis cs:16.3284115694444,-0.1)
--(axis cs:11.4746338194444,-0.1)
--cycle;
\addplot [thick, black]
table {%
11.4746338194444 0
8.58513663888888 0
};
\addplot [thick, black]
table {%
16.3284115694444 0
18.4848753888889 0
};
\addplot [thick, black]
table {%
8.58513663888888 -0.05
8.58513663888888 0.05
};
\addplot [thick, black]
table {%
18.4848753888889 -0.05
18.4848753888889 0.05
};
\path [draw=darkslategray48, fill=peru20814477, thick]
(axis cs:8.05974969444445,0.9)
--(axis cs:8.05974969444445,1.1)
--(axis cs:8.91859658333333,1.1)
--(axis cs:8.91859658333333,0.9)
--(axis cs:8.05974969444445,0.9)
--cycle;
\addplot [thick, black]
table {%
8.05974969444445 1
7.26423222222223 1
};
\addplot [thick, black]
table {%
8.91859658333333 1
9.62724522222221 1
};
\addplot [thick, black]
table {%
7.26423222222223 0.95
7.26423222222223 1.05
};
\addplot [thick, black]
table {%
9.62724522222221 0.95
9.62724522222221 1.05
};
\path [draw=darkslategray48, fill=peru20814477, thick]
(axis cs:17.0068685277778,1.9)
--(axis cs:17.0068685277778,2.1)
--(axis cs:22.3958245138889,2.1)
--(axis cs:22.3958245138889,1.9)
--(axis cs:17.0068685277778,1.9)
--cycle;
\addplot [thick, black]
table {%
17.0068685277778 2
13.7468850833333 2
};
\addplot [thick, black]
table {%
22.3958245138889 2
24.7058904166667 2
};
\addplot [thick, black]
table {%
13.7468850833333 1.95
13.7468850833333 2.05
};
\addplot [thick, black]
table {%
24.7058904166667 1.95
24.7058904166667 2.05
};
\path [draw=darkslategray48, fill=teal3189128, thick]
(axis cs:5.39101352777778,0.1)
--(axis cs:5.39101352777778,0.3)
--(axis cs:7.451992875,0.3)
--(axis cs:7.451992875,0.1)
--(axis cs:5.39101352777778,0.1)
--cycle;
\addplot [thick, black]
table {%
5.39101352777778 0.2
4.31669516666665 0.2
};
\addplot [thick, black]
table {%
7.451992875 0.2
8.39344486111112 0.2
};
\addplot [thick, black]
table {%
4.31669516666665 0.15
4.31669516666665 0.25
};
\addplot [thick, black]
table {%
8.39344486111112 0.15
8.39344486111112 0.25
};
\path [draw=darkslategray48, fill=teal3189128, thick]
(axis cs:-7.105427357601e-15,1.1)
--(axis cs:-7.105427357601e-15,1.3)
--(axis cs:7.105427357601e-15,1.3)
--(axis cs:7.105427357601e-15,1.1)
--(axis cs:-7.105427357601e-15,1.1)
--cycle;
\addplot [thick, black]
table {%
-7.105427357601e-15 1.2
-7.105427357601e-15 1.2
};
\addplot [thick, black]
table {%
7.105427357601e-15 1.2
1.4210854715202e-14 1.2
};
\addplot [thick, black]
table {%
-7.105427357601e-15 1.15
-7.105427357601e-15 1.25
};
\addplot [thick, black]
table {%
1.4210854715202e-14 1.15
1.4210854715202e-14 1.25
};
\path [draw=darkslategray48, fill=teal3189128, thick]
(axis cs:4.46295361111111,2.1)
--(axis cs:4.46295361111111,2.3)
--(axis cs:6.33870795833333,2.3)
--(axis cs:6.33870795833333,2.1)
--(axis cs:4.46295361111111,2.1)
--cycle;
\addplot [thick, black]
table {%
4.46295361111111 2.2
3.37186644444444 2.2
};
\addplot [thick, black]
table {%
6.33870795833333 2.2
7.28218044444444 2.2
};
\addplot [thick, black]
table {%
3.37186644444444 2.15
3.37186644444444 2.25
};
\addplot [thick, black]
table {%
7.28218044444444 2.15
7.28218044444444 2.25
};
\addplot [thick, black]
table {%
7.80058994444444 -0.3
7.80058994444444 -0.1
};
\addplot [thick, black]
table {%
8.5760523611111 0.7
8.5760523611111 0.9
};
\addplot [thick, black]
table {%
14.3928419722222 1.7
14.3928419722222 1.9
};
\addplot [thick, black]
table {%
14.3916070277778 -0.1
14.3916070277778 0.1
};
\addplot [thick, black]
table {%
8.57605236111112 0.9
8.57605236111112 1.1
};
\addplot [thick, black]
table {%
20.4055278888889 1.9
20.4055278888889 2.1
};
\addplot [thick, black]
table {%
6.62567883333333 0.1
6.62567883333333 0.3
};
\addplot [thick, black]
table {%
0 1.1
0 1.3
};
\addplot [thick, black]
table {%
5.79785813888889 2.1
5.79785813888889 2.3
};
\end{axis}

\end{tikzpicture}

%% file: figures/experiments/solution_structure/legend.tex
\hspace{20pt}
\begin{minipage}{\textwidth}
\definecolor{PaleSteelBlue}{RGB}{27,79,114}    
\definecolor{WarmOrange}{RGB}{209,131,31}      
\definecolor{SlateGreyBlue}{RGB}{166,200,227}  
\definecolor{SoftNavy}{RGB}{91,141,184}        

\centering
\fbox{
\begin{tabular}{cccccccccccccccc}
\tikz\draw[draw=black, fill=PaleSteelBlue] (0,0) rectangle (0.5,0.3); & \hspace{4pt} \texttt{RDUO} \hspace{10pt} &
\tikz\draw[draw=black, fill=WarmOrange] (0,0) rectangle (0.5,0.3); & \hspace{4pt} \texttt{BAL} \hspace{10pt} &
\tikz\draw[draw=black, fill=SlateGreyBlue] (0,0) rectangle (0.5,0.3); & \hspace{4pt} \texttt{STAG} \hspace{10pt} &
\tikz\draw[draw=black, fill=SoftNavy] (0,0) rectangle (0.5,0.3); & \hspace{4pt} \texttt{INTEG} \hspace{4pt} \\
\end{tabular}
}
\vspace{0pt}
\end{minipage}

%% file: figures/experiments/solution_structure/user_equilibrium_route_distribution.tex
\begin{tikzpicture}

\definecolor{chocolate20913131}{RGB}{209,131,31}
\definecolor{darkgray176}{RGB}{176,176,176}
\definecolor{darkslategray2779114}{RGB}{27,79,114}
\definecolor{steelblue91141184}{RGB}{91,141,184}

\begin{axis}[
height=\SolutionStructureHeight,
log basis x=10,
scaled x ticks=base 10:-4,
tick align=outside,
tick label style={font=\scriptsize},
tick pos=left,
width=\SolutionStructureWidth,
x grid style={darkgray176},
xmajorgrids,
xmin=0, xmax=171848.25,
xtick style={color=black},
xticklabel style={/pgf/number format/fixed},
y dir=reverse,
y grid style={darkgray176},
ylabel={route},
ylabel style={at={(axis description cs:\SolutionStructureXYRouteLabel,0.5)}},
ymin=-0.5, ymax=4.5,
ytick style={color=black},
ytick={0,1,2,3,4},
yticklabels={1,2,3,4,5}
]
\draw[draw=black,fill=darkslategray2779114] (axis cs:0,-0.4) rectangle (axis cs:163665,-0.133333333333333);
\addlegendimage{ybar,ybar legend,draw=black,fill=darkslategray2779114}

\draw[draw=black,fill=darkslategray2779114] (axis cs:0,0.6) rectangle (axis cs:20771,0.866666666666667);
\draw[draw=black,fill=darkslategray2779114] (axis cs:0,1.6) rectangle (axis cs:3289,1.86666666666667);
\draw[draw=black,fill=darkslategray2779114] (axis cs:0,2.6) rectangle (axis cs:471,2.86666666666667);
\draw[draw=black,fill=darkslategray2779114] (axis cs:0,3.6) rectangle (axis cs:44,3.86666666666667);
\draw[draw=black,fill=chocolate20913131] (axis cs:0,-0.133333333333333) rectangle (axis cs:141613,0.133333333333333);
\addlegendimage{ybar,ybar legend,draw=black,fill=chocolate20913131}

\draw[draw=black,fill=chocolate20913131] (axis cs:0,0.866666666666667) rectangle (axis cs:33155,1.13333333333333);
\draw[draw=black,fill=chocolate20913131] (axis cs:0,1.86666666666667) rectangle (axis cs:9530,2.13333333333333);
\draw[draw=black,fill=chocolate20913131] (axis cs:0,2.86666666666667) rectangle (axis cs:3323,3.13333333333333);
\draw[draw=black,fill=chocolate20913131] (axis cs:0,3.86666666666667) rectangle (axis cs:619,4.13333333333333);
\draw[draw=black,fill=steelblue91141184] (axis cs:0,0.133333333333333) rectangle (axis cs:144109,0.4);
\addlegendimage{ybar,ybar legend,draw=black,fill=steelblue91141184}

\draw[draw=black,fill=steelblue91141184] (axis cs:0,1.13333333333333) rectangle (axis cs:31863,1.4);
\draw[draw=black,fill=steelblue91141184] (axis cs:0,2.13333333333333) rectangle (axis cs:8833,2.4);
\draw[draw=black,fill=steelblue91141184] (axis cs:0,3.13333333333333) rectangle (axis cs:2934,3.4);
\draw[draw=black,fill=steelblue91141184] (axis cs:0,4.13333333333333) rectangle (axis cs:501,4.4);
\end{axis}

\end{tikzpicture}

%% file: figures/experiments/solution_structure/staggering_frequency_barplot.tex
\begin{tikzpicture}

\definecolor{darkgray176}{RGB}{176,176,176}
\definecolor{lightsteelblue166200227}{RGB}{166,200,227}
\definecolor{steelblue91141184}{RGB}{91,141,184}

\begin{axis}[
height=\SolutionStructureHeight,
tick align=outside,
tick label style={font=\scriptsize},
tick pos=left,
width=\SolutionStructureWidth,
x grid style={darkgray176},
xmajorgrids,
xmin=0, xmax=54770.1,
xtick style={color=black},
y dir=reverse,
y grid style={darkgray176},
ylabel={[sec]},
ylabel style={at={(axis description cs:-0.8,0.5)}},
ymin=-0.6, ymax=6.6,
ytick style={color=black},
ytick={0,1,2,3,4,5,6},
yticklabels={{0},{(0-20]},{(20-40]},{(40-60]},{(60-80]},{(80-100]},{(100-120]}}
]
\draw[draw=black,fill=lightsteelblue166200227] (axis cs:0,-0.35) rectangle (axis cs:45988,0);
\addlegendimage{ybar,ybar legend,draw=black,fill=lightsteelblue166200227}

\draw[draw=black,fill=lightsteelblue166200227] (axis cs:0,0.65) rectangle (axis cs:52162,1);
\draw[draw=black,fill=lightsteelblue166200227] (axis cs:0,1.65) rectangle (axis cs:45329,2);
\draw[draw=black,fill=lightsteelblue166200227] (axis cs:0,2.65) rectangle (axis cs:29413,3);
\draw[draw=black,fill=lightsteelblue166200227] (axis cs:0,3.65) rectangle (axis cs:11949,4);
\draw[draw=black,fill=lightsteelblue166200227] (axis cs:0,4.65) rectangle (axis cs:2460,5);
\draw[draw=black,fill=lightsteelblue166200227] (axis cs:0,5.65) rectangle (axis cs:701,6);
\draw[draw=black,fill=steelblue91141184] (axis cs:0,-5.55111512312578e-17) rectangle (axis cs:47749,0.35);
\addlegendimage{ybar,ybar legend,draw=black,fill=steelblue91141184}

\draw[draw=black,fill=steelblue91141184] (axis cs:0,1) rectangle (axis cs:52027,1.35);
\draw[draw=black,fill=steelblue91141184] (axis cs:0,2) rectangle (axis cs:45025,2.35);
\draw[draw=black,fill=steelblue91141184] (axis cs:0,3) rectangle (axis cs:29011,3.35);
\draw[draw=black,fill=steelblue91141184] (axis cs:0,4) rectangle (axis cs:11365,4.35);
\draw[draw=black,fill=steelblue91141184] (axis cs:0,5) rectangle (axis cs:2234,5.35);
\draw[draw=black,fill=steelblue91141184] (axis cs:0,6) rectangle (axis cs:661,6.35);
\end{axis}

\end{tikzpicture}

%% file: figures/experiments/solution_structure/travel_time_diff_frequency_barplot.tex
\begin{tikzpicture}

\definecolor{chocolate20913131}{RGB}{209,131,31}
\definecolor{darkgray176}{RGB}{176,176,176}
\definecolor{gray}{RGB}{128,128,128}
\definecolor{lightsteelblue166200227}{RGB}{166,200,227}
\definecolor{steelblue91141184}{RGB}{91,141,184}

\begin{axis}[
height=\SolutionStructureHeight,
log basis x=10,
scaled x ticks=base 10:-4,
tick align=outside,
tick label style={font=\scriptsize},
tick pos=left,
width=\SolutionStructureWidth,
x grid style={darkgray176},
xmajorgrids,
xmin=0, xmax=95369.4,
xtick style={color=black},
xticklabel style={/pgf/number format/fixed},
y dir=reverse,
y grid style={darkgray176},
ylabel={[sec]},
ylabel style={at={(axis description cs:\SolutionStructureXYLabel,0.5)}},
ymin=-0.7, ymax=8.7,
ytick style={color=black},
ytick={0,1,2,3,4,5,6,7,8},
yticklabels={
  {(-60, -45]},
  {(-45, -30]},
  {(-30, -15]},
  {(-15, 0]},
  0,
  {(0, 15]},
  {(15, 30]},
  {(30, 45]},
  {(45, 60]}
}
]
\draw[draw=black,fill=steelblue91141184] (axis cs:0,-0.45) rectangle (axis cs:5974,-0.15);
\addlegendimage{ybar,ybar legend,draw=black,fill=steelblue91141184}

\draw[draw=black,fill=steelblue91141184] (axis cs:0,0.55) rectangle (axis cs:14185,0.85);
\draw[draw=black,fill=steelblue91141184] (axis cs:0,1.55) rectangle (axis cs:35311,1.85);
\draw[draw=black,fill=steelblue91141184] (axis cs:0,2.55) rectangle (axis cs:77795,2.85);
\draw[draw=black,fill=steelblue91141184] (axis cs:0,3.55) rectangle (axis cs:1375,3.85);
\draw[draw=black,fill=steelblue91141184] (axis cs:0,4.55) rectangle (axis cs:42806,4.85);
\draw[draw=black,fill=steelblue91141184] (axis cs:0,5.55) rectangle (axis cs:5142,5.85);
\draw[draw=black,fill=steelblue91141184] (axis cs:0,6.55) rectangle (axis cs:1009,6.85);
\draw[draw=black,fill=steelblue91141184] (axis cs:0,7.55) rectangle (axis cs:329,7.85);
\draw[draw=black,fill=chocolate20913131] (axis cs:0,-0.15) rectangle (axis cs:4038,0.15);
\addlegendimage{ybar,ybar legend,draw=black,fill=chocolate20913131}

\draw[draw=black,fill=chocolate20913131] (axis cs:0,0.85) rectangle (axis cs:9237,1.15);
\draw[draw=black,fill=chocolate20913131] (axis cs:0,1.85) rectangle (axis cs:24576,2.15);
\draw[draw=black,fill=chocolate20913131] (axis cs:0,2.85) rectangle (axis cs:81497,3.15);
\draw[draw=black,fill=chocolate20913131] (axis cs:0,3.85) rectangle (axis cs:7053,4.15);
\draw[draw=black,fill=chocolate20913131] (axis cs:0,4.85) rectangle (axis cs:51188,5.15);
\draw[draw=black,fill=chocolate20913131] (axis cs:0,5.85) rectangle (axis cs:5750,6.15);
\draw[draw=black,fill=chocolate20913131] (axis cs:0,6.85) rectangle (axis cs:1561,7.15);
\draw[draw=black,fill=chocolate20913131] (axis cs:0,7.85) rectangle (axis cs:463,8.15);
\draw[draw=black,fill=lightsteelblue166200227] (axis cs:0,0.15) rectangle (axis cs:1102,0.45);
\addlegendimage{ybar,ybar legend,draw=black,fill=lightsteelblue166200227}

\draw[draw=black,fill=lightsteelblue166200227] (axis cs:0,1.15) rectangle (axis cs:5448,1.45);
\draw[draw=black,fill=lightsteelblue166200227] (axis cs:0,2.15) rectangle (axis cs:26645,2.45);
\draw[draw=black,fill=lightsteelblue166200227] (axis cs:0,3.15) rectangle (axis cs:90828,3.45);
\draw[draw=black,fill=lightsteelblue166200227] (axis cs:0,4.15) rectangle (axis cs:2034,4.45);
\draw[draw=black,fill=lightsteelblue166200227] (axis cs:0,5.15) rectangle (axis cs:56375,5.45);
\draw[draw=black,fill=lightsteelblue166200227] (axis cs:0,6.15) rectangle (axis cs:5236,6.45);
\draw[draw=black,fill=lightsteelblue166200227] (axis cs:0,7.15) rectangle (axis cs:302,7.45);
\draw[draw=black,fill=lightsteelblue166200227] (axis cs:0,8.15) rectangle (axis cs:17,8.45);
\addplot [semithick, gray, opacity=0.7]
table {%
0 8.7
0 -0.699999999999999
};
\end{axis}

\end{tikzpicture}

%% file: figures/experiments/solution_structure/arrival_diff_frequency_barplot.tex
\begin{tikzpicture}

\definecolor{chocolate20913131}{RGB}{209,131,31}
\definecolor{darkgray176}{RGB}{176,176,176}
\definecolor{gray}{RGB}{128,128,128}
\definecolor{lightsteelblue166200227}{RGB}{166,200,227}
\definecolor{steelblue91141184}{RGB}{91,141,184}

\begin{axis}[
height=\SolutionStructureHeight,
log basis x=10,
scaled x ticks=base 10:-4,
tick align=outside,
tick label style={font=\scriptsize},
tick pos=left,
width=\SolutionStructureWidth,
x grid style={darkgray176},
xmajorgrids,
xmin=0, xmax=111454.35,
xtick style={color=black},
xticklabel style={/pgf/number format/fixed},
y dir=reverse,
y grid style={darkgray176},
ylabel={[sec]},
ylabel style={at={(axis description cs:\SolutionStructureXYLabel,0.5)}},
ymin=-0.7, ymax=6.7,
ytick style={color=black},
ytick={0,1,2,3,4,5,6},
yticklabels={{(-90, -60]},{(-60, -30]},{(-30, 0]},0,{(0, 30]},{(30, 60]},{(60, 90]}}
]
\draw[draw=black,fill=steelblue91141184] (axis cs:0,-0.45) rectangle (axis cs:1253,-0.15);
\addlegendimage{ybar,ybar legend,draw=black,fill=steelblue91141184}

\draw[draw=black,fill=steelblue91141184] (axis cs:0,0.55) rectangle (axis cs:7586,0.85);
\draw[draw=black,fill=steelblue91141184] (axis cs:0,1.55) rectangle (axis cs:51546,1.85);
\draw[draw=black,fill=steelblue91141184] (axis cs:0,2.55) rectangle (axis cs:1020,2.85);
\draw[draw=black,fill=steelblue91141184] (axis cs:0,3.55) rectangle (axis cs:78185,3.85);
\draw[draw=black,fill=steelblue91141184] (axis cs:0,4.55) rectangle (axis cs:38038,4.85);
\draw[draw=black,fill=steelblue91141184] (axis cs:0,5.55) rectangle (axis cs:9158,5.85);
\draw[draw=black,fill=chocolate20913131] (axis cs:0,-0.15) rectangle (axis cs:2399,0.15);
\addlegendimage{ybar,ybar legend,draw=black,fill=chocolate20913131}

\draw[draw=black,fill=chocolate20913131] (axis cs:0,0.85) rectangle (axis cs:13272,1.15);
\draw[draw=black,fill=chocolate20913131] (axis cs:0,1.85) rectangle (axis cs:106147,2.15);
\draw[draw=black,fill=chocolate20913131] (axis cs:0,2.85) rectangle (axis cs:6888,3.15);
\draw[draw=black,fill=chocolate20913131] (axis cs:0,3.85) rectangle (axis cs:57031,4.15);
\draw[draw=black,fill=chocolate20913131] (axis cs:0,4.85) rectangle (axis cs:2024,5.15);
\draw[draw=black,fill=chocolate20913131] (axis cs:0,5.85) rectangle (axis cs:236,6.15);
\draw[draw=black,fill=lightsteelblue166200227] (axis cs:0,0.15) rectangle (axis cs:43,0.45);
\addlegendimage{ybar,ybar legend,draw=black,fill=lightsteelblue166200227}

\draw[draw=black,fill=lightsteelblue166200227] (axis cs:0,1.15) rectangle (axis cs:1472,1.45);
\draw[draw=black,fill=lightsteelblue166200227] (axis cs:0,2.15) rectangle (axis cs:43774,2.45);
\draw[draw=black,fill=lightsteelblue166200227] (axis cs:0,3.15) rectangle (axis cs:1600,3.45);
\draw[draw=black,fill=lightsteelblue166200227] (axis cs:0,4.15) rectangle (axis cs:90174,4.45);
\draw[draw=black,fill=lightsteelblue166200227] (axis cs:0,5.15) rectangle (axis cs:40608,5.45);
\draw[draw=black,fill=lightsteelblue166200227] (axis cs:0,6.15) rectangle (axis cs:9368,6.45);
\addplot [semithick, gray, opacity=0.7]
table {%
0 6.7
0 -0.7
};
\end{axis}

\end{tikzpicture}

%% file: figures/experiments/flow_control/legend_plot.tex
\hspace{15pt}
\begin{minipage}{\textwidth}
\definecolor{TotalDelayRemoved}{RGB}{217,217,217}     
\definecolor{TotalDelayFleet}{RGB}{217,130,27}        
\definecolor{TotalDelayExogenous}{RGB}{59,141,189}    

\centering
\fbox{
\begin{tabular}{cccccccccccc}
    \tikz\draw[draw=black, fill=TotalDelayRemoved] (0,0) rectangle (0.5,0.3); & $\TotalDelayRemoved$ & 
    \tikz\draw[draw=black, fill=TotalDelayExogenous] (0,0) rectangle (0.5,0.3); & $\TotalDelayExogenous$ & 
    \tikz\draw[draw=black, fill=TotalDelayFleet] (0,0) rectangle (0.5,0.3); & $\TotalDelayFleet$ &

    \tikz{
        \draw[draw=black, line width=1pt] (0,0.15) -- (0.5,0.15); 
        \draw[draw=black, fill=white] (0.25,0.15) circle (0.1);    
    } & $\AvgTotalDelay$ & 
    \tikz{
        \draw[draw=black, line width=1pt] (0,0.15) -- (0.5,0.15); 
        \draw[draw=black, fill=white] (0.25,0.3) -- (0.1,0.05) -- (0.4,0.05) -- cycle; 
    }  & $\AvgTotalDelayExogenous$ &
    \tikz{
        \draw[draw=black, line width=1pt] (0,0.15) -- (0.5,0.15); 
        \draw[draw=black, fill=white] (0.15,0.05) rectangle (0.35,0.25); 
    } & $\AvgTotalDelayFleet$
\end{tabular}
}
\vspace{0pt}
\end{minipage}

%% file: figures/experiments/flow_control/welfare_goal_plot.tex
\begin{tikzpicture}

\definecolor{chocolate21713027}{RGB}{217,130,27}
\definecolor{darkgray176}{RGB}{176,176,176}
\definecolor{gainsboro217}{RGB}{217,217,217}
\definecolor{gray}{RGB}{128,128,128}
\definecolor{steelblue59141189}{RGB}{59,141,189}

\begin{axis}[
height=\HeightFlowControl,
tick align=outside,
tick pos=left,
width=\WidthFlowControl,
x grid style={darkgray176},
xlabel={\(\displaystyle \mathrm{\FractionControlled}\) [\%]},
xmin=-5, xmax=105,
xtick style={color=black},
y grid style={gray},
ylabel={[hrs]},
ymajorgrids,
ymin=0, ymax=94.9566825555555,
ytick style={color=black},
axis line style={-},
]
\draw[draw=black,fill=gainsboro217] (axis cs:-5,0) rectangle (axis cs:5,94.9566825555555);
\draw[draw=black,fill=gainsboro217] (axis cs:5,0) rectangle (axis cs:15,94.9566825555555);
\draw[draw=black,fill=gainsboro217] (axis cs:15,0) rectangle (axis cs:25,94.9566825555555);
\draw[draw=black,fill=gainsboro217] (axis cs:25,0) rectangle (axis cs:35,94.9566825555555);
\draw[draw=black,fill=gainsboro217] (axis cs:35,0) rectangle (axis cs:45,94.9566825555555);
\draw[draw=black,fill=gainsboro217] (axis cs:45,0) rectangle (axis cs:55,94.9566825555555);
\draw[draw=black,fill=gainsboro217] (axis cs:55,0) rectangle (axis cs:65,94.9566825555555);
\draw[draw=black,fill=gainsboro217] (axis cs:65,0) rectangle (axis cs:75,94.9566825555555);
\draw[draw=black,fill=gainsboro217] (axis cs:75,0) rectangle (axis cs:85,94.9566825555555);
\draw[draw=black,fill=gainsboro217] (axis cs:85,0) rectangle (axis cs:95,94.9566825555555);
\draw[draw=black,fill=gainsboro217] (axis cs:95,0) rectangle (axis cs:105,94.9566825555555);
\draw[draw=black,fill=chocolate21713027] (axis cs:-5,0) rectangle (axis cs:5,94.9566825555555);
\draw[draw=black,fill=chocolate21713027] (axis cs:5,0) rectangle (axis cs:15,86.1715320277778);
\draw[draw=black,fill=chocolate21713027] (axis cs:15,0) rectangle (axis cs:25,80.4924434166667);
\draw[draw=black,fill=chocolate21713027] (axis cs:25,0) rectangle (axis cs:35,75.6935336666667);
\draw[draw=black,fill=chocolate21713027] (axis cs:35,0) rectangle (axis cs:45,72.2941901666667);
\draw[draw=black,fill=chocolate21713027] (axis cs:45,0) rectangle (axis cs:55,69.1581452222222);
\draw[draw=black,fill=chocolate21713027] (axis cs:55,0) rectangle (axis cs:65,66.8394303333333);
\draw[draw=black,fill=chocolate21713027] (axis cs:65,0) rectangle (axis cs:75,65.3055638611111);
\draw[draw=black,fill=chocolate21713027] (axis cs:75,0) rectangle (axis cs:85,63.1848929444444);
\draw[draw=black,fill=chocolate21713027] (axis cs:85,0) rectangle (axis cs:95,61.7645139444444);
\draw[draw=black,fill=chocolate21713027] (axis cs:95,0) rectangle (axis cs:105,60.6841322222222);
\draw[draw=black,fill=steelblue59141189] (axis cs:-5,0) rectangle (axis cs:5,94.9566825555555);
\draw[draw=black,fill=steelblue59141189] (axis cs:5,8.10714280555556) rectangle (axis cs:15,86.1715320277778);
\draw[draw=black,fill=steelblue59141189] (axis cs:15,15.5926771944444) rectangle (axis cs:25,80.4924434166667);
\draw[draw=black,fill=steelblue59141189] (axis cs:25,21.7702643888889) rectangle (axis cs:35,75.6935336666667);
\draw[draw=black,fill=steelblue59141189] (axis cs:35,28.18512975) rectangle (axis cs:45,72.2941901666667);
\draw[draw=black,fill=steelblue59141189] (axis cs:45,33.9910597222222) rectangle (axis cs:55,69.1581452222222);
\draw[draw=black,fill=steelblue59141189] (axis cs:55,40.5296876666667) rectangle (axis cs:65,66.8394303333333);
\draw[draw=black,fill=steelblue59141189] (axis cs:65,46.0044338888889) rectangle (axis cs:75,65.3055638611111);
\draw[draw=black,fill=steelblue59141189] (axis cs:75,50.6862158055556) rectangle (axis cs:85,63.1848929444444);
\draw[draw=black,fill=steelblue59141189] (axis cs:85,55.6611718333333) rectangle (axis cs:95,61.7645139444444);
\draw[draw=black,fill=steelblue59141189] (axis cs:95,60.6841322222222) rectangle (axis cs:105,60.6841322222222);
\end{axis}

\begin{axis}[
axis y line=right,
height=\HeightFlowControl,
tick align=outside,
width=\WidthFlowControl,
x grid style={darkgray176},
xmin=-5, xmax=105,
xtick pos=left,
xtick style={color=black},
y grid style={darkgray176},
ylabel={[min/veh]},
ymin=0.499, ymax=0.901,
ytick pos=right,
ytick style={color=black},
axis line style={-},
ytick={0.5, 0.6, 0.7, 0.8, 0.9},
yticklabels={0.5, 0.6, 0.7, 0.8, 0.9},
yticklabel style={anchor=west},
]
\addplot [thick, black, mark=square*, mark size=2.5, mark options={solid,fill=white}]
table {%
10 0.697888907221425
20 0.680407732121212
30 0.64889014571949
40 0.633848495127436
50 0.611533308345827
60 0.601182017305315
70 0.58841740211753
80 0.569081764283932
90 0.559221418285332
100 0.549841125541126
};
\addplot [thick, black, mark=*, mark size=2.5, mark options={solid,fill=white}]
table {%
0 0.860374653176281
10 0.780774980620155
20 0.729318424192087
30 0.685836910299003
40 0.655036455753549
50 0.626621672203765
60 0.605612476593174
70 0.591714562317527
80 0.572499785059901
90 0.559630147488171
100 0.549841125541126
};
\addplot [thick, black, mark=triangle*, mark size=3.5, mark options={solid,fill=white}]
table {%
0 0.860374653176281
10 0.7905254604782
20 0.742135691506257
30 0.701973564041368
40 0.669333238492666
50 0.641930371159112
60 0.612566767559177
70 0.59972439064388
80 0.586792353938445
90 0.563385425641026
};
\end{axis}

\end{tikzpicture}

%% file: figures/experiments/flow_control/profit_goal_plot.tex
\begin{tikzpicture}

\definecolor{chocolate21713027}{RGB}{217,130,27}
\definecolor{darkgray176}{RGB}{176,176,176}
\definecolor{gainsboro217}{RGB}{217,217,217}
\definecolor{gray}{RGB}{128,128,128}
\definecolor{steelblue59141189}{RGB}{59,141,189}

\begin{axis}[
height=\HeightFlowControl,
tick align=outside,
tick pos=left,
width=\WidthFlowControl,
x grid style={darkgray176},
xlabel={\(\displaystyle \mathrm{\FractionControlled}\) [\%]},
xmin=-5, xmax=105,
xtick style={color=black},
y grid style={gray},
ylabel={[hrs]},
ymajorgrids,
ymin=0, ymax=94.9566825555555,
ytick style={color=black},
axis line style={-},
]
\draw[draw=black,fill=gainsboro217] (axis cs:-5,0) rectangle (axis cs:5,94.9566825555555);
\draw[draw=black,fill=gainsboro217] (axis cs:5,0) rectangle (axis cs:15,94.9566825555555);
\draw[draw=black,fill=gainsboro217] (axis cs:15,0) rectangle (axis cs:25,94.9566825555555);
\draw[draw=black,fill=gainsboro217] (axis cs:25,0) rectangle (axis cs:35,94.9566825555555);
\draw[draw=black,fill=gainsboro217] (axis cs:35,0) rectangle (axis cs:45,94.9566825555555);
\draw[draw=black,fill=gainsboro217] (axis cs:45,0) rectangle (axis cs:55,94.9566825555555);
\draw[draw=black,fill=gainsboro217] (axis cs:55,0) rectangle (axis cs:65,94.9566825555555);
\draw[draw=black,fill=gainsboro217] (axis cs:65,0) rectangle (axis cs:75,94.9566825555555);
\draw[draw=black,fill=gainsboro217] (axis cs:75,0) rectangle (axis cs:85,94.9566825555555);
\draw[draw=black,fill=gainsboro217] (axis cs:85,0) rectangle (axis cs:95,94.9566825555555);
\draw[draw=black,fill=gainsboro217] (axis cs:95,0) rectangle (axis cs:105,94.9566825555555);
\draw[draw=black,fill=chocolate21713027] (axis cs:-5,0) rectangle (axis cs:5,94.9566825555555);
\draw[draw=black,fill=chocolate21713027] (axis cs:5,0) rectangle (axis cs:15,88.5492049166667);
\draw[draw=black,fill=chocolate21713027] (axis cs:15,0) rectangle (axis cs:25,83.4229662777778);
\draw[draw=black,fill=chocolate21713027] (axis cs:25,0) rectangle (axis cs:35,79.0110194444444);
\draw[draw=black,fill=chocolate21713027] (axis cs:35,0) rectangle (axis cs:45,75.3805807222222);
\draw[draw=black,fill=chocolate21713027] (axis cs:45,0) rectangle (axis cs:55,71.9060611944444);
\draw[draw=black,fill=chocolate21713027] (axis cs:55,0) rectangle (axis cs:65,68.7763936111111);
\draw[draw=black,fill=chocolate21713027] (axis cs:65,0) rectangle (axis cs:75,66.4369185833333);
\draw[draw=black,fill=chocolate21713027] (axis cs:75,0) rectangle (axis cs:85,64.39818325);
\draw[draw=black,fill=chocolate21713027] (axis cs:85,0) rectangle (axis cs:95,62.5361924166667);
\draw[draw=black,fill=chocolate21713027] (axis cs:95,0) rectangle (axis cs:105,60.6841322222222);
\draw[draw=black,fill=steelblue59141189] (axis cs:-5,0) rectangle (axis cs:5,94.9566825555555);
\draw[draw=black,fill=steelblue59141189] (axis cs:5,6.50473425) rectangle (axis cs:15,88.5492049166667);
\draw[draw=black,fill=steelblue59141189] (axis cs:15,13.1294416388889) rectangle (axis cs:25,83.4229662777778);
\draw[draw=black,fill=steelblue59141189] (axis cs:25,19.4329360833333) rectangle (axis cs:35,79.0110194444444);
\draw[draw=black,fill=steelblue59141189] (axis cs:35,25.6706296111111) rectangle (axis cs:45,75.3805807222222);
\draw[draw=black,fill=steelblue59141189] (axis cs:45,31.6326480277778) rectangle (axis cs:55,71.9060611944444);
\draw[draw=black,fill=steelblue59141189] (axis cs:55,38.4705229722222) rectangle (axis cs:65,68.7763936111111);
\draw[draw=black,fill=steelblue59141189] (axis cs:65,43.9969150833333) rectangle (axis cs:75,66.4369185833333);
\draw[draw=black,fill=steelblue59141189] (axis cs:75,50.0868016944444) rectangle (axis cs:85,64.39818325);
\draw[draw=black,fill=steelblue59141189] (axis cs:85,55.42160175) rectangle (axis cs:95,62.5361924166667);
\draw[draw=black,fill=steelblue59141189] (axis cs:95,60.6841322222222) rectangle (axis cs:105,60.6841322222222);
\end{axis}

\begin{axis}[
axis y line=right,
height=\HeightFlowControl,
tick align=outside,
width=\WidthFlowControl,
x grid style={darkgray176},
xmin=-5, xmax=105,
xtick pos=left,
xtick style={color=black},
y grid style={darkgray176},
ylabel={[min/veh]},
ymin=0.499, ymax=0.901,
ytick pos=right,
ytick style={color=black},
axis line style={-},
ytick={0.5, 0.6, 0.7, 0.8, 0.9},
yticklabels={0.5, 0.6, 0.7, 0.8, 0.9},
yticklabel style={anchor=west},
]
\addplot [thick, black, mark=square*, mark size=2.5, mark options={solid,fill=white}]
table {%
10 0.559948428981349
20 0.57292108969697
30 0.579223132141083
40 0.577300515992004
50 0.569103112943528
60 0.570638165224557
70 0.562740333617566
80 0.562351815431637
90 0.55681448509712
100 0.549841125541126
};
\addplot [thick, black, mark=*, mark size=2.5, mark options={solid,fill=white}]
table {%
0 0.860374653176281
10 0.802318377378435
20 0.755871032417195
30 0.715895676029397
40 0.683001335447498
50 0.651519732960838
60 0.623162732809826
70 0.601965435668982
80 0.583493052703111
90 0.566622099818786
100 0.549841125541126
};
\addplot [thick, black, mark=triangle*, mark size=3.5, mark options={solid,fill=white}]
table {%
0 0.860374653176281
10 0.830830082700422
20 0.803813889524173
30 0.775587980400665
40 0.754323992581352
50 0.735139881350776
60 0.705608163885655
70 0.697255416882444
80 0.671895847678665
90 0.656731446153846
};
\end{axis}

\end{tikzpicture}

%% file: contents/5.Conclusion.tex
\section{Conclusion} \label{sec:conclusion}

In this work, we introduced a novel framework for integrating staggered and balanced routing into \gls{amod} system operations. Our approach minimizes fleet travel times while enhancing overall network efficiency by strategically assigning trips to optimized routes and departure time patterns. We described the problem through a mathematical formulation of the problem and showed that, under mild assumptions, our congestion model provides an unbiased approximation of a discretized Vickrey bottleneck model.

To solve large-scale instances of the problem, we developed a custom metaheuristic based on a \gls{lns} framework. In a case study using real-world data from taxi trips in the Manhattan street network, we addressed a realistic volume of peak-hour demand approximating real-world traffic congestion levels. To generate route alternatives, we employed a $k$-shortest path algorithm with overlap constraints, enabling a diverse yet efficient set of routing options while maintaining computational tractability. Our results show that, under full control, our method can reduce total delays by up to 25\% compared to selfish routing. When focusing on network congestion, improvements reach up to 35\%. These gains result from redistributing trips onto slightly longer routes and applying modest staggering to departure times. While a small number of trips experience modest travel time increases, the broader gains across the network outweigh these losses. Spatially, the improvements result from congestion relief at bottlenecks, without generating new hotspots elsewhere in the network.

We also examined how the system responds to varying fractions of controlled traffic under two paradigms: a welfare-oriented authority optimizing total system travel time, and a profit-driven \gls{amod} operator focused solely on fleet performance. In both scenarios, staggered and balanced routing consistently improve outcomes for both AMoD and exogenous trips across all control levels.

This research opens several avenues for future work. First, replacing the assumption of precomputed routes with methods that compute globally optimal ones could unlock further improvements in balancing operations. Second, advancing departure time selection policies may further enhance overall system performance. Third, evaluating the framework under higher congestion levels would provide a more realistic assessment of algorithm performance in real-world conditions. Additionally, incorporating scenarios with incomplete demand information would improve realism and allow benchmarking the cost of limited information. Finally, future work should explore how baseload traffic dynamically adapts to the behavior of \gls{amod} services.

%% file: contents/6.Appendix.tex
\section{Theorem~\ref{thm:vickrey} proof}\label{appendix:vickrey_proof}
\vspace*{-\baselineskip}  
\vspace{.3cm}
\proof{Proof}
Consider the travel time function defined in equation \eqref{eq:travel_time_linear}. Taking expectations, we have:
\[
\mathbb{E}[T(t)] = \tau_a + (\phi \cdot \tau_a)\mathbb{E}[f(t)].
\]
To find the expected number of trips \(\mathbb{E}[f(t)]\) in the system, we apply Little's Law, which relates the expected number of trips in a stable system to the arrival rate and the expected time a trip spends in the system:
\[
\mathbb{E}[f(t)] = \lambda \mathbb{E}[T(t)].
\]
Substituting the expression for $\mathbb{E}[T(t)]$ into $\mathbb{E}[f(t)]$, we have:
\begin{align*}
&\mathbb{E}[f(t)] = \rho + (\phi \cdot \rho) \, \mathbb{E}[f(t)] \\
&\qquad \Longrightarrow \mathbb{E}[f(t)] - (\phi \cdot \rho) \, \mathbb{E}[f(t)] = \rho \\
& \qquad \Longrightarrow \mathbb{E}[f(t)] = {\rho}\cdot{(1 - \phi \cdot \rho)^{-1}}.
\end{align*}
The expected number of trips \(\mathbb{E}[f(t)]\) is finite when the condition \(\phi \cdot \rho < 1\) holds. Since \(\phi \in (0, 1]\) and \(\rho \in (0, 1)\), this inequality is always satisfied, ensuring system stability for any valid choice of \(\phi\).
Substituting \(\mathbb{E}[f(t)]\) back into the expression for $\mathbb{E}[T(t)]$:
\[
\mathbb{E}[T(t)] \!=\! \tau_a + (\phi \cdot \tau_a) \mathbb{E}[f(t)] = \tau_a + ({\phi \cdot \tau_a \cdot \rho})({1 - \phi \cdot \rho})^{-1}\!.
\]
To match the expected travel time with that of Vickrey's model in Equation~\eqref{eq:expected_travel_time_vickrey}, we set:
\begin{align*}
&\tau_a + \left(\phi \cdot \tau_a \cdot \rho\right) \cdot (1 - \phi \cdot \rho)^{-1}
= \tau_a + \left(\tau_a \cdot \rho\right) \cdot \left(2 (1 - \rho)\right)^{-1} \\
& \qquad\Longrightarrow \left(\phi \cdot \tau_a \cdot \rho\right) \cdot (1 - \phi \cdot \rho)^{-1}
= \left(\tau_a \cdot \rho\right) \cdot \left(2 (1 - \rho)\right)^{-1} \\
& \qquad\Longrightarrow \phi \cdot (1 - \phi \cdot \rho)^{-1}
= 2 (1 - \rho)^{-1} \\
& \qquad \Longrightarrow 2 \phi \cdot (1 - \rho)
= 1 - \phi \cdot \rho.
\end{align*}
Rearranging the expression to solve for $\phi$ gives:
\[
\phi = {(2 - \rho)^{-1}}.
\]
Since $0 < \rho \leq 1$, it follows that $\phi \in (0, 1]$.
This demonstrates that, for this choice of $\phi$, the expected travel time $\mathbb{E}[T(t)]$ matches the expected travel time in Vickrey's model.
\hfill\halmos

\section{Solution-finding modules pseudocode }\label{appendix:baseline-pseudocode}
{
We provide pseudocode for the two baseline assignment procedures used to generate initial solutions: the reactive dynamic user optimum (\textsc{buildRDUO}) in Algorithm~\ref{algo:rduo}, and the greedy assignment heuristic in Algorithm~\ref{algo:greedy}.}

\begin{algorithm}[!t]
\footnotesize
\caption{\textsc{buildRDUO}}\label{algo:rduo}
\hspace{\algorithmicindent}\textbf{Input:} Set of trips $\SetTrips$
\begin{algorithmic}[1]
\itemsep-0em
    \State $\pi \gets (\emptyset, \emptyset)$ \Comment{Initialize empty solution.}
    \State $\pi \gets \textsc{insertAtEarliestTimes}(\pi, \SetTrips, \texttt{EARLIEST})$ \Comment{Insert trips by increasing $\TripEarliestDeparture{\Trip}$.}
    \State \textbf{return} $\pi$
\end{algorithmic}
\end{algorithm}

\begin{algorithm}[!t]
\footnotesize
\caption{\textsc{greedyAssignment}}\label{algo:greedy}
\hspace{\algorithmicindent}\textbf{Input:} Set of trips $\SetTrips$
\vspace{-0em}
\begin{algorithmic}[1]
\itemsep-0em
    \State $\pi \gets (\emptyset, \emptyset)$ \Comment{Initialize empty solution.}
    \State $\pi \gets \textsc{insert}(\pi, \SetTrips, \texttt{DEADLINE})$ \Comment{Insert trips with tightest $\TripLatestArrival{\Trip}$ first.}
    \State $\SetTrips^{I} \gets \textsc{getInfeasibleTrips}(\pi)$
    \If{$\SetTrips^{I} \neq \emptyset$}
        \State $\pi \gets \textsc{localSearch}(\pi, \SetTrips^{I})$
    \EndIf
    \State \textbf{return} $\pi$
\end{algorithmic}
\end{algorithm}

\section{Description of \texttt{constructSchedule}}
\label{appendix:construct_schedule}
Algorithm~\ref{alg:construct_schedule} details the procedure to construct the schedule $\Schedule$ associated with $\Solution$.
We start by sorting the trip start times $\StartTime{\Trip}$ into a priority queue $\Queue$ and initialize an empty set $\SetOfArcArrivals{\Arc}$ for each arc $\Arc$ to track when each trip completes the traversal of $\Arc$. Then, we construct an initial solution iteratively: 
\begin{enumerate}
\item[i)] We start by extracting the departure time $\Departure{\Arc}{\Trip}$ with the highest priority from $\Queue$. 
\item[ii)] We then set $\Flow{\Arc}{\Trip} = |\{ \Arrival{\Arc}{\SecondTrip} \in \SetOfArcArrivals{\Arc} \,|\, \Departure{\Arc}{\Trip} < \Arrival{\Arc}{\SecondTrip}\}|$,
that is the number of trips that have not completed their traversal of arc $\Arc$ by time $\Departure{\Arc}{\Trip}$, and compute {$\DelayOnArc{\Arc}{\Trip} = \DelayFunction$}. 
\item[iii)] Finally, we add the time in which the arc traversal is completed, which is $\Departure{\Arc}{\Trip} + \DelayOnArc{\Arc}{\Trip} + \NominalTravelTime{\Arc}$, to $\SetOfArcArrivals{\Arc}$ by setting $\SetOfArcArrivals{\Arc} =  \SetOfArcArrivals{\Arc} \cup (\Departure{\Arc}{\Trip} + \DelayOnArc{\Arc}{\Trip} + \NominalTravelTime{\Arc})$. If $\Arc$ is not the last arc of route $\TripPath{\Trip}$, we insert $\Departure{\Arc}{\Trip} + \DelayOnArc{\Arc}{\Trip} + \NominalTravelTime{\Arc}$ into $\Queue$. 
\end{enumerate}
We repeat these steps until $\Queue$ is empty.

\section{Alternative route design}\label{appendix:route_design}
In this section, we explain the rationale for selecting alternative routes for trips by solving a \gls{kspwlo} problem. The aim is to find $k$ routes connecting an origin-destination pair, each as short as possible, while ensuring that any two have at most $\SimilarityThs$ \% of shared path length.
An important consideration is that the route selection influences both our baseline \gls{rduo} solution and our proposed algorithmic solutions since both draw routes from the same candidate set. Thus, choosing appropriate parameters $k$ and $\SimilarityThs$ is crucial.
\begin{algorithm}[!t]
    \footnotesize
    \caption{\textsc{constructSchedule}}\label{alg:construct_schedule}
    \hspace{\algorithmicindent}\textbf{Input:} Solution $\Solution$
    \vspace{-0em}
    \begin{algorithmic}[1]
    \itemsep-0em  
    \State $\Queue \gets \textsc{initializePriorityQueue}(\Solution)$
    \State $\Schedule \gets \textsc{initializeSchedule}(\Solution)$
    \While{$\Queue$ is not empty}
        \State $(\Trip, \Arc) \gets \Queue.\textsc{pop}()$
        \State $\Flow{\Arc}{\Trip} \gets \textsc{computeFlowOnArc}(\Departure{\Arc}{\Trip})$
        \State $\DelayOnArc{\Arc}{\Trip} \gets \textsc{computeDelay}(\Flow{\Arc}{\Trip})$
        \State $\Arc' \gets \textsc{getNextArc}(\Trip, \Arc, \Solution)$
        \State $\Departure{\Arc'}{\Trip} \gets \Departure{\Arc}{\Trip} + \NominalTravelTime{\Arc} + \DelayOnArc{\Arc}{\Trip}$
        \State $\Schedule \gets \textsc{addDeparture}(\Departure{\Arc'}{\Trip})$
        \If{$\Arc' \neq \LastArc{\Trip}$}
            \State $\Queue \gets \textsc{pushTrip}(\Trip, \Arc')$
        \EndIf
    \EndWhile
    \vspace{0.7em}
    \State \textbf{return} $\Schedule$
    \end{algorithmic}
\end{algorithm}
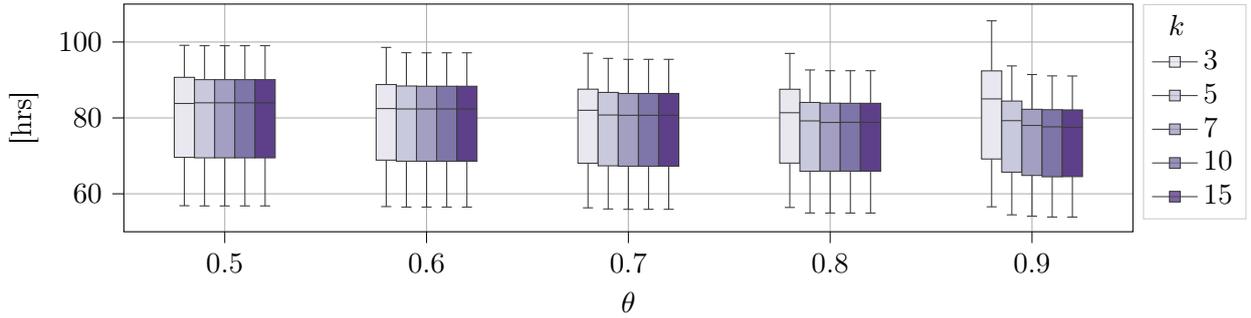
\begin{figure}[!t]
    \centering
    \import{figures/appendix/route_design/quantitative}{user_equilibrium_total_delay_hours.tex}
    \caption{Total delay in the \gls{rduo} for various combinations of \( k \) and \( \SimilarityThs \).}
    \label{fig:route_design_total_delay}
\end{figure}
Figure~\ref{fig:route_design_total_delay} illustrates the total delays in the \gls{rduo} across the 31 instances considered for various combinations of $k$ and $\SimilarityThs$. Evidently, higher values of $k$ in conjunction with larger $\SimilarityThs$ values provide users with greater flexibility in their route choices, resulting in the smallest median total delay (approximately 78 hours) observed for configurations where $\SimilarityThs = 0.9$ and $k \geq 7$. However, since our algorithm aims to evenly distribute traffic across alternative network segments to enhance overall efficiency, diversity in route selection becomes highly beneficial. In particular, the scenario with $\SimilarityThs = 0.6$ and $k = 5$ leads to an increase in median delay of approximately five hours compared to the optimal flexible scenario ($\SimilarityThs = 0.9$ and $k \geq 7$). We consider this trade-off acceptable given the clear advantage of promoting more diverse route choices and consequently improved network balance.

\begin{figure}[!t]
    \centering
    \subfigcapskip = -10pt
    \begin{minipage}{\textwidth}
        \subfigure[\(k=5\), \( \SimilarityThs = 0.6 \)\label{fig:route_design_5_0.6}]{
            \includegraphics[width=0.23\linewidth]{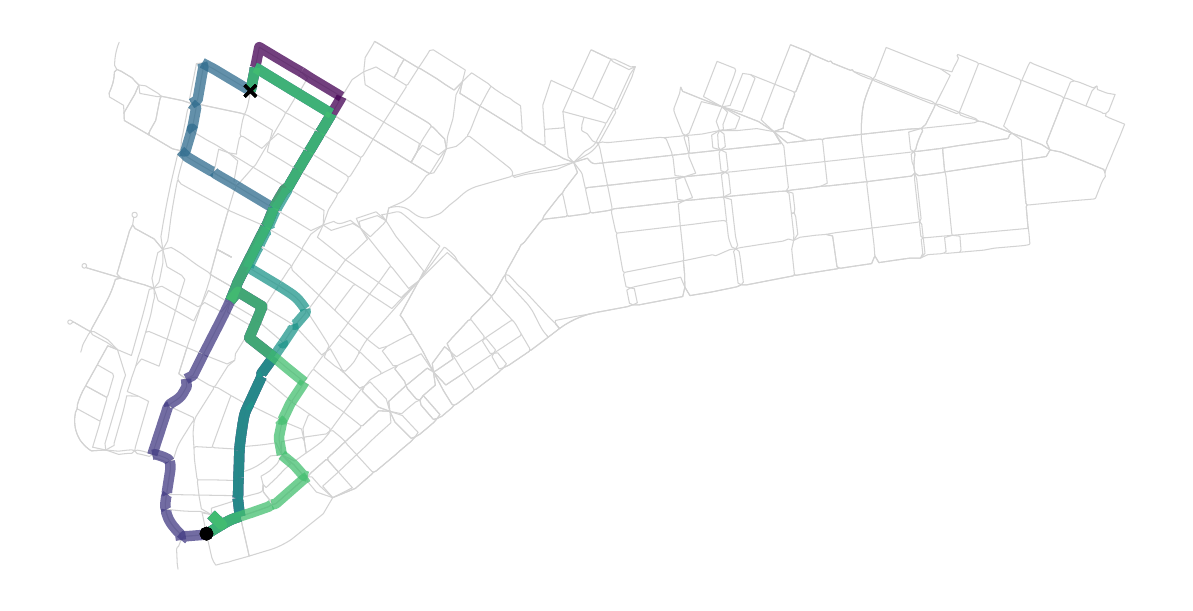}
            \hfill
            \includegraphics[width=0.23\linewidth]{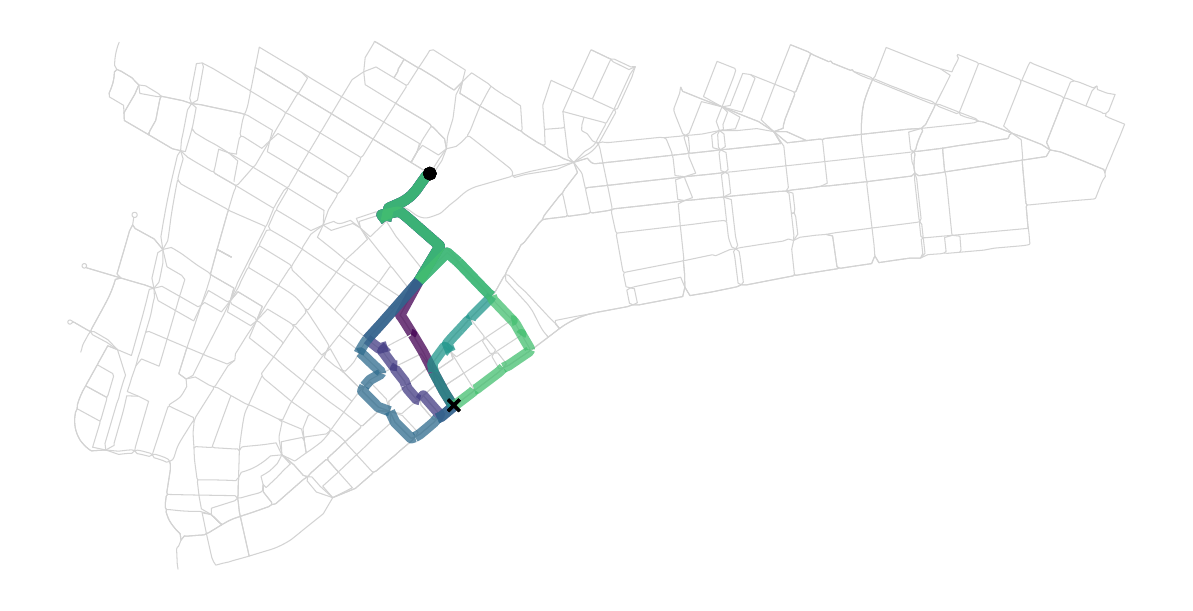}
        }
        \hfill
        \subfigure[\(k=10\), \( \SimilarityThs = 0.6 \)\label{fig:route_design_10_0.6}]{
            \includegraphics[width=0.23\linewidth]{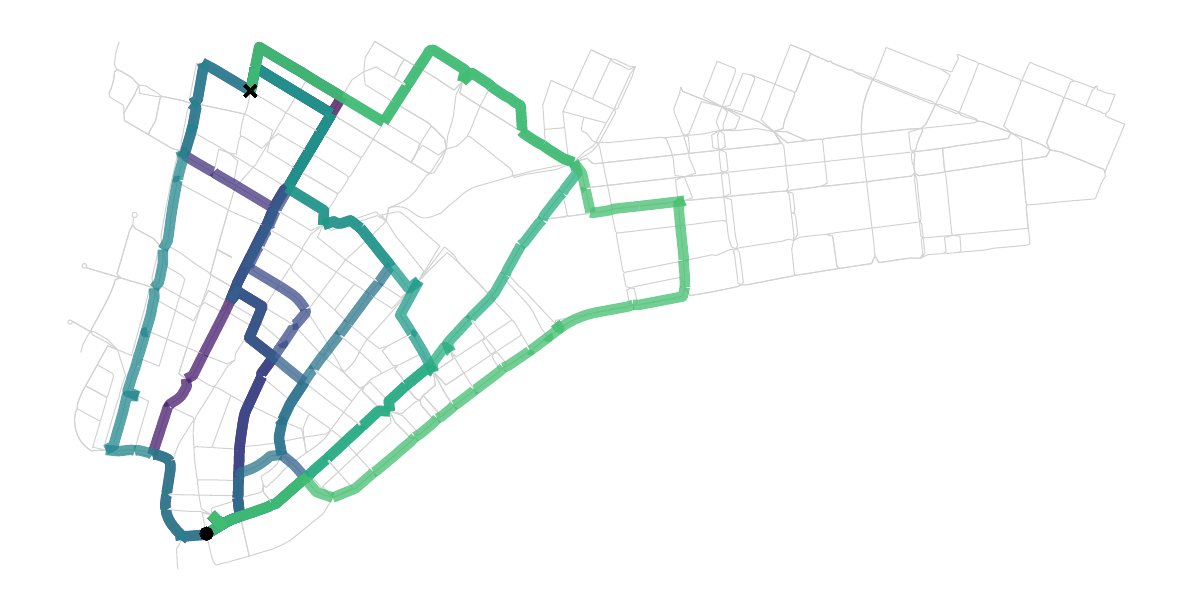}
            \hfill
            \includegraphics[width=0.23\linewidth]{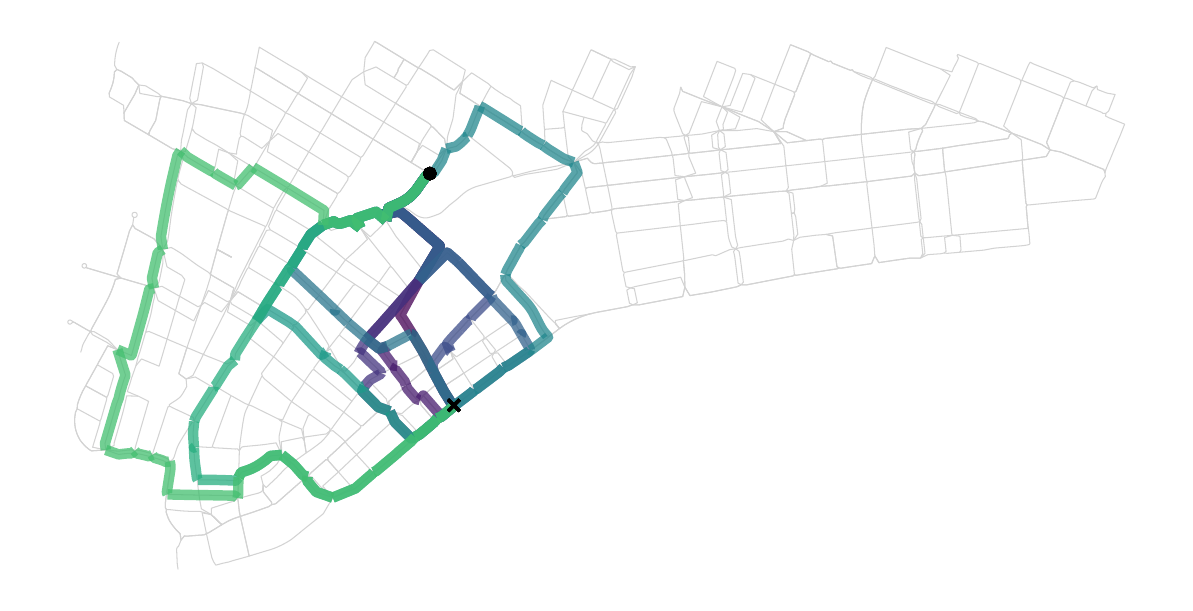}
        }

        \subfigure[\(k=5\), \( \SimilarityThs = 0.9 \)\label{fig:route_design_5_0.9}]{
            \includegraphics[width=0.23\linewidth]{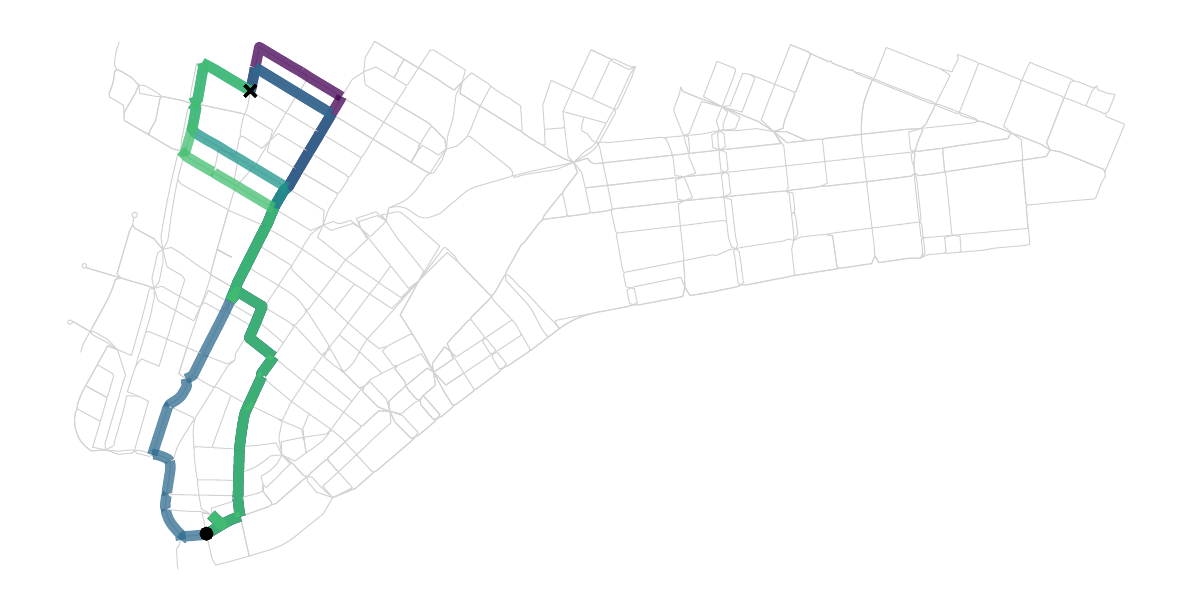}
            \hfill
            \includegraphics[width=0.23\linewidth]{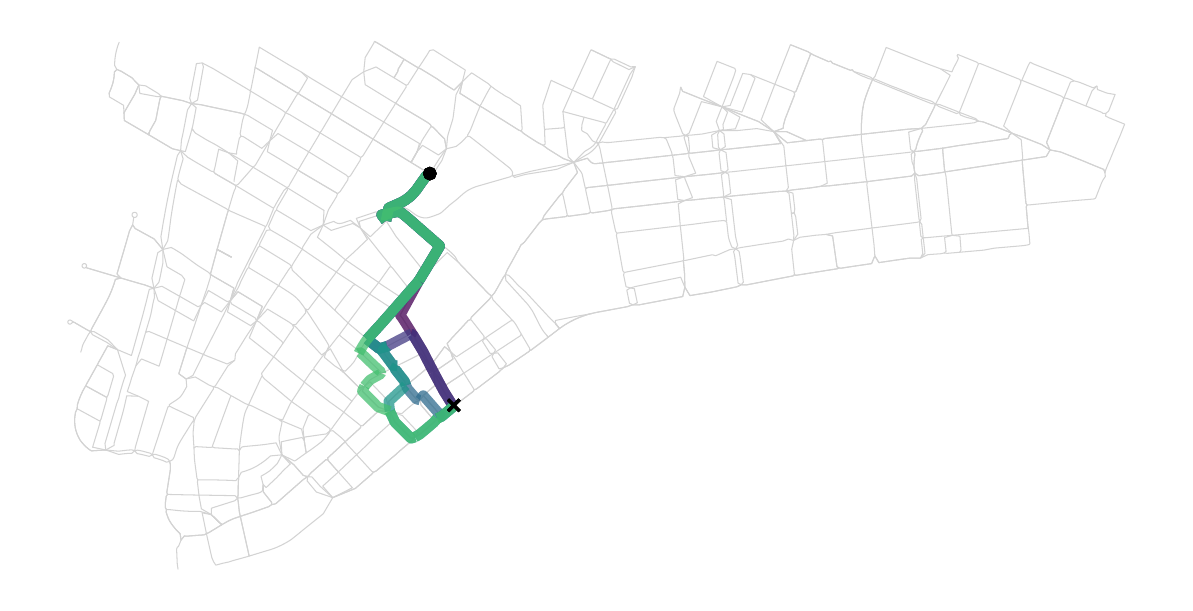}
        }
        \hfill
        \subfigure[\(k=10\), \( \SimilarityThs = 0.9 \)\label{fig:route_design_10_0.9}]{
            \includegraphics[width=0.23\linewidth]{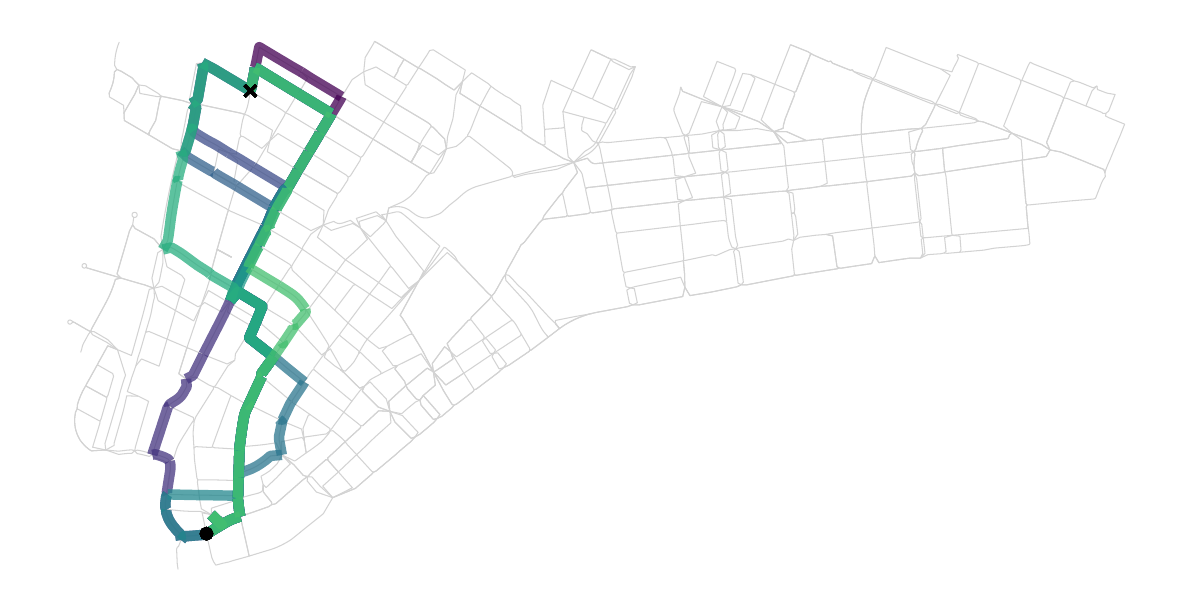}
            \hfill
            \includegraphics[width=0.23\linewidth]{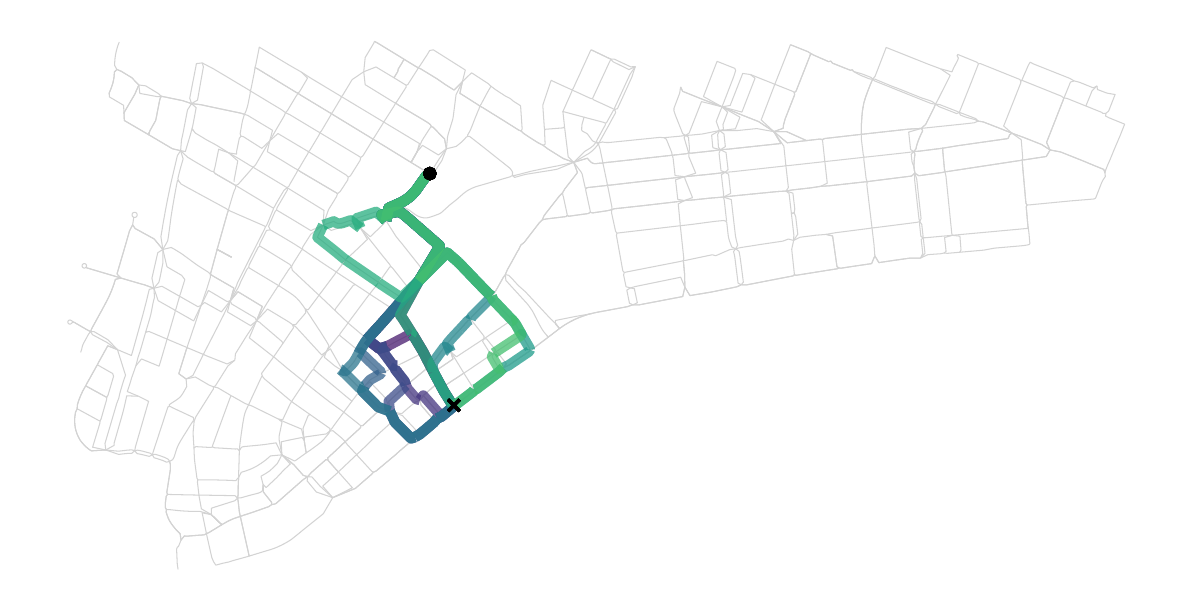}
        }
    \end{minipage}
    \caption{Alternative routes for a long and a short trip with different \( k \) and \( \theta \) values.}
    \label{fig:route_design_visual}
\end{figure}

To illustrate the visual differences among alternative route configurations, Figure~\ref{fig:route_design_visual} presents examples of route variations for both a long and a short trip under various parameter settings. High values of $k$ combined with low $\SimilarityThs$ generate routes excessively long in distance, which frequently must be discarded due to violating the latest arrival constraints even under free-flow conditions. Conversely, configurations with low $k$ and high $\SimilarityThs$ produce highly similar routes that cover limited network areas, complicating effective traffic distribution. Optimal choices thus generally fall into two categories: either low $k$ with low $\SimilarityThs$, or high $k$ with high $\SimilarityThs$. Although higher values may achieve more comprehensive network coverage, we select the smaller parameter combination ($k=5$, $\SimilarityThs=0.6$) to substantially reduce computational complexity while maintaining sufficient route diversity.

\section{Detailed instances description}\label{appendix:instances_description}
\begin{figure}[!t]
    \centering
    \subfigcapskip = -10pt

    \subfigure[Free-flow travel times\label{fig:instances_histograms_a}]{
        \import{figures/appendix/instances_description/}{travel_times_histogram.tex}
    }
    \subfigure[Relative max. detours\label{fig:instances_histograms_b}]{
        \import{figures/appendix/instances_description/}{relative_max_deviation_histogram.tex}
    }
    \subfigure[Time window lengths\label{fig:instances_histograms_c}]{
        \import{figures/appendix/instances_description/}{max_trip_durations_histogram.tex}
    }
    \subfigure[Max. staggering\label{fig:instances_histograms_d}]{
        \import{figures/appendix/instances_description/}{max_staggering_histogram.tex}
    }

    \caption{Frequency histograms summarizing key trip features across all instances: 
    (a) free-flow travel times, 
    (b) relative maximum detours between shortest and longest assigned routes, 
    (c) lengths of time windows, and 
    (d) maximum staggering allowed per trip.}
    \label{fig:instances_histograms}
\end{figure}
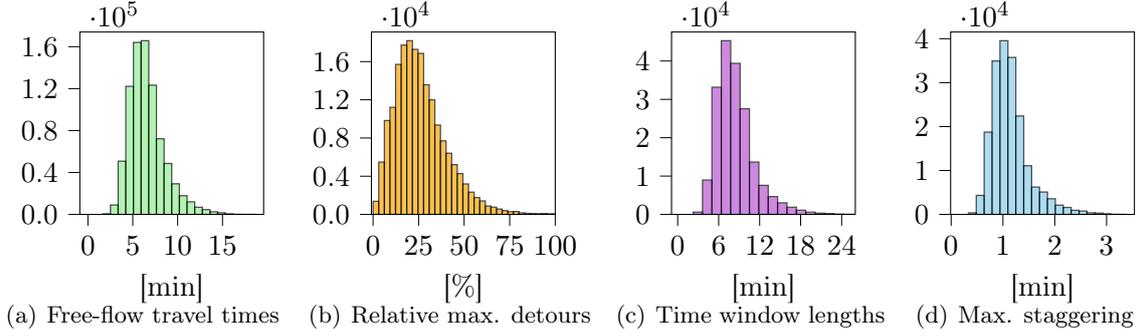

In this section, we provide a detailed overview of the characteristics of the instances under consideration. Figure~\ref{fig:instances_histograms} illustrates various key attributes of these instances aggregated together. The free-flow travel time for the routes ranges from approximately 3 to 15 minutes, including alternative route options. The relative maximum detour, defined as the percentage increase in free-flow travel time between the shortest and longest (in distance) available route for a given trip, most frequently centers around 25\%, although in some cases the longest route covers twice the distance of the shortest. Most trips have a service time window between six and twelve minutes, with a few extending up to twenty minutes. Lastly, most trips allow for a maximum staggering period of roughly one minute, with a few cases extending up to three minutes.
\input{tables/instance_summary}
Table~\ref{table:instance_summary} summarizes the main characteristics of the \gls{rduo} solution for each day considered. For each day, we report the number of trips $|\SetTrips|$, the total travel time $\TotalTravelTime{}$, the total delay $D$ in the \gls{rduo} and their ratio, as well as statistics related to the delay experienced by individual trips $\TotalDelayTrip{\Trip}$ and the ratio between their actual travel time $\TravelTimeRoute{\Trip}$ and the corresponding free-flow travel time $\FreeFlowRoute{\Trip}$. 
On average, the instances contain approximately 6000 trips, and the \gls{rduo} total delay accounts for roughly 12\% of the total travel time. The individual trip delays pass basic sanity checks: no trip experiences more than four minutes of delay, and travel times never exceed twice the corresponding free-flow travel time.

\section{Parameter sensitivity analysis of the \gls{lns}}\label{appendix:parameters_search}
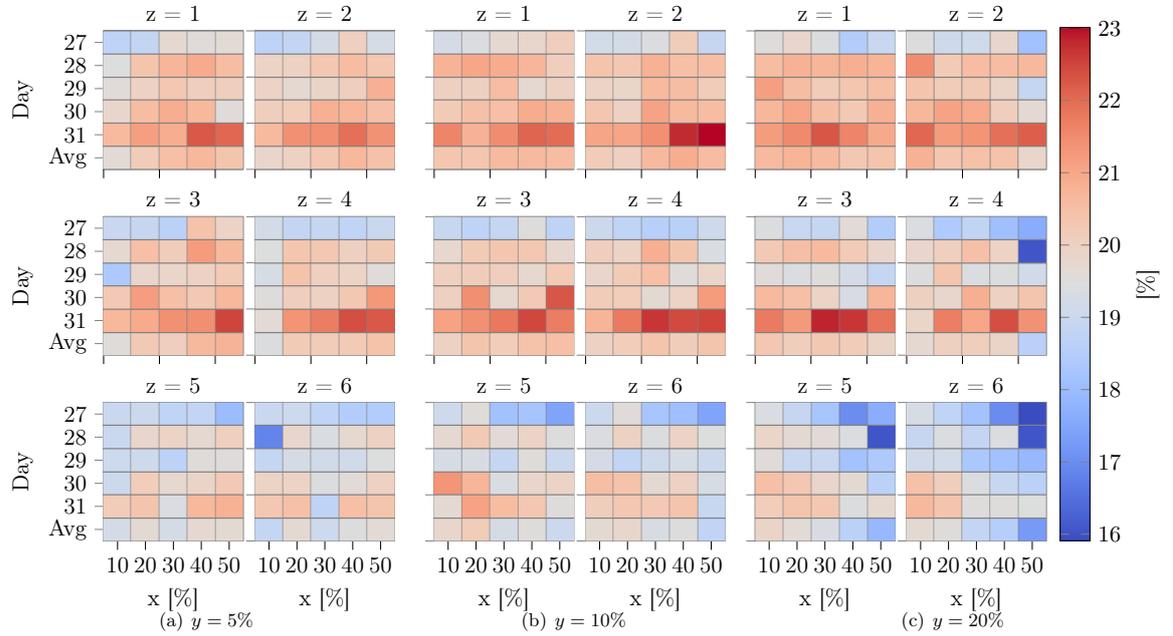
\begin{figure}[!t]
    \centering
    \subfigcapskip = -10pt
    \resizebox{.65\textwidth}{!}{
        \hspace{-3cm}
        \begin{minipage}{0.25\linewidth}
            \centering
            \subfigure[\( y = 5\% \)]{
                \import{figures/appendix/parameter_search/heatmap-grid-relative-improvement_dp_5_tex}{heatmap-grid-relative-improvement_dp_5.tex}
            }
        \end{minipage}
        \hspace{2.4cm}
        \begin{minipage}{0.25\linewidth}
            \centering
            \subfigure[\( y = 10\% \)]{
                \import{figures/appendix/parameter_search/heatmap-grid-relative-improvement_dp_10_tex}{heatmap-grid-relative-improvement_dp_10.tex}
            }
        \end{minipage}
        \hspace{1cm}
        \begin{minipage}{0.25\linewidth}
            \centering
            \subfigure[\( y = 20\% \)]{
                \import{figures/appendix/parameter_search/heatmap-grid-relative-improvement_dp_20_tex}{heatmap-grid-relative-improvement_dp_20.tex}
            }
        \end{minipage}
    }

    \caption{\raggedright Total delay reduction across five representative days and averaged values for different parameter configurations of the \gls{lns}. We vary the initial pool size \( x\% \), the sample size \( y\% \) drawn from the pool, and the maximum number of destroy-repair cycles \( z \) performed before invoking the local search operator for intensification.}
    \label{fig:parameter_search}
\end{figure}

To evaluate the impact of parameter settings on the performance of the \gls{lns} algorithm, we conduct a sensitivity analysis by varying three key parameters: the initial pool size $x$, the sample size $y$ drawn from the pool, and the number of destroy-repair cycles $z$ performed before invoking the local search for intensification. Figure~\ref{fig:parameter_search} presents heatmaps illustrating how different combinations of these parameters affect total delay reduction across five representative days—specifically, days 27 to 31 of the instances described in Appendix~\ref{appendix:instances_description}.
We also report the average delay reductions across these days. Each run is subject to a time limit of three hours. 

First, increasing the number of destroy-repair cycles $z$ tends to slow down the algorithm and may reduce its effectiveness, with the highest delay reductions often observed between $z = 2$ and $z = 4$. This suggests that it is more advantageous to intensify early using local search and then rebuild a new pool of trips, rather than repeatedly destroying and repairing the same set. Second, a moderate initial pool size, especially $x = 30\%$ or $40\%$, tends to give better results in most settings. Starting from $x = 50\%$, we generally observe a drop in performance. A larger pool appears beneficial as it includes both highly congested trips - prime candidates for improvement - and less congested ones, thereby supporting more effective exploration of the solution space.

Regarding the sample size $y$, the configuration with $y = 10\%$ offers the most robust performance, achieving high average delay reductions with relatively low variability across days. In contrast, sample sizes of $y = 5\%$ and $y = 20\%$ lead to greater variability and more frequent underperformance. This suggests that while the sampling procedure is helpful, once a promising pool of trips has been identified, it is more effective to focus on optimizing a moderately sized subset of trips rather than attempting to destroy and reinsert larger portions at once.

\section{Experimental setup for the \texttt{STAG} and \texttt{MATH} comparison}\label{appendix:comparison_with_math}

This appendix describes the experimental setup used to compare the staggering-only variant of our algorithm, \texttt{STAG}, with the matheuristic \texttt{MATH} developed for the same problem.

We run the experiments on the full Manhattan network, preprocessed with contraction hierarchies using a maximum shortcut length of 500 meters \citep[cf.][]{GeisbergerSandersEtAl2008}, and compute nominal travel times assuming a constant speed of 20~km/h.
We extract trip origins, destinations, and earliest departure times from the NYC taxi dataset collected in January 2015. We construct 31 instances, each consisting of 5000 trips sampled between 4:00 and 5:00~p.m. on a different day of the month.
As a baseline, we simulate an uncontrolled scenario in which each trip departs as early as possible and follows the shortest-distance route.
To model congestion, we define arc travel times using piecewise linear functions:
\begin{align*}
\tau_{\Arc}^{\Trip} =  \NominalTravelTime{\Arc} + \max_{k \in \NumPieces} \left\{ 0, \sum_{1 \leq k' \leq k} \PWLSlope{\Arc}{k'} \cdot \left(\Flow{\Arc}{\Trip} - \ThFlow{\Arc}{k'}\right) \right\},
\end{align*}
where $\NominalTravelTime{\Arc} > 0$, $\PWLSlope{\Arc}{k} > 0$, $\ThFlow{\Arc}{k} > 0$, and $\NumPieces = 3$ denotes the number of non-flat segments.
We control traffic intensity by restricting arc capacities to allow one \gls{amod} trip every 15 seconds. Accordingly, we set the first breakpoint to $\ThFlow{\Arc}{1} = \left\lceil \NominalTravelTime{\Arc} / 15 \right\rceil$, and define the remaining breakpoints as $\ThFlow{\Arc}{2} = 2 \cdot \ThFlow{\Arc}{1}$ and $\ThFlow{\Arc}{3} = 3 \cdot \ThFlow{\Arc}{1}$.
The corresponding slopes are:
\[
\PWLSlope{\Arc}{1} = \frac{0.5 \cdot \NominalTravelTime{\Arc}}{\ThFlow{\Arc}{1}}, \quad 
\PWLSlope{\Arc}{2} = \frac{\NominalTravelTime{\Arc}}{\ThFlow{\Arc}{1}}, \quad 
\PWLSlope{\Arc}{3} = \frac{1.5 \cdot \NominalTravelTime{\Arc}}{\ThFlow{\Arc}{1}}.
\]
Under the chosen parameter configuration, total delays in the baseline setting span from two to thirty hours, with a median of approximately six hours.
To explore the impact of departure-time flexibility, we generate two sets of instances — $\SigmaDown$ and $\SigmaUp$ — by setting the maximum allowable start time shift $\TripMaxStaggering{\Trip}$ to 10\% and 20\% of the nominal travel time along each trip’s assigned route, respectively.
The latest arrival time $\TripLatestArrival{\Trip}$ allows for a maximum travel time equal to 125\% of the free-flow time of the assigned route.

To conclude, we configure \texttt{STAG} with parameters $x = 50\%$, $y = 10\%$, and $z = 25$, chosen via a calibration procedure similar to that described in Appendix~\ref{appendix:parameters_search}. For \texttt{MATH}, we set a time limit of two hours.  {Further implementation details are available at \url{https://github.com/tumBAIS/integ_bal_stag}.}

%% file: figures/appendix/route_design/quantitative/user_equilibrium_total_delay_hours.tex
\begin{tikzpicture}

\definecolor{darkgray162159194}{RGB}{162,159,194}
\definecolor{darkgray176}{RGB}{176,176,176}
\definecolor{darkslateblue9364137}{RGB}{93,64,137}
\definecolor{darkslategray60}{RGB}{60,60,60}
\definecolor{lavender233232240}{RGB}{233,232,240}
\definecolor{lightgray201201221}{RGB}{201,201,221}
\definecolor{lightgray204}{RGB}{204,204,204}
\definecolor{lightslategray126118169}{RGB}{126,118,169}

\begin{axis}[
height=\PathAnalysisBoxplotsHeight,
legend cell align={left},
  legend style={
    fill opacity=0.8,
    draw opacity=1,
    text opacity=1,
    at={(1.01,1)},
    anchor=north west,
    draw=lightgray204,
  },
tick align=outside,
tick pos=left,
width=\PathAnalysisBoxplotsWidth,
x grid style={darkgray176},
xlabel={$\SimilarityThs$},
xmajorgrids,
xmin=-0.5, xmax=4.5,
xtick style={color=black},
xtick={0,1,2,3,4},
xticklabels={0.5,0.6,0.7,0.8,0.9},
y grid style={darkgray176},
ylabel={[hrs]},
ymajorgrids,
ymin=50, ymax=110,
ytick style={color=black}
]
\addlegendimage{empty legend}
\addlegendentry{\hspace{-.6cm}$k$}
\draw[draw=darkslategray60,fill=lavender233232240] (axis cs:0,0) rectangle (axis cs:0,0);

\addlegendimage{mark=square*,draw=darkslategray60,fill=lavender233232240}
\addlegendentry{3}

\draw[draw=darkslategray60,fill=lightgray201201221] (axis cs:0,0) rectangle (axis cs:0,0);
\addlegendimage{mark=square*,draw=darkslategray60,fill=lightgray201201221}
\addlegendentry{5}

\draw[draw=darkslategray60,fill=darkgray162159194] (axis cs:0,0) rectangle (axis cs:0,0);
\addlegendimage{mark=square*,draw=darkslategray60,fill=darkgray162159194}
\addlegendentry{7}

\draw[draw=darkslategray60,fill=lightslategray126118169] (axis cs:0,0) rectangle (axis cs:0,0);
\addlegendimage{mark=square*,draw=darkslategray60,fill=lightslategray126118169}
\addlegendentry{10}

\draw[draw=darkslategray60,fill=darkslateblue9364137] (axis cs:0,0) rectangle (axis cs:0,0);
\addlegendimage{mark=square*,draw=darkslategray60,fill=darkslateblue9364137}
\addlegendentry{15}

\addplot [darkslategray60, forget plot]
table {%
-0.2 69.6146754166667
-0.2 56.8571316666667
};
\addplot [darkslategray60, forget plot]
table {%
-0.2 90.6964313888889
-0.2 99.134602
};
\addplot [darkslategray60, forget plot]
table {%
-0.225 56.8571316666667
-0.175 56.8571316666667
};
\addplot [darkslategray60, forget plot]
table {%
-0.225 99.134602
-0.175 99.134602
};
\addplot [darkslategray60, forget plot]
table {%
0.8 68.8738959305556
0.8 56.6341188611111
};
\addplot [darkslategray60, forget plot]
table {%
0.8 88.7990773055556
0.8 98.5824380277778
};
\addplot [darkslategray60, forget plot]
table {%
0.775 56.6341188611111
0.825 56.6341188611111
};
\addplot [darkslategray60, forget plot]
table {%
0.775 98.5824380277778
0.825 98.5824380277778
};
\addplot [darkslategray60, forget plot]
table {%
1.8 68.0589218333333
1.8 56.3034505277778
};
\addplot [darkslategray60, forget plot]
table {%
1.8 87.5748437777778
1.8 97.0596931944444
};
\addplot [darkslategray60, forget plot]
table {%
1.775 56.3034505277778
1.825 56.3034505277778
};
\addplot [darkslategray60, forget plot]
table {%
1.775 97.0596931944444
1.825 97.0596931944444
};
\addplot [darkslategray60, forget plot]
table {%
2.8 68.0886926388889
2.8 56.39787125
};
\addplot [darkslategray60, forget plot]
table {%
2.8 87.5580899027778
2.8 96.9981003333333
};
\addplot [darkslategray60, forget plot]
table {%
2.775 56.39787125
2.825 56.39787125
};
\addplot [darkslategray60, forget plot]
table {%
2.775 96.9981003333333
2.825 96.9981003333333
};
\addplot [darkslategray60, forget plot]
table {%
3.8 69.1732885972222
3.8 56.573854
};
\addplot [darkslategray60, forget plot]
table {%
3.8 92.3959620833333
3.8 105.588286277778
};
\addplot [darkslategray60, forget plot]
table {%
3.775 56.573854
3.825 56.573854
};
\addplot [darkslategray60, forget plot]
table {%
3.775 105.588286277778
3.825 105.588286277778
};
\addplot [darkslategray60, forget plot]
table {%
-0.1 69.4826950555556
-0.1 56.807397
};
\addplot [darkslategray60, forget plot]
table {%
-0.1 90.0857726388889
-0.1 99.0350186944444
};
\addplot [darkslategray60, forget plot]
table {%
-0.125 56.807397
-0.075 56.807397
};
\addplot [darkslategray60, forget plot]
table {%
-0.125 99.0350186944444
-0.075 99.0350186944444
};
\addplot [darkslategray60, forget plot]
table {%
0.9 68.59195825
0.9 56.5052017222222
};
\addplot [darkslategray60, forget plot]
table {%
0.9 88.4159416666667
0.9 97.1821825833333
};
\addplot [darkslategray60, forget plot]
table {%
0.875 56.5052017222222
0.925 56.5052017222222
};
\addplot [darkslategray60, forget plot]
table {%
0.875 97.1821825833333
0.925 97.1821825833333
};
\addplot [darkslategray60, forget plot]
table {%
1.9 67.3736641805556
1.9 56.0078345
};
\addplot [darkslategray60, forget plot]
table {%
1.9 86.7222888194445
1.9 95.6887335833333
};
\addplot [darkslategray60, forget plot]
table {%
1.875 56.0078345
1.925 56.0078345
};
\addplot [darkslategray60, forget plot]
table {%
1.875 95.6887335833333
1.925 95.6887335833333
};
\addplot [darkslategray60, forget plot]
table {%
2.9 65.9507895
2.9 54.9511013333333
};
\addplot [darkslategray60, forget plot]
table {%
2.9 84.0788508472222
2.9 92.6347391944444
};
\addplot [darkslategray60, forget plot]
table {%
2.875 54.9511013333333
2.925 54.9511013333333
};
\addplot [darkslategray60, forget plot]
table {%
2.875 92.6347391944444
2.925 92.6347391944444
};
\addplot [darkslategray60, forget plot]
table {%
3.9 65.6964739305556
3.9 54.4586294722222
};
\addplot [darkslategray60, forget plot]
table {%
3.9 84.4441539583333
3.9 93.7082342222222
};
\addplot [darkslategray60, forget plot]
table {%
3.875 54.4586294722222
3.925 54.4586294722222
};
\addplot [darkslategray60, forget plot]
table {%
3.875 93.7082342222222
3.925 93.7082342222222
};
\addplot [darkslategray60, forget plot]
table {%
0 69.4826950555556
0 56.807397
};
\addplot [darkslategray60, forget plot]
table {%
0 90.0857726388889
0 99.0350186944444
};
\addplot [darkslategray60, forget plot]
table {%
-0.025 56.807397
0.025 56.807397
};
\addplot [darkslategray60, forget plot]
table {%
-0.025 99.0350186944444
0.025 99.0350186944444
};
\addplot [darkslategray60, forget plot]
table {%
1 68.5784814861111
1 56.5052017222222
};
\addplot [darkslategray60, forget plot]
table {%
1 88.3475236805556
1 97.1590301111111
};
\addplot [darkslategray60, forget plot]
table {%
0.975 56.5052017222222
1.025 56.5052017222222
};
\addplot [darkslategray60, forget plot]
table {%
0.975 97.1590301111111
1.025 97.1590301111111
};
\addplot [darkslategray60, forget plot]
table {%
2 67.2746969027778
2 55.91310675
};
\addplot [darkslategray60, forget plot]
table {%
2 86.4452754861111
2 95.43748075
};
\addplot [darkslategray60, forget plot]
table {%
1.975 55.91310675
2.025 55.91310675
};
\addplot [darkslategray60, forget plot]
table {%
1.975 95.43748075
2.025 95.43748075
};
\addplot [darkslategray60, forget plot]
table {%
3 65.9648401111111
3 54.9335945555556
};
\addplot [darkslategray60, forget plot]
table {%
3 83.8914472916667
3 92.4399796111111
};
\addplot [darkslategray60, forget plot]
table {%
2.975 54.9335945555556
3.025 54.9335945555556
};
\addplot [darkslategray60, forget plot]
table {%
2.975 92.4399796111111
3.025 92.4399796111111
};
\addplot [darkslategray60, forget plot]
table {%
4 64.8416863333333
4 54.1047118333333
};
\addplot [darkslategray60, forget plot]
table {%
4 82.2755946111111
4 91.4341448055555
};
\addplot [darkslategray60, forget plot]
table {%
3.975 54.1047118333333
4.025 54.1047118333333
};
\addplot [darkslategray60, forget plot]
table {%
3.975 91.4341448055555
4.025 91.4341448055555
};
\addplot [darkslategray60, forget plot]
table {%
0.1 69.4826950555556
0.1 56.807397
};
\addplot [darkslategray60, forget plot]
table {%
0.1 90.0857726388889
0.1 99.0350186944444
};
\addplot [darkslategray60, forget plot]
table {%
0.075 56.807397
0.125 56.807397
};
\addplot [darkslategray60, forget plot]
table {%
0.075 99.0350186944444
0.125 99.0350186944444
};
\addplot [darkslategray60, forget plot]
table {%
1.1 68.5784814861111
1.1 56.5052017222222
};
\addplot [darkslategray60, forget plot]
table {%
1.1 88.3475236805556
1.1 97.1590301111111
};
\addplot [darkslategray60, forget plot]
table {%
1.075 56.5052017222222
1.125 56.5052017222222
};
\addplot [darkslategray60, forget plot]
table {%
1.075 97.1590301111111
1.125 97.1590301111111
};
\addplot [darkslategray60, forget plot]
table {%
2.1 67.2746969027778
2.1 55.91310675
};
\addplot [darkslategray60, forget plot]
table {%
2.1 86.4452754861111
2.1 95.43748075
};
\addplot [darkslategray60, forget plot]
table {%
2.075 55.91310675
2.125 55.91310675
};
\addplot [darkslategray60, forget plot]
table {%
2.075 95.43748075
2.125 95.43748075
};
\addplot [darkslategray60, forget plot]
table {%
3.1 65.9640744027778
3.1 54.92808
};
\addplot [darkslategray60, forget plot]
table {%
3.1 83.85416125
3.1 92.4399796111111
};
\addplot [darkslategray60, forget plot]
table {%
3.075 54.92808
3.125 54.92808
};
\addplot [darkslategray60, forget plot]
table {%
3.075 92.4399796111111
3.125 92.4399796111111
};
\addplot [darkslategray60, forget plot]
table {%
4.1 64.4875244305556
4.1 53.8785926944444
};
\addplot [darkslategray60, forget plot]
table {%
4.1 82.1692360277778
4.1 91.0777395833333
};
\addplot [darkslategray60, forget plot]
table {%
4.075 53.8785926944444
4.125 53.8785926944444
};
\addplot [darkslategray60, forget plot]
table {%
4.075 91.0777395833333
4.125 91.0777395833333
};
\addplot [darkslategray60, forget plot]
table {%
0.2 69.4826950555556
0.2 56.807397
};
\addplot [darkslategray60, forget plot]
table {%
0.2 90.0857726388889
0.2 99.0350186944444
};
\addplot [darkslategray60, forget plot]
table {%
0.175 56.807397
0.225 56.807397
};
\addplot [darkslategray60, forget plot]
table {%
0.175 99.0350186944444
0.225 99.0350186944444
};
\addplot [darkslategray60, forget plot]
table {%
1.2 68.5784814861111
1.2 56.5052017222222
};
\addplot [darkslategray60, forget plot]
table {%
1.2 88.3475236805556
1.2 97.1590301111111
};
\addplot [darkslategray60, forget plot]
table {%
1.175 56.5052017222222
1.225 56.5052017222222
};
\addplot [darkslategray60, forget plot]
table {%
1.175 97.1590301111111
1.225 97.1590301111111
};
\addplot [darkslategray60, forget plot]
table {%
2.2 67.2746969027778
2.2 55.91310675
};
\addplot [darkslategray60, forget plot]
table {%
2.2 86.4452754861111
2.2 95.43748075
};
\addplot [darkslategray60, forget plot]
table {%
2.175 55.91310675
2.225 55.91310675
};
\addplot [darkslategray60, forget plot]
table {%
2.175 95.43748075
2.225 95.43748075
};
\addplot [darkslategray60, forget plot]
table {%
3.2 65.9640744027778
3.2 54.92808
};
\addplot [darkslategray60, forget plot]
table {%
3.2 83.85416125
3.2 92.4399796111111
};
\addplot [darkslategray60, forget plot]
table {%
3.175 54.92808
3.225 54.92808
};
\addplot [darkslategray60, forget plot]
table {%
3.175 92.4399796111111
3.225 92.4399796111111
};
\addplot [darkslategray60, forget plot]
table {%
4.2 64.5311711111111
4.2 53.8663811666667
};
\addplot [darkslategray60, forget plot]
table {%
4.2 82.1123924722222
4.2 91.0516282222222
};
\addplot [darkslategray60, forget plot]
table {%
4.175 53.8663811666667
4.225 53.8663811666667
};
\addplot [darkslategray60, forget plot]
table {%
4.175 91.0516282222222
4.225 91.0516282222222
};
\path [draw=darkslategray60, fill=lavender233232240]
(axis cs:-0.25,69.6146754166667)
--(axis cs:-0.15,69.6146754166667)
--(axis cs:-0.15,90.6964313888889)
--(axis cs:-0.25,90.6964313888889)
--(axis cs:-0.25,69.6146754166667)
--cycle;
\path [draw=darkslategray60, fill=lavender233232240]
(axis cs:0.75,68.8738959305556)
--(axis cs:0.85,68.8738959305556)
--(axis cs:0.85,88.7990773055556)
--(axis cs:0.75,88.7990773055556)
--(axis cs:0.75,68.8738959305556)
--cycle;
\path [draw=darkslategray60, fill=lavender233232240]
(axis cs:1.75,68.0589218333333)
--(axis cs:1.85,68.0589218333333)
--(axis cs:1.85,87.5748437777778)
--(axis cs:1.75,87.5748437777778)
--(axis cs:1.75,68.0589218333333)
--cycle;
\path [draw=darkslategray60, fill=lavender233232240]
(axis cs:2.75,68.0886926388889)
--(axis cs:2.85,68.0886926388889)
--(axis cs:2.85,87.5580899027778)
--(axis cs:2.75,87.5580899027778)
--(axis cs:2.75,68.0886926388889)
--cycle;
\path [draw=darkslategray60, fill=lavender233232240]
(axis cs:3.75,69.1732885972222)
--(axis cs:3.85,69.1732885972222)
--(axis cs:3.85,92.3959620833333)
--(axis cs:3.75,92.3959620833333)
--(axis cs:3.75,69.1732885972222)
--cycle;
\path [draw=darkslategray60, fill=lightgray201201221]
(axis cs:-0.15,69.4826950555556)
--(axis cs:-0.05,69.4826950555556)
--(axis cs:-0.05,90.0857726388889)
--(axis cs:-0.15,90.0857726388889)
--(axis cs:-0.15,69.4826950555556)
--cycle;
\path [draw=darkslategray60, fill=lightgray201201221]
(axis cs:0.85,68.59195825)
--(axis cs:0.95,68.59195825)
--(axis cs:0.95,88.4159416666667)
--(axis cs:0.85,88.4159416666667)
--(axis cs:0.85,68.59195825)
--cycle;
\path [draw=darkslategray60, fill=lightgray201201221]
(axis cs:1.85,67.3736641805556)
--(axis cs:1.95,67.3736641805556)
--(axis cs:1.95,86.7222888194445)
--(axis cs:1.85,86.7222888194445)
--(axis cs:1.85,67.3736641805556)
--cycle;
\path [draw=darkslategray60, fill=lightgray201201221]
(axis cs:2.85,65.9507895)
--(axis cs:2.95,65.9507895)
--(axis cs:2.95,84.0788508472222)
--(axis cs:2.85,84.0788508472222)
--(axis cs:2.85,65.9507895)
--cycle;
\path [draw=darkslategray60, fill=lightgray201201221]
(axis cs:3.85,65.6964739305556)
--(axis cs:3.95,65.6964739305556)
--(axis cs:3.95,84.4441539583333)
--(axis cs:3.85,84.4441539583333)
--(axis cs:3.85,65.6964739305556)
--cycle;
\path [draw=darkslategray60, fill=darkgray162159194]
(axis cs:-0.05,69.4826950555556)
--(axis cs:0.05,69.4826950555556)
--(axis cs:0.05,90.0857726388889)
--(axis cs:-0.05,90.0857726388889)
--(axis cs:-0.05,69.4826950555556)
--cycle;
\path [draw=darkslategray60, fill=darkgray162159194]
(axis cs:0.95,68.5784814861111)
--(axis cs:1.05,68.5784814861111)
--(axis cs:1.05,88.3475236805556)
--(axis cs:0.95,88.3475236805556)
--(axis cs:0.95,68.5784814861111)
--cycle;
\path [draw=darkslategray60, fill=darkgray162159194]
(axis cs:1.95,67.2746969027778)
--(axis cs:2.05,67.2746969027778)
--(axis cs:2.05,86.4452754861111)
--(axis cs:1.95,86.4452754861111)
--(axis cs:1.95,67.2746969027778)
--cycle;
\path [draw=darkslategray60, fill=darkgray162159194]
(axis cs:2.95,65.9648401111111)
--(axis cs:3.05,65.9648401111111)
--(axis cs:3.05,83.8914472916667)
--(axis cs:2.95,83.8914472916667)
--(axis cs:2.95,65.9648401111111)
--cycle;
\path [draw=darkslategray60, fill=darkgray162159194]
(axis cs:3.95,64.8416863333333)
--(axis cs:4.05,64.8416863333333)
--(axis cs:4.05,82.2755946111111)
--(axis cs:3.95,82.2755946111111)
--(axis cs:3.95,64.8416863333333)
--cycle;
\path [draw=darkslategray60, fill=lightslategray126118169]
(axis cs:0.05,69.4826950555556)
--(axis cs:0.15,69.4826950555556)
--(axis cs:0.15,90.0857726388889)
--(axis cs:0.05,90.0857726388889)
--(axis cs:0.05,69.4826950555556)
--cycle;
\path [draw=darkslategray60, fill=lightslategray126118169]
(axis cs:1.05,68.5784814861111)
--(axis cs:1.15,68.5784814861111)
--(axis cs:1.15,88.3475236805556)
--(axis cs:1.05,88.3475236805556)
--(axis cs:1.05,68.5784814861111)
--cycle;
\path [draw=darkslategray60, fill=lightslategray126118169]
(axis cs:2.05,67.2746969027778)
--(axis cs:2.15,67.2746969027778)
--(axis cs:2.15,86.4452754861111)
--(axis cs:2.05,86.4452754861111)
--(axis cs:2.05,67.2746969027778)
--cycle;
\path [draw=darkslategray60, fill=lightslategray126118169]
(axis cs:3.05,65.9640744027778)
--(axis cs:3.15,65.9640744027778)
--(axis cs:3.15,83.85416125)
--(axis cs:3.05,83.85416125)
--(axis cs:3.05,65.9640744027778)
--cycle;
\path [draw=darkslategray60, fill=lightslategray126118169]
(axis cs:4.05,64.4875244305556)
--(axis cs:4.15,64.4875244305556)
--(axis cs:4.15,82.1692360277778)
--(axis cs:4.05,82.1692360277778)
--(axis cs:4.05,64.4875244305556)
--cycle;
\path [draw=darkslategray60, fill=darkslateblue9364137]
(axis cs:0.15,69.4826950555556)
--(axis cs:0.25,69.4826950555556)
--(axis cs:0.25,90.0857726388889)
--(axis cs:0.15,90.0857726388889)
--(axis cs:0.15,69.4826950555556)
--cycle;
\path [draw=darkslategray60, fill=darkslateblue9364137]
(axis cs:1.15,68.5784814861111)
--(axis cs:1.25,68.5784814861111)
--(axis cs:1.25,88.3475236805556)
--(axis cs:1.15,88.3475236805556)
--(axis cs:1.15,68.5784814861111)
--cycle;
\path [draw=darkslategray60, fill=darkslateblue9364137]
(axis cs:2.15,67.2746969027778)
--(axis cs:2.25,67.2746969027778)
--(axis cs:2.25,86.4452754861111)
--(axis cs:2.15,86.4452754861111)
--(axis cs:2.15,67.2746969027778)
--cycle;
\path [draw=darkslategray60, fill=darkslateblue9364137]
(axis cs:3.15,65.9640744027778)
--(axis cs:3.25,65.9640744027778)
--(axis cs:3.25,83.85416125)
--(axis cs:3.15,83.85416125)
--(axis cs:3.15,65.9640744027778)
--cycle;
\path [draw=darkslategray60, fill=darkslateblue9364137]
(axis cs:4.15,64.5311711111111)
--(axis cs:4.25,64.5311711111111)
--(axis cs:4.25,82.1123924722222)
--(axis cs:4.15,82.1123924722222)
--(axis cs:4.15,64.5311711111111)
--cycle;
\addplot [darkslategray60, forget plot]
table {%
-0.25 83.8097275833333
-0.15 83.8097275833333
};
\addplot [darkslategray60, forget plot]
table {%
0.75 82.4689609444444
0.85 82.4689609444444
};
\addplot [darkslategray60, forget plot]
table {%
1.75 82.0125673611111
1.85 82.0125673611111
};
\addplot [darkslategray60, forget plot]
table {%
2.75 81.3708856111111
2.85 81.3708856111111
};
\addplot [darkslategray60, forget plot]
table {%
3.75 85.01331675
3.85 85.01331675
};
\addplot [darkslategray60, forget plot]
table {%
-0.15 83.98777775
-0.05 83.98777775
};
\addplot [darkslategray60, forget plot]
table {%
0.85 82.3642316944444
0.95 82.3642316944444
};
\addplot [darkslategray60, forget plot]
table {%
1.85 80.7751995833333
1.95 80.7751995833333
};
\addplot [darkslategray60, forget plot]
table {%
2.85 79.2172404166667
2.95 79.2172404166667
};
\addplot [darkslategray60, forget plot]
table {%
3.85 79.2786439722222
3.95 79.2786439722222
};
\addplot [darkslategray60, forget plot]
table {%
-0.05 83.98777775
0.05 83.98777775
};
\addplot [darkslategray60, forget plot]
table {%
0.95 82.3642316944444
1.05 82.3642316944444
};
\addplot [darkslategray60, forget plot]
table {%
1.95 80.694329
2.05 80.694329
};
\addplot [darkslategray60, forget plot]
table {%
2.95 78.8068023611111
3.05 78.8068023611111
};
\addplot [darkslategray60, forget plot]
table {%
3.95 78.0154178055556
4.05 78.0154178055556
};
\addplot [darkslategray60, forget plot]
table {%
0.05 83.98777775
0.15 83.98777775
};
\addplot [darkslategray60, forget plot]
table {%
1.05 82.3642316944444
1.15 82.3642316944444
};
\addplot [darkslategray60, forget plot]
table {%
2.05 80.694329
2.15 80.694329
};
\addplot [darkslategray60, forget plot]
table {%
3.05 78.8518369444444
3.15 78.8518369444444
};
\addplot [darkslategray60, forget plot]
table {%
4.05 77.6466328611111
4.15 77.6466328611111
};
\addplot [darkslategray60, forget plot]
table {%
0.15 83.98777775
0.25 83.98777775
};
\addplot [darkslategray60, forget plot]
table {%
1.15 82.3642316944444
1.25 82.3642316944444
};
\addplot [darkslategray60, forget plot]
table {%
2.15 80.694329
2.25 80.694329
};
\addplot [darkslategray60, forget plot]
table {%
3.15 78.8518369444444
3.25 78.8518369444444
};
\addplot [darkslategray60, forget plot]
table {%
4.15 77.5234067777778
4.25 77.5234067777778
};
\end{axis}

\end{tikzpicture}

%% file: tables/instance_summary.tex
\begin{table}[!htbp]
\centering
\renewcommand{\arraystretch}{1.15} 
\caption{\raggedright Summary of \gls{rduo} characteristics for the instances considered. We report the day, the number of trips \( |\SetTrips| \), the total travel time \( \TotalTravelTime{} \), the total delay \( D \), their ratio, and the average, median, and maximum values of the total delay per trip \( \TotalDelayTrip{\Trip} \). We also include the corresponding statistics for the ratio between each trip's travel time and its free-flow travel time (with a value of 1 indicating free-flow conditions).}
\label{table:instance_summary}
\begin{tabularx}{\textwidth}{>{\centering\arraybackslash}X >{\centering\arraybackslash}X >{\centering\arraybackslash}X >{\centering\arraybackslash}X >{\centering\arraybackslash}X >{\centering\arraybackslash}X >{\centering\arraybackslash}X >{\centering\arraybackslash}X >{\centering\arraybackslash}X >{\centering\arraybackslash}X >{\centering\arraybackslash}X >{\centering\arraybackslash}X}
\toprule
Day & \( |\SetTrips| \) & \( \TotalTravelTime{} \) [hrs] & \( D \) [hrs] & \( \TotalTravelTime{}/D \) [\%] & \multicolumn{3}{c}{\( \TotalDelayTrip{\Trip} \) [min]} & \multicolumn{3}{c}{\( \TravelTimeRoute{\Trip}/\FreeFlowRoute{\Trip} \)} \\
\cmidrule(lr){6-8} \cmidrule(lr){9-11}
& & & & & Avg & Med & Max & Avg & Med & Max \\
\midrule
1 & 5283 & 561 & 54 & \textcolor{white}{0}9.63 & 0.65 & 0.55 & 2.90 & 1.11 & 1.09 & 1.62 \\
2 & 5428 & 578 & 60 & 10.38 & 0.70 & 0.56 & 2.77 & 1.12 & 1.10 & 1.59 \\
3 & 5557 & 598 & 64 & 10.70 & 0.73 & 0.59 & 3.60 & 1.12 & 1.10 & 1.68 \\
4 & 5448 & 575 & 56 & \textcolor{white}{0}9.74 & 0.65 & 0.55 & 2.80 & 1.11 & 1.09 & 1.55 \\
5 & 5399 & 575 & 57 & \textcolor{white}{0}9.91 & 0.67 & 0.58 & 2.78 & 1.11 & 1.10 & 1.54 \\
6 & 6309 & 689 & 81 & 11.76 & 0.81 & 0.69 & 3.18 & 1.14 & 1.11 & 1.70 \\
7 & 5857 & 626 & 70 & 11.18 & 0.76 & 0.65 & 3.04 & 1.13 & 1.11 & 1.68 \\
8 & 5937 & 635 & 70 & 11.02 & 0.75 & 0.63 & 3.00 & 1.13 & 1.11 & 1.61 \\
9 & 6325 & 690 & 83 & 12.03 & 0.83 & 0.72 & 3.45 & 1.14 & 1.12 & 1.70 \\
10 & 5807 & 628 & 67 & 10.67 & 0.73 & 0.62 & 2.97 & 1.12 & 1.10 & 1.83 \\
11 & 5574 & 599 & 62 & 10.35 & 0.71 & 0.60 & 2.96 & 1.12 & 1.10 & 1.50 \\
12 & 6142 & 656 & 72 & 10.98 & 0.75 & 0.64 & 3.12 & 1.13 & 1.11 & 1.57 \\
13 & 6572 & 715 & 88 & 12.31 & 0.85 & 0.74 & 3.14 & 1.14 & 1.12 & 1.71 \\
14 & 6275 & 675 & 78 & 11.56 & 0.79 & 0.67 & 3.52 & 1.13 & 1.11 & 1.82 \\
15 & 6281 & 679 & 81 & 11.93 & 0.82 & 0.69 & 3.33 & 1.14 & 1.12 & 1.83 \\
16 & 6627 & 725 & 90 & 12.41 & 0.86 & 0.73 & 3.64 & 1.14 & 1.12 & 1.76 \\
17 & 6344 & 688 & 83 & 12.06 & 0.84 & 0.68 & 3.37 & 1.14 & 1.12 & 1.85 \\
18 & 5747 & 624 & 69 & 11.06 & 0.75 & 0.63 & 3.04 & 1.12 & 1.11 & 1.61 \\
19 & 5611 & 603 & 63 & 10.45 & 0.71 & 0.61 & 3.59 & 1.12 & 1.10 & 1.59 \\
20 & 6478 & 708 & 85 & 12.01 & 0.84 & 0.71 & 3.09 & 1.14 & 1.12 & 1.87 \\
21 & 6622 & 722 & 90 & 12.47 & 0.87 & 0.75 & 3.54 & 1.15 & 1.13 & 1.78 \\
22 & 6194 & 673 & 76 & 11.29 & 0.78 & 0.68 & 3.02 & 1.13 & 1.11 & 1.65 \\
23 & 6471 & 707 & 88 & 12.45 & 0.86 & 0.72 & 3.38 & 1.14 & 1.12 & 1.82 \\
24 & 6109 & 668 & 78 & 11.68 & 0.81 & 0.67 & 3.31 & 1.13 & 1.11 & 1.67 \\
25 & 5355 & 573 & 58 & 10.12 & 0.68 & 0.56 & 3.22 & 1.11 & 1.10 & 1.54 \\
26 & 6240 & 673 & 79 & 11.74 & 0.80 & 0.67 & 3.56 & 1.14 & 1.11 & 1.96 \\
27 & 6593 & 722 & 88 & 12.19 & 0.85 & 0.72 & 3.56 & 1.14 & 1.12 & 1.71 \\
28 & 6443 & 698 & 84 & 12.03 & 0.84 & 0.70 & 3.59 & 1.14 & 1.12 & 1.72 \\
29 & 6353 & 680 & 79 & 11.62 & 0.79 & 0.67 & 3.07 & 1.14 & 1.11 & 1.64 \\
30 & 6715 & 736 & 92 & 12.50 & 0.87 & 0.74 & 3.82 & 1.15 & 1.13 & 1.90 \\
31 & 6144 & 680 & 81 & 11.91 & 0.84 & 0.68 & 3.60 & 1.14 & 1.11 & 1.64 \\
\midrule
Avg & 6072 & 657 & 76 & 11.57 & 0.78 & 0.66 & 3.26 & 1.13 & 1.11 & 1.70 \\
\bottomrule
\end{tabularx}
\end{table}

%% file: figures/appendix/instances_description/travel_times_histogram.tex
\begin{tikzpicture}

\definecolor{darkgray176}{RGB}{176,176,176}
\definecolor{lightgreen}{RGB}{144,238,144}

\begin{axis}[
height=\InstancesHistogramsHeight,
minor xtick={},
minor ytick={},
scaled y ticks=true,
tick align=outside,
tick pos=left,
width=\InstancesHistogramsWidth,
x grid style={darkgray176},
xlabel={[min]},
xmin=-0.9, xmax=19.58665,
xtick style={color=black},
xtick={-5,0,5,10,15,20},
xticklabel style={/pgf/number format/fixed, /pgf/number format/precision=0},
y grid style={darkgray176},
y tick label style={/pgf/number format/sci, /pgf/number format/sci zerofill=false},
ymin=0, ymax=173968.2,
ytick style={color=black},
ytick={0,40000,80000,120000,160000,200000},
yticklabel style={/pgf/number format/sci, /pgf/number format/precision=0},
yticklabels={0.0,0.4,0.8,1.2,1.6,2.0}
]
\draw[draw=black,fill=lightgreen,opacity=0.7,line width=0.48pt] (axis cs:1.6635,0) rectangle (axis cs:2.51698333333333,548);
\draw[draw=black,fill=lightgreen,opacity=0.7,line width=0.48pt] (axis cs:2.51698333333333,0) rectangle (axis cs:3.37046666666667,8988);
\draw[draw=black,fill=lightgreen,opacity=0.7,line width=0.48pt] (axis cs:3.37046666666667,0) rectangle (axis cs:4.22395,50786);
\draw[draw=black,fill=lightgreen,opacity=0.7,line width=0.48pt] (axis cs:4.22395,0) rectangle (axis cs:5.07743333333333,122151);
\draw[draw=black,fill=lightgreen,opacity=0.7,line width=0.48pt] (axis cs:5.07743333333333,0) rectangle (axis cs:5.93091666666667,164063);
\draw[draw=black,fill=lightgreen,opacity=0.7,line width=0.48pt] (axis cs:5.93091666666667,0) rectangle (axis cs:6.7844,165684);
\draw[draw=black,fill=lightgreen,opacity=0.7,line width=0.48pt] (axis cs:6.7844,0) rectangle (axis cs:7.63788333333333,123289);
\draw[draw=black,fill=lightgreen,opacity=0.7,line width=0.48pt] (axis cs:7.63788333333333,0) rectangle (axis cs:8.49136666666667,72168);
\draw[draw=black,fill=lightgreen,opacity=0.7,line width=0.48pt] (axis cs:8.49136666666667,0) rectangle (axis cs:9.34485,48593);
\draw[draw=black,fill=lightgreen,opacity=0.7,line width=0.48pt] (axis cs:9.34485,0) rectangle (axis cs:10.1983333333333,29223);
\draw[draw=black,fill=lightgreen,opacity=0.7,line width=0.48pt] (axis cs:10.1983333333333,0) rectangle (axis cs:11.0518166666667,17886);
\draw[draw=black,fill=lightgreen,opacity=0.7,line width=0.48pt] (axis cs:11.0518166666667,0) rectangle (axis cs:11.9053,11995);
\draw[draw=black,fill=lightgreen,opacity=0.7,line width=0.48pt] (axis cs:11.9053,0) rectangle (axis cs:12.7587833333333,6883);
\draw[draw=black,fill=lightgreen,opacity=0.7,line width=0.48pt] (axis cs:12.7587833333333,0) rectangle (axis cs:13.6122666666667,4578);
\draw[draw=black,fill=lightgreen,opacity=0.7,line width=0.48pt] (axis cs:13.6122666666667,0) rectangle (axis cs:14.46575,2518);
\draw[draw=black,fill=lightgreen,opacity=0.7,line width=0.48pt] (axis cs:14.46575,0) rectangle (axis cs:15.3192333333333,1556);
\draw[draw=black,fill=lightgreen,opacity=0.7,line width=0.48pt] (axis cs:15.3192333333333,0) rectangle (axis cs:16.1727166666667,524);
\draw[draw=black,fill=lightgreen,opacity=0.7,line width=0.48pt] (axis cs:16.1727166666667,0) rectangle (axis cs:17.0262,75);
\draw[draw=black,fill=lightgreen,opacity=0.7,line width=0.48pt] (axis cs:17.0262,0) rectangle (axis cs:17.8796833333333,9);
\draw[draw=black,fill=lightgreen,opacity=0.7,line width=0.48pt] (axis cs:17.8796833333333,0) rectangle (axis cs:18.7331666666667,1);
\end{axis}

\end{tikzpicture}

%% file: figures/appendix/instances_description/relative_max_deviation_histogram.tex
\begin{tikzpicture}

\definecolor{darkgray176}{RGB}{176,176,176}
\definecolor{orange}{RGB}{255,165,0}

\begin{axis}[
height=\InstancesHistogramsHeight,
minor xtick={},
minor ytick={},
scaled y ticks=true,
tick align=outside,
tick pos=left,
width=\InstancesHistogramsWidth,
x grid style={darkgray176},
xlabel={[\%]},
xmin=-0.9, xmax=100,
xtick style={color=black},
xtick={-25,0,25,50,75,100},
xticklabel style={/pgf/number format/fixed, /pgf/number format/precision=0},
y grid style={darkgray176},
y tick label style={/pgf/number format/sci, /pgf/number format/sci zerofill=false},
ymin=0, ymax=19093.2,
ytick style={color=black},
ytick={0,4000,8000,12000,16000,20000},
yticklabel style={/pgf/number format/sci, /pgf/number format/precision=0},
yticklabels={0.0,0.4,0.8,1.2,1.6,2.0}
]
\draw[draw=black,fill=orange,opacity=0.7,line width=0.48pt] (axis cs:0.04,0) rectangle (axis cs:3.1192,1374);
\draw[draw=black,fill=orange,opacity=0.7,line width=0.48pt] (axis cs:3.1192,0) rectangle (axis cs:6.1984,5475);
\draw[draw=black,fill=orange,opacity=0.7,line width=0.48pt] (axis cs:6.1984,0) rectangle (axis cs:9.2776,9828);
\draw[draw=black,fill=orange,opacity=0.7,line width=0.48pt] (axis cs:9.2776,0) rectangle (axis cs:12.3568,11209);
\draw[draw=black,fill=orange,opacity=0.7,line width=0.48pt] (axis cs:12.3568,0) rectangle (axis cs:15.436,15705);
\draw[draw=black,fill=orange,opacity=0.7,line width=0.48pt] (axis cs:15.436,0) rectangle (axis cs:18.5152,17725);
\draw[draw=black,fill=orange,opacity=0.7,line width=0.48pt] (axis cs:18.5152,0) rectangle (axis cs:21.5944,18184);
\draw[draw=black,fill=orange,opacity=0.7,line width=0.48pt] (axis cs:21.5944,0) rectangle (axis cs:24.6736,17277);
\draw[draw=black,fill=orange,opacity=0.7,line width=0.48pt] (axis cs:24.6736,0) rectangle (axis cs:27.7528,16855);
\draw[draw=black,fill=orange,opacity=0.7,line width=0.48pt] (axis cs:27.7528,0) rectangle (axis cs:30.832,13569);
\draw[draw=black,fill=orange,opacity=0.7,line width=0.48pt] (axis cs:30.832,0) rectangle (axis cs:33.9112,11996);
\draw[draw=black,fill=orange,opacity=0.7,line width=0.48pt] (axis cs:33.9112,0) rectangle (axis cs:36.9904,9353);
\draw[draw=black,fill=orange,opacity=0.7,line width=0.48pt] (axis cs:36.9904,0) rectangle (axis cs:40.0696,7698);
\draw[draw=black,fill=orange,opacity=0.7,line width=0.48pt] (axis cs:40.0696,0) rectangle (axis cs:43.1488,6403);
\draw[draw=black,fill=orange,opacity=0.7,line width=0.48pt] (axis cs:43.1488,0) rectangle (axis cs:46.228,5209);
\draw[draw=black,fill=orange,opacity=0.7,line width=0.48pt] (axis cs:46.228,0) rectangle (axis cs:49.3072,4272);
\draw[draw=black,fill=orange,opacity=0.7,line width=0.48pt] (axis cs:49.3072,0) rectangle (axis cs:52.3864,3189);
\draw[draw=black,fill=orange,opacity=0.7,line width=0.48pt] (axis cs:52.3864,0) rectangle (axis cs:55.4656,2331);
\draw[draw=black,fill=orange,opacity=0.7,line width=0.48pt] (axis cs:55.4656,0) rectangle (axis cs:58.5448,1765);
\draw[draw=black,fill=orange,opacity=0.7,line width=0.48pt] (axis cs:58.5448,0) rectangle (axis cs:61.624,1414);
\draw[draw=black,fill=orange,opacity=0.7,line width=0.48pt] (axis cs:61.624,0) rectangle (axis cs:64.7032,880);
\draw[draw=black,fill=orange,opacity=0.7,line width=0.48pt] (axis cs:64.7032,0) rectangle (axis cs:67.7824,762);
\draw[draw=black,fill=orange,opacity=0.7,line width=0.48pt] (axis cs:67.7824,0) rectangle (axis cs:70.8616,564);
\draw[draw=black,fill=orange,opacity=0.7,line width=0.48pt] (axis cs:70.8616,0) rectangle (axis cs:73.9408,398);
\draw[draw=black,fill=orange,opacity=0.7,line width=0.48pt] (axis cs:73.9408,0) rectangle (axis cs:77.02,289);
\draw[draw=black,fill=orange,opacity=0.7,line width=0.48pt] (axis cs:77.02,0) rectangle (axis cs:80.0992,303);
\draw[draw=black,fill=orange,opacity=0.7,line width=0.48pt] (axis cs:80.0992,0) rectangle (axis cs:83.1784,165);
\draw[draw=black,fill=orange,opacity=0.7,line width=0.48pt] (axis cs:83.1784,0) rectangle (axis cs:86.2576,109);
\draw[draw=black,fill=orange,opacity=0.7,line width=0.48pt] (axis cs:86.2576,0) rectangle (axis cs:89.3368,95);
\draw[draw=black,fill=orange,opacity=0.7,line width=0.48pt] (axis cs:89.3368,0) rectangle (axis cs:92.416,67);
\draw[draw=black,fill=orange,opacity=0.7,line width=0.48pt] (axis cs:92.416,0) rectangle (axis cs:95.4952,89);
\draw[draw=black,fill=orange,opacity=0.7,line width=0.48pt] (axis cs:95.4952,0) rectangle (axis cs:98.5744,57);
\draw[draw=black,fill=orange,opacity=0.7,line width=0.48pt] (axis cs:98.5744,0) rectangle (axis cs:101.6536,43);
\draw[draw=black,fill=orange,opacity=0.7,line width=0.48pt] (axis cs:101.6536,0) rectangle (axis cs:104.7328,22);
\draw[draw=black,fill=orange,opacity=0.7,line width=0.48pt] (axis cs:104.7328,0) rectangle (axis cs:107.812,20);
\draw[draw=black,fill=orange,opacity=0.7,line width=0.48pt] (axis cs:107.812,0) rectangle (axis cs:110.8912,26);
\draw[draw=black,fill=orange,opacity=0.7,line width=0.48pt] (axis cs:110.8912,0) rectangle (axis cs:113.9704,18);
\draw[draw=black,fill=orange,opacity=0.7,line width=0.48pt] (axis cs:113.9704,0) rectangle (axis cs:117.0496,4);
\draw[draw=black,fill=orange,opacity=0.7,line width=0.48pt] (axis cs:117.0496,0) rectangle (axis cs:120.1288,0);
\draw[draw=black,fill=orange,opacity=0.7,line width=0.48pt] (axis cs:120.1288,0) rectangle (axis cs:123.208,25);
\draw[draw=black,fill=orange,opacity=0.7,line width=0.48pt] (axis cs:123.208,0) rectangle (axis cs:126.2872,8);
\draw[draw=black,fill=orange,opacity=0.7,line width=0.48pt] (axis cs:126.2872,0) rectangle (axis cs:129.3664,1);
\draw[draw=black,fill=orange,opacity=0.7,line width=0.48pt] (axis cs:129.3664,0) rectangle (axis cs:132.4456,1);
\draw[draw=black,fill=orange,opacity=0.7,line width=0.48pt] (axis cs:132.4456,0) rectangle (axis cs:135.5248,3);
\draw[draw=black,fill=orange,opacity=0.7,line width=0.48pt] (axis cs:135.5248,0) rectangle (axis cs:138.604,9);
\draw[draw=black,fill=orange,opacity=0.7,line width=0.48pt] (axis cs:138.604,0) rectangle (axis cs:141.6832,6);
\draw[draw=black,fill=orange,opacity=0.7,line width=0.48pt] (axis cs:141.6832,0) rectangle (axis cs:144.7624,0);
\draw[draw=black,fill=orange,opacity=0.7,line width=0.48pt] (axis cs:144.7624,0) rectangle (axis cs:147.8416,0);
\draw[draw=black,fill=orange,opacity=0.7,line width=0.48pt] (axis cs:147.8416,0) rectangle (axis cs:150.9208,1);
\draw[draw=black,fill=orange,opacity=0.7,line width=0.48pt] (axis cs:150.9208,0) rectangle (axis cs:154,3);
\end{axis}

\end{tikzpicture}

%% file: figures/appendix/instances_description/max_trip_durations_histogram.tex
\begin{tikzpicture}

\definecolor{darkgray176}{RGB}{176,176,176}
\definecolor{mediumorchid}{RGB}{186,85,211}

\begin{axis}[
height=\InstancesHistogramsHeight,
minor xtick={},
minor ytick={},
scaled y ticks=true,
tick align=outside,
tick pos=left,
width=\InstancesHistogramsWidth,
x grid style={darkgray176},
xlabel={[min]},
xmin=-0.9, xmax=26,
xtick style={color=black},
xtick={-6,0,6,12,18,24,30},
xticklabel style={/pgf/number format/fixed, /pgf/number format/precision=0},
y grid style={darkgray176},
y tick label style={/pgf/number format/sci, /pgf/number format/sci zerofill=false},
ymin=0, ymax=47561.85,
ytick style={color=black},
ytick={0,10000,20000,30000,40000,50000},
yticklabel style={/pgf/number format/sci, /pgf/number format/precision=0},
yticklabels={0,1,2,3,4,5}
]
\draw[draw=black,fill=mediumorchid,opacity=0.7,line width=0.48pt] (axis cs:2.26666666666667,0) rectangle (axis cs:3.64,610);
\draw[draw=black,fill=mediumorchid,opacity=0.7,line width=0.48pt] (axis cs:3.64,0) rectangle (axis cs:5.01333333333333,8965);
\draw[draw=black,fill=mediumorchid,opacity=0.7,line width=0.48pt] (axis cs:5.01333333333333,0) rectangle (axis cs:6.38666666666667,33162);
\draw[draw=black,fill=mediumorchid,opacity=0.7,line width=0.48pt] (axis cs:6.38666666666667,0) rectangle (axis cs:7.76,45297);
\draw[draw=black,fill=mediumorchid,opacity=0.7,line width=0.48pt] (axis cs:7.76,0) rectangle (axis cs:9.13333333333333,39403);
\draw[draw=black,fill=mediumorchid,opacity=0.7,line width=0.48pt] (axis cs:9.13333333333333,0) rectangle (axis cs:10.5066666666667,27537);
\draw[draw=black,fill=mediumorchid,opacity=0.7,line width=0.48pt] (axis cs:10.5066666666667,0) rectangle (axis cs:11.88,13657);
\draw[draw=black,fill=mediumorchid,opacity=0.7,line width=0.48pt] (axis cs:11.88,0) rectangle (axis cs:13.2533333333333,7584);
\draw[draw=black,fill=mediumorchid,opacity=0.7,line width=0.48pt] (axis cs:13.2533333333333,0) rectangle (axis cs:14.6266666666667,4661);
\draw[draw=black,fill=mediumorchid,opacity=0.7,line width=0.48pt] (axis cs:14.6266666666667,0) rectangle (axis cs:16,3019);
\draw[draw=black,fill=mediumorchid,opacity=0.7,line width=0.48pt] (axis cs:16,0) rectangle (axis cs:17.3733333333333,1884);
\draw[draw=black,fill=mediumorchid,opacity=0.7,line width=0.48pt] (axis cs:17.3733333333333,0) rectangle (axis cs:18.7466666666667,1101);
\draw[draw=black,fill=mediumorchid,opacity=0.7,line width=0.48pt] (axis cs:18.7466666666667,0) rectangle (axis cs:20.12,596);
\draw[draw=black,fill=mediumorchid,opacity=0.7,line width=0.48pt] (axis cs:20.12,0) rectangle (axis cs:21.4933333333333,328);
\draw[draw=black,fill=mediumorchid,opacity=0.7,line width=0.48pt] (axis cs:21.4933333333333,0) rectangle (axis cs:22.8666666666667,215);
\draw[draw=black,fill=mediumorchid,opacity=0.7,line width=0.48pt] (axis cs:22.8666666666667,0) rectangle (axis cs:24.24,78);
\draw[draw=black,fill=mediumorchid,opacity=0.7,line width=0.48pt] (axis cs:24.24,0) rectangle (axis cs:25.6133333333333,18);
\draw[draw=black,fill=mediumorchid,opacity=0.7,line width=0.48pt] (axis cs:25.6133333333333,0) rectangle (axis cs:26.9866666666667,15);
\draw[draw=black,fill=mediumorchid,opacity=0.7,line width=0.48pt] (axis cs:26.9866666666667,0) rectangle (axis cs:28.36,8);
\draw[draw=black,fill=mediumorchid,opacity=0.7,line width=0.48pt] (axis cs:28.36,0) rectangle (axis cs:29.7333333333333,10);
\draw[draw=black,fill=mediumorchid,opacity=0.7,line width=0.48pt] (axis cs:29.7333333333333,0) rectangle (axis cs:31.1066666666667,11);
\draw[draw=black,fill=mediumorchid,opacity=0.7,line width=0.48pt] (axis cs:31.1066666666667,0) rectangle (axis cs:32.48,9);
\draw[draw=black,fill=mediumorchid,opacity=0.7,line width=0.48pt] (axis cs:32.48,0) rectangle (axis cs:33.8533333333333,3);
\draw[draw=black,fill=mediumorchid,opacity=0.7,line width=0.48pt] (axis cs:33.8533333333333,0) rectangle (axis cs:35.2266666666667,4);
\draw[draw=black,fill=mediumorchid,opacity=0.7,line width=0.48pt] (axis cs:35.2266666666667,0) rectangle (axis cs:36.6,7);
\draw[draw=black,fill=mediumorchid,opacity=0.7,line width=0.48pt] (axis cs:36.6,0) rectangle (axis cs:37.9733333333333,5);
\draw[draw=black,fill=mediumorchid,opacity=0.7,line width=0.48pt] (axis cs:37.9733333333333,0) rectangle (axis cs:39.3466666666667,0);
\draw[draw=black,fill=mediumorchid,opacity=0.7,line width=0.48pt] (axis cs:39.3466666666667,0) rectangle (axis cs:40.72,2);
\draw[draw=black,fill=mediumorchid,opacity=0.7,line width=0.48pt] (axis cs:40.72,0) rectangle (axis cs:42.0933333333333,1);
\draw[draw=black,fill=mediumorchid,opacity=0.7,line width=0.48pt] (axis cs:42.0933333333333,0) rectangle (axis cs:43.4666666666667,1);
\draw[draw=black,fill=mediumorchid,opacity=0.7,line width=0.48pt] (axis cs:43.4666666666667,0) rectangle (axis cs:44.84,5);
\draw[draw=black,fill=mediumorchid,opacity=0.7,line width=0.48pt] (axis cs:44.84,0) rectangle (axis cs:46.2133333333333,3);
\draw[draw=black,fill=mediumorchid,opacity=0.7,line width=0.48pt] (axis cs:46.2133333333333,0) rectangle (axis cs:47.5866666666667,2);
\draw[draw=black,fill=mediumorchid,opacity=0.7,line width=0.48pt] (axis cs:47.5866666666667,0) rectangle (axis cs:48.96,4);
\draw[draw=black,fill=mediumorchid,opacity=0.7,line width=0.48pt] (axis cs:48.96,0) rectangle (axis cs:50.3333333333333,4);
\draw[draw=black,fill=mediumorchid,opacity=0.7,line width=0.48pt] (axis cs:50.3333333333333,0) rectangle (axis cs:51.7066666666667,1);
\draw[draw=black,fill=mediumorchid,opacity=0.7,line width=0.48pt] (axis cs:51.7066666666667,0) rectangle (axis cs:53.08,1);
\draw[draw=black,fill=mediumorchid,opacity=0.7,line width=0.48pt] (axis cs:53.08,0) rectangle (axis cs:54.4533333333333,1);
\draw[draw=black,fill=mediumorchid,opacity=0.7,line width=0.48pt] (axis cs:54.4533333333333,0) rectangle (axis cs:55.8266666666667,3);
\draw[draw=black,fill=mediumorchid,opacity=0.7,line width=0.48pt] (axis cs:55.8266666666667,0) rectangle (axis cs:57.2,5);
\draw[draw=black,fill=mediumorchid,opacity=0.7,line width=0.48pt] (axis cs:57.2,0) rectangle (axis cs:58.5733333333333,2);
\draw[draw=black,fill=mediumorchid,opacity=0.7,line width=0.48pt] (axis cs:58.5733333333333,0) rectangle (axis cs:59.9466666666667,1);
\draw[draw=black,fill=mediumorchid,opacity=0.7,line width=0.48pt] (axis cs:59.9466666666667,0) rectangle (axis cs:61.32,3);
\draw[draw=black,fill=mediumorchid,opacity=0.7,line width=0.48pt] (axis cs:61.32,0) rectangle (axis cs:62.6933333333333,2);
\draw[draw=black,fill=mediumorchid,opacity=0.7,line width=0.48pt] (axis cs:62.6933333333333,0) rectangle (axis cs:64.0666666666667,0);
\draw[draw=black,fill=mediumorchid,opacity=0.7,line width=0.48pt] (axis cs:64.0666666666667,0) rectangle (axis cs:65.44,5);
\draw[draw=black,fill=mediumorchid,opacity=0.7,line width=0.48pt] (axis cs:65.44,0) rectangle (axis cs:66.8133333333333,3);
\draw[draw=black,fill=mediumorchid,opacity=0.7,line width=0.48pt] (axis cs:66.8133333333333,0) rectangle (axis cs:68.1866666666667,1);
\draw[draw=black,fill=mediumorchid,opacity=0.7,line width=0.48pt] (axis cs:68.1866666666667,0) rectangle (axis cs:69.56,2);
\draw[draw=black,fill=mediumorchid,opacity=0.7,line width=0.48pt] (axis cs:69.56,0) rectangle (axis cs:70.9333333333333,1);
\end{axis}

\end{tikzpicture}

%% file: figures/appendix/instances_description/max_staggering_histogram.tex
\begin{tikzpicture}

\definecolor{darkgray176}{RGB}{176,176,176}
\definecolor{skyblue}{RGB}{135,206,235}

\begin{axis}[
height=\InstancesHistogramsHeight,
minor xtick={},
minor ytick={},
scaled y ticks=true,
tick align=outside,
tick pos=left,
width=\InstancesHistogramsWidth,
x grid style={darkgray176},
xlabel={[min]},
xmin=-0.01, xmax=3.53744166666667,
xtick style={color=black},
xtick={-1,0,1,2,3,4},
xticklabel style={/pgf/number format/fixed, /pgf/number format/precision=0},
y grid style={darkgray176},
y tick label style={/pgf/number format/sci, /pgf/number format/sci zerofill=false},
ymin=0, ymax=41522.25,
ytick style={color=black},
ytick={0,10000,20000,30000,40000,50000},
yticklabel style={/pgf/number format/sci, /pgf/number format/precision=0},
yticklabels={0,1,2,3,4,5}
]
\draw[draw=black,fill=skyblue,opacity=0.7,line width=0.48pt] (axis cs:0.332666666666667,0) rectangle (axis cs:0.485275000000001,357);
\draw[draw=black,fill=skyblue,opacity=0.7,line width=0.48pt] (axis cs:0.485275000000001,0) rectangle (axis cs:0.637883333333334,4332);
\draw[draw=black,fill=skyblue,opacity=0.7,line width=0.48pt] (axis cs:0.637883333333334,0) rectangle (axis cs:0.790491666666667,18734);
\draw[draw=black,fill=skyblue,opacity=0.7,line width=0.48pt] (axis cs:0.790491666666667,0) rectangle (axis cs:0.9431,34926);
\draw[draw=black,fill=skyblue,opacity=0.7,line width=0.48pt] (axis cs:0.9431,0) rectangle (axis cs:1.09570833333333,39545);
\draw[draw=black,fill=skyblue,opacity=0.7,line width=0.48pt] (axis cs:1.09570833333333,0) rectangle (axis cs:1.24831666666667,35739);
\draw[draw=black,fill=skyblue,opacity=0.7,line width=0.48pt] (axis cs:1.24831666666667,0) rectangle (axis cs:1.400925,22405);
\draw[draw=black,fill=skyblue,opacity=0.7,line width=0.48pt] (axis cs:1.400925,0) rectangle (axis cs:1.55353333333333,11081);
\draw[draw=black,fill=skyblue,opacity=0.7,line width=0.48pt] (axis cs:1.55353333333333,0) rectangle (axis cs:1.70614166666667,6224);
\draw[draw=black,fill=skyblue,opacity=0.7,line width=0.48pt] (axis cs:1.70614166666667,0) rectangle (axis cs:1.85875,5097);
\draw[draw=black,fill=skyblue,opacity=0.7,line width=0.48pt] (axis cs:1.85875,0) rectangle (axis cs:2.01135833333333,3565);
\draw[draw=black,fill=skyblue,opacity=0.7,line width=0.48pt] (axis cs:2.01135833333333,0) rectangle (axis cs:2.16396666666667,2120);
\draw[draw=black,fill=skyblue,opacity=0.7,line width=0.48pt] (axis cs:2.16396666666667,0) rectangle (axis cs:2.316575,1590);
\draw[draw=black,fill=skyblue,opacity=0.7,line width=0.48pt] (axis cs:2.316575,0) rectangle (axis cs:2.46918333333333,931);
\draw[draw=black,fill=skyblue,opacity=0.7,line width=0.48pt] (axis cs:2.46918333333333,0) rectangle (axis cs:2.62179166666667,763);
\draw[draw=black,fill=skyblue,opacity=0.7,line width=0.48pt] (axis cs:2.62179166666667,0) rectangle (axis cs:2.7744,368);
\draw[draw=black,fill=skyblue,opacity=0.7,line width=0.48pt] (axis cs:2.7744,0) rectangle (axis cs:2.92700833333333,303);
\draw[draw=black,fill=skyblue,opacity=0.7,line width=0.48pt] (axis cs:2.92700833333333,0) rectangle (axis cs:3.07961666666667,136);
\draw[draw=black,fill=skyblue,opacity=0.7,line width=0.48pt] (axis cs:3.07961666666667,0) rectangle (axis cs:3.232225,21);
\draw[draw=black,fill=skyblue,opacity=0.7,line width=0.48pt] (axis cs:3.232225,0) rectangle (axis cs:3.38483333333333,3);
\end{axis}

\end{tikzpicture}

%% file: figures/appendix/parameter_search/heatmap-grid-relative-improvement_dp_5_tex/heatmap-grid-relative-improvement_dp_5.tex
\begin{tikzpicture}

\definecolor{darkgray176}{RGB}{176,176,176}

\begin{groupplot}[group style={group size=2 by 3, horizontal sep=\HeatmapHorSep, vertical sep=\HeatmapVerSep}, width=\HeatmapWidth, height=\HeatmapHeight]
\nextgroupplot[
scaled x ticks=manual:{}{\pgfmathparse{#1}},
tick align=outside,
tick pos=left,
title style={yshift=-5pt}, title={z = 1},
x grid style={darkgray176},
xmin=0, xmax=5,
xtick style={color=black},
xticklabels={},
y dir=reverse,
ylabel={Day},
ymin=0, ymax=6, 
ytick={0.5,1.5,2.5,3.5,4.5,5.5},
yticklabels={27,28,29,30,31,Avg}
]
\addplot graphics [includegraphics cmd=\pgfimage,xmin=0, xmax=5, ymin=6, ymax=0] {\heatmappath{5}heatmap-grid-relative-improvement_dp_5-000.png};

\nextgroupplot[
scaled x ticks=manual:{}{\pgfmathparse{#1}},
tick align=outside,
tick pos=left,
title style={yshift=-5pt}, title={z = 2},
x grid style={darkgray176},
xmin=0, xmax=5,
xtick style={color=black},
xticklabels={},
y dir=reverse,
ymin=0, ymax=6, 
 yticklabels={},
]
\addplot graphics [includegraphics cmd=\pgfimage,xmin=0, xmax=5, ymin=6, ymax=0] {\heatmappath{5}heatmap-grid-relative-improvement_dp_5-001.png};

\nextgroupplot[
scaled x ticks=manual:{}{\pgfmathparse{#1}},
tick align=outside,
tick pos=left,
title style={yshift=-5pt}, title={z = 3},
x grid style={darkgray176},
xlabel={},
xmin=0, xmax=5,
xtick style={color=black},
xticklabels={},
y dir=reverse,
ylabel={Day},
ymin=0, ymax=6, 
ytick={0.5,1.5,2.5,3.5,4.5,5.5},
yticklabels={27,28,29,30,31,Avg}
]
\addplot graphics [includegraphics cmd=\pgfimage,xmin=0, xmax=5, ymin=6, ymax=0] {\heatmappath{5}heatmap-grid-relative-improvement_dp_5-002.png};
\nextgroupplot[
scaled x ticks=manual:{}{\pgfmathparse{#1}},
tick align=outside,
tick pos=left,
title style={yshift=-5pt}, title={z = 4},
x grid style={darkgray176},
xlabel={},
xmin=0, xmax=5,
xtick style={color=black},
xticklabels={},
y dir=reverse,
ymin=0, ymax=6, 
 yticklabels={},
]
\addplot graphics [includegraphics cmd=\pgfimage,xmin=0, xmax=5, ymin=6, ymax=0] {\heatmappath{5}heatmap-grid-relative-improvement_dp_5-003.png};
\nextgroupplot[
tick align=outside,
tick pos=left,
title style={yshift=-5pt}, title={z = 5},
x grid style={darkgray176},
xlabel={x [\%]},
xmin=0, xmax=5,
xtick={0.5,1.5,2.5,3.5,4.5},
xticklabels={10,20,30,40,50},
xtick style={color=black},
y dir=reverse,
ylabel={Day},
ymin=0, ymax=6, 
ytick={0.5,1.5,2.5,3.5,4.5,5.5},
yticklabels={27,28,29,30,31,Avg}
]
\addplot graphics [includegraphics cmd=\pgfimage,xmin=0, xmax=5, ymin=6, ymax=0] {\heatmappath{5}heatmap-grid-relative-improvement_dp_5-004.png};

\nextgroupplot[
tick align=outside,
tick pos=left,
title style={yshift=-5pt}, title={z = 6},
x grid style={darkgray176},
xlabel={x [\%]},
xmin=0, xmax=5,
xtick={0.5,1.5,2.5,3.5,4.5},
xticklabels={10,20,30,40,50},
xtick style={color=black},
y dir=reverse,
ymin=0, ymax=6, 
yticklabels={},
]
\addplot graphics [includegraphics cmd=\pgfimage,xmin=0, xmax=5, ymin=6, ymax=0] {\heatmappath{5}heatmap-grid-relative-improvement_dp_5-005.png};
\end{groupplot}

\end{tikzpicture}

%% file: figures/appendix/parameter_search/heatmap-grid-relative-improvement_dp_10_tex/heatmap-grid-relative-improvement_dp_10.tex
\begin{tikzpicture}

\definecolor{darkgray176}{RGB}{176,176,176}

\begin{groupplot}[group style={group size=2 by 3, horizontal sep=\HeatmapHorSep, vertical sep=\HeatmapVerSep}, width=\HeatmapWidth, height=\HeatmapHeight]
\nextgroupplot[
scaled x ticks=manual:{}{\pgfmathparse{#1}},
tick align=outside,
tick pos=left,
title style={yshift=-5pt}, title={z = 1},
x grid style={darkgray176},
xmin=0, xmax=5,
xtick style={color=black},
xticklabels={},
y dir=reverse,
ymin=0, ymax=6, 
 yticklabels={},
]
\addplot graphics [includegraphics cmd=\pgfimage,xmin=0, xmax=5, ymin=6, ymax=0] {\heatmappath{10}heatmap-grid-relative-improvement_dp_10-000.png};

\nextgroupplot[
scaled x ticks=manual:{}{\pgfmathparse{#1}},
tick align=outside,
tick pos=left,
title style={yshift=-5pt}, title={z = 2},
x grid style={darkgray176},
xmin=0, xmax=5,
xtick style={color=black},
xticklabels={},
y dir=reverse,
ymin=0, ymax=6, 
 yticklabels={},
]
\addplot graphics [includegraphics cmd=\pgfimage,xmin=0, xmax=5, ymin=6, ymax=0] {\heatmappath{10}heatmap-grid-relative-improvement_dp_10-001.png};

\nextgroupplot[
scaled x ticks=manual:{}{\pgfmathparse{#1}},
tick align=outside,
tick pos=left,
title style={yshift=-5pt}, title={z = 3},
x grid style={darkgray176},
xlabel={},
xmin=0, xmax=5,
xtick style={color=black},
xticklabels={},
y dir=reverse,
ymin=0, ymax=6, 
 yticklabels={},
]
\addplot graphics [includegraphics cmd=\pgfimage,xmin=0, xmax=5, ymin=6, ymax=0] {\heatmappath{10}heatmap-grid-relative-improvement_dp_10-002.png};

\nextgroupplot[
scaled x ticks=manual:{}{\pgfmathparse{#1}},
tick align=outside,
tick pos=left,
title style={yshift=-5pt}, title={z = 4},
x grid style={darkgray176},
xlabel={},
xmin=0, xmax=5,
xtick style={color=black},
xticklabels={},
y dir=reverse,
ymin=0, ymax=6, 
 yticklabels={},
]
\addplot graphics [includegraphics cmd=\pgfimage,xmin=0, xmax=5, ymin=6, ymax=0] {\heatmappath{10}heatmap-grid-relative-improvement_dp_10-003.png};

\nextgroupplot[
tick align=outside,
tick pos=left,
title style={yshift=-5pt}, title={z = 5},
x grid style={darkgray176},
xlabel={x [\%]},
xmin=0, xmax=5,
xtick={0.5,1.5,2.5,3.5,4.5},
xticklabels={10,20,30,40,50},
xtick style={color=black},
y dir=reverse,
ymin=0, ymax=6, 
yticklabels={},
]
\addplot graphics [includegraphics cmd=\pgfimage,xmin=0, xmax=5, ymin=6, ymax=0] {\heatmappath{10}heatmap-grid-relative-improvement_dp_10-004.png};

\nextgroupplot[
tick align=outside,
tick pos=left,
title style={yshift=-5pt}, title={z = 6},
x grid style={darkgray176},
xlabel={x [\%]},
xmin=0, xmax=5,
xtick={0.5,1.5,2.5,3.5,4.5},
xticklabels={10,20,30,40,50},
xtick style={color=black},
y dir=reverse,
ymin=0, ymax=6, 
yticklabels={},
]
\addplot graphics [includegraphics cmd=\pgfimage,xmin=0, xmax=5, ymin=6, ymax=0] {\heatmappath{10}heatmap-grid-relative-improvement_dp_10-005.png};

\end{groupplot}

\end{tikzpicture}

%% file: figures/appendix/parameter_search/heatmap-grid-relative-improvement_dp_20_tex/heatmap-grid-relative-improvement_dp_20.tex
\begin{tikzpicture}

\definecolor{darkgray176}{RGB}{176,176,176}

\begin{groupplot}[group style={group size=2 by 3, horizontal sep=\HeatmapHorSep, vertical sep=\HeatmapVerSep}, width=\HeatmapWidth, height=\HeatmapHeight]
\nextgroupplot[
scaled x ticks=manual:{}{\pgfmathparse{#1}},
tick align=outside,
tick pos=left,
title style={yshift=-5pt}, title={z = 1},
x grid style={darkgray176},
xmin=0, xmax=5,
xtick style={color=black},
xticklabels={},
y dir=reverse,
ymin=0, ymax=6, 
 yticklabels={},
]
\addplot graphics [includegraphics cmd=\pgfimage,xmin=0, xmax=5, ymin=6, ymax=0] {\heatmappath{20}heatmap-grid-relative-improvement_dp_20-000.png};

\nextgroupplot[
colorbar right,
colorbar style={at={(\XCbar,\YCbar)}, anchor=west,height=\HeightCbar,ylabel={[\%]}},
colormap={mymap}{[1pt]
  rgb(0pt)=(0.2298057,0.298717966,0.753683153);
  rgb(1pt)=(0.26623388,0.353094838,0.801466763);
  rgb(2pt)=(0.30386891,0.406535296,0.84495867);
  rgb(3pt)=(0.342804478,0.458757618,0.883725899);
  rgb(4pt)=(0.38301334,0.50941904,0.917387822);
  rgb(5pt)=(0.424369608,0.558148092,0.945619588);
  rgb(6pt)=(0.46666708,0.604562568,0.968154911);
  rgb(7pt)=(0.509635204,0.648280772,0.98478814);
  rgb(8pt)=(0.552953156,0.688929332,0.995375608);
  rgb(9pt)=(0.596262162,0.726149107,0.999836203);
  rgb(10pt)=(0.639176211,0.759599947,0.998151185);
  rgb(11pt)=(0.681291281,0.788964712,0.990363227);
  rgb(12pt)=(0.722193294,0.813952739,0.976574709);
  rgb(13pt)=(0.761464949,0.834302879,0.956945269);
  rgb(14pt)=(0.798691636,0.849786142,0.931688648);
  rgb(15pt)=(0.833466556,0.860207984,0.901068838);
  rgb(16pt)=(0.865395197,0.86541021,0.865395561);
  rgb(17pt)=(0.897787179,0.848937047,0.820880546);
  rgb(18pt)=(0.924127593,0.827384882,0.774508472);
  rgb(19pt)=(0.944468518,0.800927443,0.726736146);
  rgb(20pt)=(0.958852946,0.769767752,0.678007945);
  rgb(21pt)=(0.96732803,0.734132809,0.628751763);
  rgb(22pt)=(0.969954137,0.694266682,0.579375448);
  rgb(23pt)=(0.966811177,0.650421156,0.530263762);
  rgb(24pt)=(0.958003065,0.602842431,0.481775914);
  rgb(25pt)=(0.943660866,0.551750968,0.434243684);
  rgb(26pt)=(0.923944917,0.49730856,0.387970225);
  rgb(27pt)=(0.89904617,0.439559467,0.343229596);
  rgb(28pt)=(0.869186849,0.378313092,0.300267182);
  rgb(29pt)=(0.834620542,0.312874446,0.259301199);
  rgb(30pt)=(0.795631745,0.24128379,0.220525627);
  rgb(31pt)=(0.752534934,0.157246067,0.184115123);
  rgb(32pt)=(0.705673158,0.01555616,0.150232812)
},
point meta max=23.0052384545524,
point meta min=15.909170560695,
scaled x ticks=manual:{}{\pgfmathparse{#1}},
tick align=outside,
tick pos=left,
title style={yshift=-5pt}, title={z = 2},
x grid style={darkgray176},
xmin=0, xmax=5,
xtick style={color=black},
xticklabels={},
y dir=reverse,
ymin=0, ymax=6, 
 yticklabels={},
]
\addplot graphics [includegraphics cmd=\pgfimage,xmin=0, xmax=5, ymin=6, ymax=0] {\heatmappath{20}heatmap-grid-relative-improvement_dp_20-001.png};

\nextgroupplot[
scaled x ticks=manual:{}{\pgfmathparse{#1}},
tick align=outside,
tick pos=left,
title style={yshift=-5pt}, title={z = 3},
x grid style={darkgray176},
xlabel={},
xmin=0, xmax=5,
xtick style={color=black},
xticklabels={},
y dir=reverse,
ymin=0, ymax=6, 
 yticklabels={},
]
\addplot graphics [includegraphics cmd=\pgfimage,xmin=0, xmax=5, ymin=6, ymax=0] {\heatmappath{20}heatmap-grid-relative-improvement_dp_20-002.png};

\nextgroupplot[
scaled x ticks=manual:{}{\pgfmathparse{#1}},
tick align=outside,
tick pos=left,
title style={yshift=-5pt}, title={z = 4},
x grid style={darkgray176},
xlabel={},
xmin=0, xmax=5,
xtick style={color=black},
xticklabels={},
y dir=reverse,
ymin=0, ymax=6, 
 yticklabels={},
]
\addplot graphics [includegraphics cmd=\pgfimage,xmin=0, xmax=5, ymin=6, ymax=0] {\heatmappath{20}heatmap-grid-relative-improvement_dp_20-003.png};

\nextgroupplot[
tick align=outside,
tick pos=left,
title style={yshift=-5pt}, title={z = 5},
x grid style={darkgray176},
xlabel={x [\%]},
xmin=0, xmax=5,
xtick={0.5,1.5,2.5,3.5,4.5},
xticklabels={10,20,30,40,50},
xtick style={color=black},
y dir=reverse,
ymin=0, ymax=6, 
yticklabels={},
]
\addplot graphics [includegraphics cmd=\pgfimage,xmin=0, xmax=5, ymin=6, ymax=0] {\heatmappath{20}heatmap-grid-relative-improvement_dp_20-004.png};

\nextgroupplot[
tick align=outside,
tick pos=left,
title style={yshift=-5pt}, title={z = 6},
x grid style={darkgray176},
xlabel={x [\%]},
xmin=0, xmax=5,
xtick={0.5,1.5,2.5,3.5,4.5},
xticklabels={10,20,30,40,50},
xtick style={color=black},
y dir=reverse,
ymin=0, ymax=6, 
yticklabels={},
]
\addplot graphics [includegraphics cmd=\pgfimage,xmin=0, xmax=5, ymin=6, ymax=0] {\heatmappath{20}heatmap-grid-relative-improvement_dp_20-005.png};

\end{groupplot}

\end{tikzpicture}